\RequirePackage{amsfonts}
\RequirePackage{graphicx}
\RequirePackage{mathrsfs}
\RequirePackage{amsmath}
\RequirePackage{float}
\RequirePackage{rotating,color}
\RequirePackage{bm,amssymb,amsmath,mathrsfs,amsfonts}

\def\boxit#1{\vbox{\hrule\hbox{\vrule\kern6pt  \vbox{\kern6pt#1\kern6pt}\kern6pt\vrule}\hrule}}

\newcommand{\be}{\begin{equation}}
\newcommand{\ee}{\end{equation}}
\newcommand{\beaa}{\begin{eqnarray*}}
\newcommand{\eeaa}{\end{eqnarray*}}
\newcommand{\bea}{\begin{eqnarray}}
\newcommand{\eea}{\end{eqnarray}}
\newcommand{\lbl}{\label}

\newcommand{\ml}{\mathcal}

\newcommand{\bd}{\bold}

\def\diag{\mathrm {diag}}

\def\tr{\mathrm {tr}}

\def\R{{\bf R}}
\def\I{{\bf I}}

\def\D{{\bf D}}

\def\tr{\mathrm {tr}}

\def\bms{{\bm\Sigma}}

\newtheorem{theorem}{ \noindent T{\footnotesize HEOREM}}
\newtheorem{prop}{ \noindent P{\footnotesize ROPOSITION}}
\newtheorem{lemma}{ \noindent L{\footnotesize EMMA}}
\newtheorem{coro}{ \noindent C{\footnotesize OROLLARY}}
\newtheorem{remark}{ \noindent R{\footnotesize EMARK}}

\documentclass[aos,preprint]{imsart}

\RequirePackage[numbers]{natbib}
\RequirePackage[colorlinks,citecolor=blue,urlcolor=blue]{hyperref}

\startlocaldefs
\endlocaldefs

\begin{document}

\begin{frontmatter}

\title{Max-Sum tests for cross-sectional dependence of high-dimensional panel data}
\runtitle{Max-sum test for cross-sectional dependence}


\begin{aug}
\author{\fnms{Long} \snm{Feng}\thanksref{m1}\ead[label=e1]{fengl100@nenu.edu.cn}},
\author{\fnms{Tiefeng} \snm{Jiang}\thanksref{m2}\ead[label=e2]{jiang040@umn.edu}},
\author{\fnms{Binghui} \snm{Liu}\thanksref{m1}\ead[label=e3]{liubh100@nenu.edu.cn}}
\and
\author{\fnms{Wei} \snm{Xiong}\thanksref{m3}\ead[label=e4]{xiongwei@uibe.edu.cn}}
\affiliation{Northeast Normal University\thanksmark{m1}, University of Minnesota\thanksmark{m2} and University of International Business and Economics\thanksmark{m3} }
\runauthor{L. Feng et al.}

\address{Address of Long Feng and Binghui Liu\\
Key Laboratory of Applied Statistics of MOE $\&$ School of Mathematics and Statistics\\
Northeast Normal University\\
\phantom{E-mail:\ }\printead*{e1,e3}}

\address{Address of Tiefeng Jiang\\
School of Statistics\\
313 Ford Hall\\
224 Church Street SE\\
Minneapolis, MN55455\\
\printead{e2}}

\address{Address of Wei Xiong\\
School of Statistics\\
University of International Business and Economics\\
\printead{e4}}
\end{aug}

\begin{abstract}
We consider a testing problem for cross-sectional dependence  for high-dimensional panel data,
where the number of cross-sectional units is potentially much larger than the number of observations. The cross-sectional dependence is described through a linear regression model. We study three tests named  the sum test, the max test and the max-sum test, where the latter two are new. The sum test is initially proposed by Breusch and Pagan (1980). We design the max and sum tests for sparse and non-sparse residuals in the linear regressions, respectively.
And the max-sum test is devised to compromise both situations on the residuals.
Indeed, our simulation shows that the max-sum test outperforms the previous two tests. This makes the max-sum test very useful in practice where sparsity or not for a set of data  is usually vague. Towards the theoretical analysis of the three tests, we have settled two conjectures regarding the sum of squares of sample correlation coefficients asked by Pesaran (2004 and 2008). In addition, we establish the asymptotic theory for maxima of sample correlations coefficients appeared in the linear regression model for panel data, which is also the first successful attempt to our knowledge.   To study the max-sum test, we create a novel method to show asymptotic independence between maxima and sums of dependent random variables. We expect  the method itself is useful for other problems of this nature. Finally, an extensive simulation study as well as a case study are carried out. They demonstrate advantages of
our proposed methods in terms of both empirical powers and robustness for residuals regardless of sparsity or not.
%
%
%
%
%
\end{abstract}

\begin{keyword}
\kwd{high-dimensional data}
\kwd{panel data models}
\kwd{hypothesis tests}
\kwd{cross-sectional dependence}
\kwd{asymptotic normality}
\kwd{extreme-value distribution}
\kwd{asymptotic independence}
\kwd{max-sum test}.
\end{keyword}

\end{frontmatter}

\tableofcontents

\section{Introduction}\lbl{INtro}
In this paper we will study the cross-sectional dependence for the following linear regression   model for panel data
\begin{align}\lbl{haoshichengshuang}
y_{it}=x_{it}'\beta_i+\epsilon_{it}
\end{align}
for $i=1,\cdots,N$ and $t=1,\cdots,T$,
where $i$ represents households, individuals, firms, etc., and $t$ represents time. In the literature of panel data, the index $i$ stands for {\it sections}. For each section $i$, the corresponding model is a standard multiple linear regression model, where
$y_{it}\in \mathbb{R}$ is the dependent variable and $x_{it}\in \mathbb{R}^p$ is the  regressor
with slope parameter $\beta_i\in \mathbb{R}^p$. The first coordinate of $x_{it}$  is one if there is an intercept in the linear regression model \eqref{haoshichengshuang}. The value of  $\beta_i$  may vary across $i$. In \eqref{haoshichengshuang}, we assume  $\{\epsilon_{it};\, 1\leq t\leq T\}$ are independent and identically distributed (i.i.d.) for each section $i$. However, across sections the random errors may be dependent, that is, $\{\epsilon_{it};\, 1\leq i\leq N\}$ may be dependent for some $t.$ Such dependence is referred to as cross-sectional dependence. The objective of this paper is to test if there exists cross-sectional dependence by using a few of new methods. Before stating our results, we will introduce some background next.


In statistics and econometrics, panel data or longitudinal data are multi-dimensional data involving measurements
over time, which contain observations of various phenomena over multiple time periods for the same
unit, for instance, a household or a firm.  In the study of panel data models, the cross-sectional dependence is an important concept,
described as the interaction between cross-sectional units, which could arise from the behavioral
interaction between units.

Stephan \cite{Stephan1934} argues that ``in dealing with social
data, we know that by virtue of their very social character, persons, groups and their
characteristics are interrelated and not independent. " However, to make theoretical study easier, experts assume cross-sectional independence in various model setups \cite{HsiaoP2012, Pesaran2004}. If data across individuals are dependent, inferences under the assumption of cross-sectional independence would be inaccurate and misleading; see \cite{HsiaoP2012, pesaran2015testing} and the literature therein. To this end, testing the existence of cross-sectional
dependence is an important task, which has attracted more attention in recent years, see, for instance, \cite{Chudik2013, moscone2009a,  pesaran2015testing, Pesaran2015,ST12}.


 Perhaps the most widely known test for cross-sectional independence is the Lagrange Multiplier (LM) statistic proposed by Breusch and Pagan \cite{breusch1980the} in 1980 (Google records 5353 citations currently). Their  test statistic is  the sum of squares of sample correlation coefficients between the residuals from the ordinary least square
(OLS). Precisely, for each $i$, let $\hat{\beta}_i$ be the standard estimator of $\beta_i$ in the linear regression for observations  $\{(y_{it},x_{it});\, t=1,\cdots, T\}$ and  the quantity $\hat{\bd{\epsilon}}_{it}=y_{it}-x_{it}'\hat{\beta}_i$ denotes the residual.  For each $i,j=1,\cdots,N$, define the sample correlation $\hat{\rho}_{ij}$  by
\begin{align}
\hat{\rho}_{ij}=\frac{\sum_{t=1}^T \hat{\bd{\epsilon}}_{it}\hat{\bd{\epsilon}}_{jt}}{\sqrt{\sum_{t=1}^T \hat{\bd{\epsilon}}_{it}^2\sum_{t=1}^T \hat{\bd{\epsilon}}_{jt}^2}}.\lbl{liberation_free}
\end{align}
Breusch and Pagan \cite{breusch1980the}  propose the Lagrange multiplier test statistic defined by
\begin{align}\lbl{hide_she}
S_N:=\sum_{1\le i<j\le N}T\hat{\rho}_{ij}^2.
\end{align}
 To get the rejection region, we need to figure out the limiting distribution of $S_N$ as $N$ goes to infinity.
Under the null hypothesis that there is no cross-sectional dependence, that is,
$\{\epsilon_{it};\, 1\leq i\leq N, 1\leq t \leq T\}$ from \eqref{haoshichengshuang} are independent, the asymptotic distribution of $S_N$ is understood
when the cross-sectional dimension $N$ is fixed and the time dimension $T$ goes to infinity.
In fact, assuming that $\epsilon_{it}$'s are normally distributed, Breusch and Pagan \cite{breusch1980the} show that, for fixed $N$,
\begin{align}\lbl{dashfune}
S_N\to \chi^2(d)
\end{align}
in distribution as $T\to\infty$, where $d=N(N-1)/2.$ If $N$ is relatively large, the above chi-square approximation is not accurate \citep{Pesaran2004}. A natural amendment is approximating $\chi^2(d)$ by the standard normal distribution: $(\chi^2(d)-d)/\sqrt{2d}$ goes to the standard normal distribution as $N$ goes to infinity. However,
as both $N$ and $T$ are very large, taking limit by sending $T \to\infty$ followed by sending $N\to\infty$ is not legitimate mathematically, and the approximation may  not be accurate statistically (our Remark \ref{Re_5} shows such an example). For this consideration Pesaran \cite{Pesaran2004} and Pesaran {\it et al.} \cite{pesaran2008a} provide two versions of normalization of $S_N$ and conjecture that
%
both versions satisfy the central limit theorem (CLT); some of the insights why the CLTs hold  can be seen, for example, from \cite{ST12} and \cite{Pesaran2015}. In this paper we prove the two conjectures in Theorems \ref{diligent} and \ref{Fall_nice}. This enables us to carry out the test for cross-sectional dependence through  $S_N$ in \eqref{hide_she}. In the future, when $S_N$ is used to be a test statistic, we call it the sum test.

On the other hand, when data are sparse, experts in recent years  realize that a better test than sum statistics is the  maximum of sample correlation coefficients. This is confirmed in, for example,  \cite{Cai2014}; see also \cite{caiLiu2011Adaptive}, \cite{cai2013two-sample} and \cite{caiZhang2016}. With this philosophy in mind, to test the cross-sectional dependence when the residuals $\hat{\bd{\epsilon}}_{it}$ are sparse,  we propose statistic
\begin{align}\lbl{kkkl}
L_N:=\max_{1\le i<j\le N} |\hat{\rho}_{ij}|
\end{align}
 where $\hat{\rho}_{ij}$ is defined as in \eqref{liberation_free}. Later, when $L_N$ is used to be a test statistic, we refer it to as the max test. Its limiting distribution is obtained as both $N$ and $T$ go to infinity under various moment conditions (Theorems \ref{linear_case}, \ref{exponential_case} and \ref{super_exponential_case}). The corresponding rejection region based on the test statistic $L_N$ is given after Theorem \ref{super_exponential_case}.

 In practice it is hard to tell or differentiate if a set of data is sparse. We then combine the sum test $S_N$ and the max test $L_N$ to propose another test $C_N$,  which is the minimum of the $p$-values corresponding to the tests based on $S_N$ and  $L_N$. We prove in Theorem \ref{Asym_indept} that, under normalization,  $S_N$ and $L_N$ are asymptotically independent as both $N$ and $T$ go to infinity. Hence the limiting distribution of $C_N$ is identified. In further discussions, when $C_N$ is used to be a test statistic, we name it the max-sum test. From simulation we see this test, taking care of both sparsity and non-sparsity cases, is better than the sum test  and the max test. The tool of deriving asymptotic independence between the sum and the maximum of random variables is new to our knowledge. It seems a universal method to handle asymptotic independence between random variables of this nature.

To sum up, to test cross-sectional dependence for panel data models, in this paper we study three types of tests,
{\it i.e.}, the sum test, the max test and the max-sum test. To carry the test, we have solved two open problems about the CLTs for the sum of squares of residuals; the limiting distributions of the maxima of residuals are systematically studied; a new method of studying asymptotic independence is created to develop part of the above theory successfully.

%



%
%
%

%
%

\section{The proposed tests}
\label{Sec2}

\subsection{Problem description}

Review model \eqref{haoshichengshuang} that $y_{it}=x_{it}'\beta_i+\epsilon_{it}$
%
%
for $i=1,\cdots,N$ and $t=1,\cdots,T$, where $i$ indexes the cross-sectional units and $t$ indexes the observations. In this model,
$y_{it}\in \mathbb{R}$ is the dependent variable, and $x_{it}\in \mathbb{R}^p$ is the non-random, exogenous regressor
with slope parameter $\beta_i\in \mathbb{R}^p$ that are allowed to vary across $i$.
We assume  $\{\epsilon_{it};\, 1\leq t\leq T\}$ are i.i.d. real-valued random variables for each section $i$. However, across sections the random errors may be dependent, that is, $\{\epsilon_{it};\, 1\leq i\leq N\}$ may be dependent for some $t.$ Such dependence is called cross-sectional dependence.
Set
\begin{align} \label{shadongxi}
\bd{x}_i=(x_{i1},\cdots,x_{iT})',\ \bd{y}_i=(y_{i1},\cdots,y_{iT})', \ \bd{\epsilon}_{i}=(\epsilon_{i1},\cdots,\epsilon_{iT})'
\end{align}
for $i=1,2,\cdots, N.$ Then $\bd{x}_i$ is a $T\times p$ matrix; both $\bd{y}_i$ and $\bd{\epsilon}_{i}$ are $T$-dimensional vectors. Throughout the paper we assume that the $T$ entries of $\epsilon_i$ are i.i.d. with mean zero for each  $i$.
Recalling \eqref{haoshichengshuang}, the cross-sectional independence is the same as saying that 
\begin{align}\label{h1}
H_0: \bd{\epsilon}_{1}, \bd{\epsilon}_{2},\cdots, \bd{\epsilon}_{N}\ \mbox{are independent random vectors}.
\end{align}
In general, although sometimes we assume $\bd{\epsilon}_{1}$ has the normal distribution, we do not need the exact distribution of $\bd{\epsilon}_{1}$ but rather its moments.

%

\subsection{Test statistics}

First, we list some notations used in the rest of the paper. Reviewing \eqref{shadongxi}, for each $i=1,\cdots,N$, let
\begin{align}
\hat{\beta}_i=(\bd{x}_i' \bd{x}_i)^{-1}\bd{x}_i'\bd{y}_i\text{, }\bd{P}_i=\bd{I}_T-\bd{x}_i(\bd{x}_i'\bd{x}_i)^{-1}\bd{x}_i',\label{park_1}
\end{align}
where $\bd{I}_T$ is the $T\times T$ identity matrix and $\bd{P}_i$ is a $T\times T$ projection matrix with $\bd{P}_i^2=\bd{P}_i$ and the rank of $\bd{P}_i$ is $T-p$.
For each $i,j=1,\cdots,N$, let $\hat{\rho}_{ij}$ denote the sample correlation coefficient computed by the Ordinary Least Squares (OLS)
residuals $(\hat{\bd{\epsilon}}_{i1}, \cdots, \hat{\bd{\epsilon}}_{iT})^T$ and $(\hat{\bd{\epsilon}}_{j1}, \cdots, \hat{\bd{\epsilon}}_{jT})^T$ where  $\hat{\bd{\epsilon}}_{it}=y_{it}-x_{it}'\hat{\beta}_i$ for each $i$ and $t$. Under model \eqref{haoshichengshuang}, it is easy to see that
\beaa
(\hat{\bd{\epsilon}}_{i1}, \cdots, \hat{\bd{\epsilon}}_{iT})^T=\bd{P}_i\bd{\epsilon}_i
\eeaa
for each $i$. Thus, by \eqref{liberation_free},
\begin{align}
\hat{\rho}_{ij}=\frac{\sum_{t=1}^T \hat{\bd{\epsilon}}_{it}\hat{\bd{\epsilon}}_{jt}}{\sqrt{\sum_{t=1}^T \hat{\bd{\epsilon}}_{it}^2}\sqrt{\sum_{t=1}^T \hat{\bd{\epsilon}}_{jt}^2}}
=\frac{\bd{\epsilon}_i'\bd{P}_i\bd{P}_j\bd{\epsilon}_j}
{\|\bd{P}_i\bd{\epsilon}_i\|\cdot\|\bd{P}_j\bd{\epsilon}_j\|}\label{Asia_foods}.
\end{align}

In this paper, to test the null hypothesis \eqref{h1}, we will study three types of tests as follows:
\begin{align}
\textrm{sum: } S_N&=\sum_{1\le i<j\le N}T\hat{\rho}_{ij}^2,\lbl{sumtest}\\
\textrm{max: } L_N&=\max_{1\le i<j\le N} |\hat{\rho}_{ij}|,\lbl{Write}\\
\textrm{max-sum: }C_N&=\min\{P_{L_N},P_{S_N}\},\lbl{maxsumtest}
\end{align}
respectively, where
\begin{align}
P_{L_N}&=1-F(TL_N^2-4\log N+\log\log N), \nonumber\\
P_{S_N}&=1-\Phi\Big(\frac{S_N-\mu_N}{N}\Big),\nonumber\\
\mu_N&=\frac{T}{(T-p)^2}\sum_{1\le i<j \le N} \tr(\bd{P}_i\bd{P}_j).\lbl{taiyangtaila}
\end{align}
Here, $F(y)= \exp(-e^{-y/2}/\sqrt{8\pi})$ is the extreme-value distribution function of  type I,
also called the Gumble distribution in literature, and $\Phi(y)$ is the distribution function of
 $N(0,1)$.

 For the sum test in \eqref{sumtest}, we will establish that, under $H_0$ in \eqref{h1}, $(S_N-\mu_N)/N$
converges weakly to the standard normal distribution when both $N$ and $T$ go to infinity with a certain restriction (Theorem \ref{diligent}), hence a level-$\alpha$ test will
be performed through rejecting $H_0$ when $(S_N-\mu_N)/N$ is larger than the $1-\alpha$ quantile $z_{\alpha}= \Phi^{-1}(1-\alpha)$ of the standard normal distribution.

For the max test in \eqref{Write}, under $H_0$, we will establish that $TL_N^2-4\log N+\log\log N$ has an asymptotic extreme-value
distribution as both $N$ and $T$ go to infinity (Theorems \ref{linear_case}, \ref{exponential_case} and  \ref{super_exponential_case}). We do not impose normality assumptions but rather moment conditions.
Recall $F(y)$ is defined below \eqref{taiyangtaila}.  A level-$\alpha$ test will
then be performed by rejecting $H_0$ when $TL_N^2-4\log N+\log\log N$ is larger than the $1-\alpha$ quantile
$q_{\alpha}=-\log(8\pi)-2\log\log(1-\alpha)^{-1}$ of $F(y)$.

Furthermore, for the max-sum test in \eqref{maxsumtest}, its asymptotic distribution under $H_0$ is constructed
based on the asymptotic independence between  $(S_N-\mu_N)/N$ and $TL_N^2-4\log N+\log\log N$ as both $N$ and $T$ go to infinity (Theorem \ref{Asym_indept} and Corollary \ref{coro2}). So  a level-$\alpha$ test will be performed through rejecting $H_0$ when $C_N<1-\sqrt{1-\alpha}$.

\subsection{Contributions}

In this paper, for the panel data model \eqref{haoshichengshuang} we study the cross-sectional dependence.  The asymptotic distributions of three test statistics based on residuals are established. As application, three hypothesis tests are accomplished. A real data analysis by using our results is provided. We will now further elaborate below.

In the theoretical part, we have solved two open problems on the sum of squares of residuals conjectured by economists (\cite{Pesaran2004, pesaran2008a}; see also \cite{ Pesaran2015, ST12}). We have developed an extreme-value theory for the maximum of residuals. Further, a new method is developed to show the sum and the maximum are asymptotically independent. There are not many results in literature to show asymptotic independence between  sums of and maxima of random variables. Close references are \citep{BBQ98, Xu2016}. Our method,  being  different from earlier literature, provides a general and novel tool for showing asymptotic independence between sums of and maxima of random variables.

In application, we propose three tests on the cross-sectional dependence for high-dimensional panel data: the sum test, the max test and the max-sum test.
The max test is the first high-dimensional max test for cross-sectional
dependence in panel data models, which is good for sparse residuals while  existing test statistics of sum types tend to fail. The sum test is useful for non-sparse residuals, which is clearly demonstrated by simulation in, for example, \cite{Pesaran2004, Pesaran2015, pesaran2008a, ST12}. We are able to derive the limiting distribution of the sums in this paper.

Furthermore, the max-sum test is constructed based on the  asymptotic independence
between the max and the sum statistics aforementioned. It is the first  max-sum
test for studying cross-sectional dependence for high-dimensional panel data. The advantage is that the test works well for both sparse and non-sparse residuals. Comparing the pros and cons of the max test and the sum test, the max-sum test definitely overcomes both disadvantages.   Our simulations reveal this fact clearly; see Figure \ref{Fig1} and its interpretation at the last part of Section \ref{shenzhen_Maoming}. The max-sum test is particularly useful considering it is hard to quantify or determine  in practice whether a data set is sparse or not.

\section{Theoretical results}

\label{Sec3}

We now present the main theoretical results based on the three types of tests in the order of the sum test, the max test and the max-sum test. Their proofs are presented in Section \ref{Sec5.5}.

\subsection{The limiting distribution for the sum test}\lbl{sum_test_section}

Recall that the sum test described in \eqref{maxsumtest} is a classical one for testing cross-sectional dependence in panel data models. However, the asymptotic theory has not been established yet. Pesaran from \cite{Pesaran2004, pesaran2008a}   conjectures that $S_N$ satisfies the central limit theorem. Some insights on this aspect are given, for example, in \cite{ST12} and \cite{Pesaran2015}. In the following we will present our solution to the problem as well as another one in which the details are given below. The following assumption will be needed throughout the paper. Recall a random variable $V$ is said to be continuous if $P(V=v)=0$ for every $v\in \mathbb{R}.$
\begin{align}\label{assa}
&\textit{Assume}~ \bd{\epsilon}_{i}=(\epsilon_{i1},\cdots,\epsilon_{iT})', i=1,2,\cdots, N, \textit{are independent $T$-dimensional}\nonumber\\ &\textit{random vectors}, ~ \textit{and}~  \epsilon_{i1},\cdots,\epsilon_{iT}
~   \textit{are i.i.d. continuous random variables}\nonumber\\
& \textit{with}~ E\epsilon_{i1}=0 ~\textit{and}~
\mbox{Var}(\epsilon_{i1})=\sigma_i^2>0 ~\textit{for each}~i.
\end{align}
If the $T$ entries of $\bd{\epsilon}_{i}$ are i.i.d. continuous random variables, by using a conditional argument, we then trivially have $P(\bd{a}'\bd{\epsilon}_{i}=0)=0$ for any $\bd{a}\in \mathbb{R}^T\backslash\{\bd{0}\}$. This implies that $P(\bd{M}\bd{\epsilon}_{i}= \bd{0})=0$
for any $l\times T$ matrix $\bd{M}\ne 0$ and $l\geq 1$. So $\hat{\rho}_{ij}$ in \eqref{Asia_foods} is well-defined if $T>p$ because  $\mbox{rank}(\bd{P}_i)=T-p$; see the explanation below \eqref{park_1}.
%

Now we present our solutions to Pesaran's conjectures from \cite{Pesaran2004, pesaran2008a} as follows. For mathematical rigor, we assume the parameter $T$ depends on $N$. The notation $N_k(\mu, \bd{\Sigma})$ stands for the $k$-dimensional multivariate normal distribution with mean vector $\mu$ and covariance matrix $\bd{\Sigma}$. Although the linear regression model in \eqref{haoshichengshuang} requires $p\geq 1$, that is, there are at least one regressors, the following theorem applies to the case that
\begin{align}\lbl{gaoshi_Minnesota}
p=0 \  \mbox{and}\  \bd{P}_i=\bd{I}_T\  \mbox{for each}\ i;
\end{align}
see \eqref{park_1}.

\begin{theorem}\lbl{diligent}
Assume $p\geq 0$ is fixed and $T/\sqrt{N}\to \infty$  as $N\to \infty$.
Let assumption \eqref{assa} hold with $\bd{\epsilon}_{i} \sim N_T(\bd{0}, \sigma_i^2\bd{I})$ for each $i$.
Let $S_N$ be as in \eqref{sumtest} and
$\mu_N$ be as in \eqref{taiyangtaila}, respectively.
Then $(S_N-\mu_N)/N$ converges to $N(0, 1)$ in distribution as $N\to\infty$.
\end{theorem}

From the assumption, it is allowed that the number of cross-sectional units $N$ is much larger than the number of observations $T$, for example, $N$ is of order $T^{3/2}$. We apply the framework of the Lindeberg-Feller martingale CLT to study $S_N$. Although the method is simple and is easy to follow, the technical steps are very involved due to the complex nature of the sample correlation coefficients $\hat{\rho}_{ij}$ in \eqref{Asia_foods}. Many computations focus on the conditional means, variances and higher moments.

Considering a possible better convergence rate than the CLT given in Theorem \ref{diligent}, \cite{pesaran2008a} revises the statistic $S_N$ and proposes a new one as follows.
\begin{align}
\label{R_STY}
Q_N=\sqrt{\frac{2}{N(N-1)}}\sum_{i=1}^{N-1}\sum_{j=i+1}^N
\frac{(T-p)\hat{\rho}_{ij}^2-\mu_{Nij}}{v_{Nij}},
\end{align}
where
\begin{align}
\mu_{Nij}&=\frac{1}{T-p}\mbox{tr}(\bd{P}_i\bd{P}_j),\\
v_{Nij}^2&=a_{1N}\cdot   [\mbox{tr}(\bd{P}_i\bd{P}_j)]^2+2a_{2N}\cdot \mbox{tr}[(\bd{P}_i\bd{P}_j)^2],\\
a_{1N}&=a_{2N}-\frac{1}{(T-p)^2},\\
a_{2N}&=3\cdot \Big[\frac{(T-p-8)(T-p+2)+24}{(T-p+2)(T-p-2)(T-p-4)}\Big]^2.
\end{align}
Pesaran {\it et al}. \cite{pesaran2008a}  conjecture that $Q_N$ also satisfies the CLT. We confirm it in the next theorem.

\begin{theorem}\lbl{Fall_nice} Assume the setting in Theorem \ref{diligent}.
Then $Q_N$ converges to $N(0, 1)$ in distribution as $N\to \infty$.
\end{theorem}
Our simulation in Figure \ref{Fig1} shows that the effects of the two approximations in Theorems \ref{diligent} and \ref{Fall_nice} are too  close to be distinguishable. Theorem \ref{diligent} allows us to perform a level-$\alpha$ test by rejecting the null hypothesis from \eqref{h1}
when $(S_N-\mu_N)/N$ is larger than the $1-\alpha$ quantile $z_{\alpha}= \Phi^{-1}(1-\alpha)$ of $N(0,1)$.
Theorem \ref{Fall_nice} establishes the central limit theorem of $Q_N$ under the same null hypothesis. This  provides a theoretical guarantee for $Q_N$ that has been used
in econometrics; see, for example, \cite{Chudik2013, moscone2009a, pesaran2015testing}.

Now we make some comments.
%

\begin{remark} Let us see how heuristically the mean $\mu_N$ and the standard deviation $N$ in Theorem \ref{diligent} are figured out. In fact, $\mu_N$ is computed by using Lemma \ref{use_tea}(i). The variance is calculated via Lemma \ref{good_tie}(iv) by noticing  $\frac{1}{T}\mbox{tr}[(\bd{P}_i\bd{P}_j)^2]\to 1$ and  $\frac{1}{T^2}[\mbox{tr}(\bd{P}_i\bd{P}_j)]^2\to 1$ as $N\to\infty$ and by regarding  $\{\hat{\rho}_{ij}^2,\, 1\leq i<j \leq N\}$ as independent random variables although they are weakly correlated.
\end{remark}

\begin{remark}\lbl{Re_2} Take $p=0$ in Theorem \ref{diligent}. From \eqref{gaoshi_Minnesota}, $\bd{P}_i=\bd{I}_T$ for all $1\leq i \leq N$ and hence $\mbox{tr}(\bd{P}_i\bd{P}_j)=T$. Theorem \ref{diligent} then says that, assuming  $T/\sqrt{N}\to \infty$,
\begin{align}\lbl{Remark2}
\frac{S_N}{N}-\frac{1}{2}N\to N\Big(-\frac{1}{2}, 1\Big)
\end{align}
in distribution as $N\to\infty.$ This is the trivial case that no linear regression is involved. And $\hat{\rho}_{ij}
=\frac{\bd{\epsilon}_i'\bd{\epsilon}_j}
{\|\bd{\epsilon}_i\|\cdot\|\bd{\epsilon}_j\|}$ from \eqref{Asia_foods}, where $\bd{\epsilon}_i$ and $\bd{\epsilon}_j$ are independent Gaussian vectors with distributions $N_T(\bd{0}, \sigma_i^2\bd{I})$ and $N_T(\bd{0}, \sigma_j^2\bd{I})$, respectively.

Suppose $N/T\to \gamma\in (0, \infty).$ Rewrite \eqref{Remark2} by using the Slutsky lemma to see
\begin{align}\lbl{Fenga}
\sum_{1\le i<j\le N}\hat{\rho}_{ij}^2-\frac{N(N-1)}{2T} \to N(0, \gamma^2)
\end{align}
as $N\to \infty$. This recovers the result by Schott \cite{Schott2005}. A quick reminder is that our assumption ``\,$T/\sqrt{N}\to \infty$" is less stringent than ``$N/T\to \gamma\in (0, \infty)$". A further discussion about Schott's work is continued in the next remark.
\end{remark}

\begin{remark}\lbl{Re_6}  Assume $V\sim N_d(\bd{\mu}, \bd{\Sigma})$, that is, $V$ follows a $d$-dimensional multivariate normal distribution with mean vector $\bd{\mu}$ and covariance matrix $\bd{\Sigma}$. Assume a random sample of size $n$ is given.  Under the assumption $d/n\to c>0$, Schott \cite{Schott2005} studies the null hypothesis that the $d$-entries of $V$ are independent,  that is, $\bd{\Sigma}$ is diagonal. Comparing this with model \eqref{haoshichengshuang}, his model corresponds to \eqref{gaoshi_Minnesota}.
His result is stated in \eqref{Fenga}.
 Jiang \cite{Jiang2019} investigates the same testing problem through the likelihood ratio test  and obtains the CLT for a big class of alternative hypothesis. However, neither derivations of the above two CLTs help the proofs of Theorems \ref{diligent} and \ref{Fall_nice} in this paper.
\end{remark}

The scenario in the following remark is not practical. We consider it purely for mathematical purposes. They serve for further discussions.

\begin{remark}\lbl{Re_3} 
Assume $\bd{x}_1=\cdots=\bd{x}_N$. Then $\bd{P}_1=\cdots =\bd{P}_N$ and $\mbox{tr}(\bd{P}_i\bd{P}_j)=\mbox{tr}(\bd{P}_1)=T-p.$ Hence,
\beaa
\mu_N=\frac{T}{T-p}\cdot \frac{N(N-1)}{2}=N\cdot \Big[\frac{TN}{2(T-p)}-\frac{1}{2}+o(1)\Big].
\eeaa
Theorem \ref{diligent} says that
\beaa
\frac{S_N}{N}-\frac{TN}{2(T-p)}\to N\Big(-\frac{1}{2}, 1\Big)
\eeaa
in distribution as $N\to\infty$. In particular, if $\frac{N}{T}\to c\in (0, \infty)$, then $\frac{TN}{2(T-p)}=\frac{N}{2}+\frac{1}{2}c p+o(1)$. Hence
\begin{align}\lbl{Extreme_CLT}
\frac{S_N}{N}-\frac{1}{2}N\to N\Big(\frac{c p-1}{2}, 1\Big).
\end{align}
A point for this extreme example is that, as these $\bd{x}_i$ are highly correlated,  the CLT in \eqref{Extreme_CLT} is indeed different from the trivial CLT in \eqref{Remark2}. Interestingly, the next example is completely different from this one.
\end{remark}
\begin{remark}\lbl{Re_4} Assume $T=Np$. In this case, $\frac{N}{T}=\frac{1}{p}.$  Construct
\beaa
\bd{x}_1=(\bd{I}_p, \bd{0}_{p\times (T-p)})',\ \bd{x}_2=(\bd{0}_{p\times p}, \bd{I}_p, \bd{0}_{p\times (T-2p)})',\ \cdots, \bd{x}_N=(\bd{0}_{p\times (T-p)}, \bd{I}_{p})'.
\eeaa
They are $T\times p$ matrices. Then $\bd{x}_i'\bd{x}_i=\bd{I}_p$ for each $1\leq i\leq N$, and hence
\beaa
\bd{x}_i(\bd{x}_i'\bd{x}_i)^{-1}\bd{x}_i'=
\begin{pmatrix}
\bd{0}_{p\times p} & \cdots &\bd{0}_{p\times p} & \cdots & \bd{0}_{p\times p}\\
\vdots & \vdots & \vdots & \vdots\\
\bd{0}_{p\times p} & \cdots & \bd{I}_p & \cdots &\bd{0}_{p\times p}\\
\vdots & \vdots & \vdots & \vdots\\
\bd{0}_{p\times p} & \cdots &\bd{0}_{p\times p}& \cdots & \bd{0}_{p\times p}
\end{pmatrix},
\eeaa
where each $\bd{0}_{p\times p}$ is a $p\times p$ submatrix with all entries equal to zero. In other words, we may regard $\bd{x}_i(\bd{x}_i'\bd{x}_i)^{-1}\bd{x}_i'$ as an $N\times N$ matrix with each entry being a block of $p\times p$ matrix, and the only non-zero entry is the $(i, i)$-entry  $\bd{I}_p$. Then $\bd{x}_i(\bd{x}_i'\bd{x}_i)^{-1}\bd{x}_i'\cdot \bd{x}_j(\bd{x}_j'\bd{x}_j)^{-1}\bd{x}_j'=\bd{0}_{T\times T}$ for $i\ne j$. It follows from the definition of $\bd{P}_i$ in \eqref{park_1} that $\mbox{tr}(\bd{P}_i\bd{P}_j)=T-2p$. Thus,
\beaa
\mu_N=\frac{N(N-1)}{2}\cdot \frac{T(T-2p)}{(T-p)^2}=N\cdot \Big[\frac{N-1}{2}+o(1)\Big]
\eeaa
by the assumption $T=Np$. Then
\beaa
\frac{S_N}{N}-\frac{1}{2}N\to N\Big(-\frac{1}{2}, 1\Big)
\eeaa
in distribution as $N\to\infty.$ The essence for this example demonstrates that the projection matrices $\bd{P}_i$  are orthogonal to each other contrary to the highly correlated case in Remark \ref{Re_3}. We see the CLT here is more like the one in the trivial case from Remark \ref{Re_2} but is different from that in Remark \ref{Re_3}.
\end{remark}

\begin{remark}\lbl{Re_5} Review \eqref{dashfune} that $S_N\to \chi^2(d)$ as $T \to\infty$ while $N$ is fixed, where $d=\frac{1}{2}N(N-1)$. By using the approximation $(\chi^2(d)-d)/\sqrt{2d}\to N(0,1)$ as $d\to\infty$, we see that, if taking limit above were legitimate, we would have
%
%
\begin{align*}
\frac{S_N-\frac{1}{2}N(N-1)}{\sqrt{N(N-1)}}\to N(0, 1).
\end{align*}
By the Slutsky lemma, this entails
\begin{align*}
\frac{S_N}{N}-\frac{1}{2}N\to N\Big(-\frac{1}{2}, 1\Big)
\end{align*}
in distribution. It is interesting to see this weak convergence in Remarks \ref{Re_2} and \ref{Re_4}, but not in Remark \ref{Re_3}. In fact, for big data with the feature that two or more parameters are large, to study a statistic of interest, it is not always valid to send parameters to infinity one by one; see such examples in, for instance, \cite{JiangQi2015, JiangYang2013, ZhengBaiYao15}.
\end{remark}



\subsection{The limiting distribution for the max test} Recall model \eqref{haoshichengshuang} and notations in \eqref{park_1}. As in Section \ref{sum_test_section}, we  assume that
$\bd{x}_i\bd{x}_i'$ is invertible for each $1\leq i\leq N$ and the quantity $T$ depends on $N$.
From \eqref{Asia_foods} and assumption \eqref{assa}, we know $\{\hat{\rho}_{ij}:1\leq i,j \leq N\}$ are invariant of $\sigma^2_1,\cdots,\sigma^2_N$, so
we are able to assume, without loss of generality, $\sigma_1=\cdots=\sigma_N\equiv 1$. Review $\epsilon_{11}$ in assumption \eqref{assa} and $L_N$ in \eqref{Write}.  As explained in \eqref{gaoshi_Minnesota}, we will also consider the case $p=0$. The main results in this section are presented as follows.

\begin{theorem}\lbl{linear_case} Assume $p\geq 0$ is fixed and $\lim_{N\to\infty}T/N=c\in (0, \infty)$. Let $\bd{\epsilon}_{1},\cdots, \bd{\epsilon}_{N}$ be i.i.d. and  assumption \eqref{assa} hold with $E|\epsilon_{11}|^{\tau}<\infty$ for some $\tau>8$.   Then, as $N\to\infty$,  $TL_N^2-4\log N +\log\log N$ converges weakly to the distribution function $F(y)= \exp(-e^{-y/2}/\sqrt{8\pi})$, $y \in \mathbb{R}$.
\end{theorem}

\begin{theorem}\lbl{exponential_case} Assume $p\geq 0$ is fixed and $\log N=o(T^{1/5})$ as $N\to \infty$. Let $\bd{\epsilon}_{1},\cdots, \bd{\epsilon}_{N}$ be i.i.d. and assumption \eqref{assa} hold with $Ee^{\omega|\epsilon_{11}|}<\infty$ for some $\omega>0.$  Then, as $N\to\infty$,  $TL_N^2-4\log N +\log\log N$ converges weakly to the distribution function $F(y)= \exp(-e^{-y/2}/\sqrt{8\pi})$, $y \in \mathbb{R}$.
\end{theorem}

We say $\xi$ is a {\it subgaussian} random variable if there exists $\sigma>0$ such that $Ee^{t\xi}\leq e^{\sigma^2t^2/2}$ for all $t\in \mathbb{R}.$ By the Markov inequality, it is easy  to see $P(|\xi|\geq x)\leq 2e^{-x^2/(2\sigma^2)}$ for all $x>0$. As a consequence, $Ee^{\theta \xi^2} < \infty$ for all $\theta < \frac{1}{2\sigma^2}$. Obviously, bounded random variables and Gaussian random variables are all subgaussian random variables.

\begin{theorem}\lbl{super_exponential_case} Assume $p\geq 0$ is fixed and $\log N=o(T^{1/3})$ as $N\to \infty$. Let $\bd{\epsilon}_{1},\cdots, \bd{\epsilon}_{N}$ be i.i.d. and  assumption \eqref{assa} hold with  $\epsilon_{11}$ being a subgaussian random variable. Then, as $N\to\infty$,  $TL_N^2-4\log N +\log\log N$ converges weakly to the distribution function $F(y)= \exp(-e^{-y/2}/\sqrt{8\pi})$, $y \in \mathbb{R}$.
\end{theorem}

The strategy of the proofs of Theorems \ref{linear_case}-\ref{super_exponential_case} is to approximate $L_N=\max_{1\le i<j\le N} |\hat{\rho}_{ij}|$ for any $p\geq 0$ by $L_N$ for the case $p=0$, in which the limiting behavior is understood in \cite{CJF13}. We have to show the difference between the two versions of $L_N$ is negligible.

Theorems \ref{linear_case}-\ref{super_exponential_case} indicate that we get the same asymptotic distribution of the max statistic $TL_N^2-4\log N +\log\log N$
under different moment assumptions. This allows us to have a flexibility to work on different pairs of $(N, T)$. Under null hypothesis \eqref{h1} and the assumptions imposed in the above three theorems, we conclude that  a level-$\alpha$ test by rejecting the null hypothesis when
$TL_N^2-4\log N+\log\log N$ is larger than the $1-\alpha$ quantile $q_{\alpha}= -\log(8\pi)-2\log\log(1-\alpha)^{-1}$ of $F(y)$.


\subsection{The limiting distribution for the max-sum test} Review the accounts before the statement of Theorem \ref{diligent}. We have the following conclusion on asymptotic independence.

\begin{theorem}\lbl{Asym_indept}
Let $S_N$, $L_N$ and $\mu_N$ be as in \eqref{sumtest}, \eqref{Write} and \eqref{taiyangtaila}, respectively. Under the same assumptions as in Theorem \ref{diligent}, we have
that $(S_N-\mu_N)/N$ and $TL_N^2-4\log N +\log\log N$ are asymptotically independent as $N\to \infty$.
\end{theorem}

By employing a new trick we prove the asymptotic independence between the maximum and the sum of random variables in Theorem \ref{Asym_indept}. This method is expected to be used in many of such type of problems. In fact, there are few literature to prove asymptotic independence between  sums of and maxima of random variables. Some close references are \citep{BBQ98, Xu2016}. The method here is new.  It gives a novel tool to establish asymptotic independence between sums of and maxima of random variables.

 To understand the idea quickly, we start with the set-up \eqref{cuba}. The first observation is that the maximum of many random variables, seemingly a global property, can be understood from their local property, that is, the maxima of subsets of random variables with fixed sizes. This step is done through the inclusion-exclusion formula. The second observation is that any such subset of random variables and the sum are  independent with very high probability. Consequently  the probability of the intersection of  the events related to local maxima and the sum can be written as the product of two individual probabilities. Then we use the inclusion-exclusion formula one more time to get the product of two individual probabilities up to negligible errors. Details are shown at the beginning of Section \ref{Key_asym_indept}.


%


An immediate application is given below. By Theorems \ref{diligent} and \ref{super_exponential_case}, we know that
\begin{align}
&\frac{1}{N}(S_N-\mu_N)\to N(0, 1)\ \mbox{weakly};\ \ \ \ \ \ \lbl{alp2}\\
&TL_N^2-4\log N +\log\log N \to F(y)=\exp(-e^{-y/2}/\sqrt{8\pi})\ \mbox{weakly}.\lbl{alp1}
\end{align}
Let $\Phi(x)$ be the distribution function of $N(0,1).$ Trivially, both $F(y)$ and $\Phi(x)$ are continuous functions.
Set $P_{S_N}=1-\Phi\{(S_N-\mu_N)/N\}$ and $P_{L_N}=1-F(TL_N^2-4\log N+\log\log N)$. By Theorem \ref{Asym_indept}, \eqref{alp2} and \eqref{alp1}, we see that $P_{L_N}$ and $P_{S_N}$ are asymptotically independent and each limit is $U[0, 1]$, the uniform distribution over $[0, 1].$ So the following holds easily.

\begin{coro}\lbl{coro2} Set $C_N=\min\{P_{S_N}, P_{L_N}\}$. Assume the setting in Theorem \ref{Asym_indept}. Then
$C_N$ converges to $W:=\min\{U,V\}$ in distribution as $N\to\infty$, where $U$ and $V$ are i.i.d. random variables with distribution $U[0, 1].$ The distribution function of $W$ is given by $G(w)=2w-w^2$ for $w \in [0, 1]$.
\end{coro}

According to Corollary \ref{coro2}, the proposed max-sum test in \eqref{maxsumtest} allows us to perform a level-$\alpha$ test by rejecting the null hypothesis \eqref{h1} if
 $C_N<1-\sqrt{1-\alpha}$.

\section{Simulation studies}\label{Sec4}

We now conduct simulations to compare the finite sample performance of the tests studied in this paper and another test in literature. The tests we have worked on in this paper are based on
$S_N$, $L_N$,    $C_N$, $Q_N$   in    \eqref{sumtest}, \eqref{Write},  \eqref{maxsumtest}, \eqref{R_STY},  respectively.
The other one is  based on
$CD$ from \cite{Pesaran2004} defined by
\begin{align}\lbl{cotton_air}
CD=\sqrt{\frac{2T}{N(N-1)}}\sum_{1\leq i<j \leq N}\hat{\rho}_{ij}.
\end{align}
Notice $Q_N$ here is the same as the notation  LM$_{\textrm{adj}}$ from \cite{pesaran2008a}.
%
%
%
In the following we will explain our simulation designs and state our simulation findings.

\subsection{Simulation designs}\lbl{S_design}

We consider the data generating process used in \cite{pesaran2008a}, which is specified as
\begin{align*}
y_{it}=\alpha_i+\sum_{l=2}^p x_{lit}\beta_{li}+\epsilon_{it}
\end{align*}
for $i=1, \cdots, N$ and $t=1, \cdots, T$. Comparing the notations in  model \eqref{haoshichengshuang}, we have $x_{it}=(1, x_{2it}, \cdots, x_{pit})^T\in \mathbb{R}^{p}$ and $\beta_i=(\alpha_i, \beta_{2i}, \cdots, \beta_{pi})\in \mathbb{R}^{p}.$

Now we independently generate  $\alpha_i\sim N(0,1)$ and $\beta_{li}\sim N(1,0.04)$. The covariates are generated by
\begin{align*}
x_{lit}=0.6x_{li t-1}+v_{lit}
\end{align*}
for $i=1,\cdots, N,~t=-50, -49, \cdots, T$ and $l=2,\cdots, p$ with $x_{li,-51}=0$, where $v_{lit}\sim N(0,\zeta_{li}^2/(1-0.6^2))$ and $\zeta_{li}^2\sim \chi^2_6/6$. In this case,  $\zeta_{li}^2$'s  are independently sampled first, then $v_{lit}$'s are independently generated by conditioning on the values of  $\zeta_{li}^2$.

Now we generate $\epsilon_{it}$'s under null hypothesis \eqref{h1}. Let $\epsilon_{it}=\sigma_i w_{it}$, where $w_{it}$'s are generated from three different distributions: (i) $N(0,1)$; (ii) $t_6/\sqrt{6/4}$; (iii) $(\chi^2_5-5)/\sqrt{10}$. Here $t_d$ is the $t$-distribution of degree $d$ and $\chi^2_d$ is the chi-square distribution of degree $d.$ The normalization in (ii) and (iii) is such that each new random variable has mean zero and variance one. Let $\sigma_i^2\sim \frac{p}{2}\chi^2_2$, as in the dynamic setup of \cite{pesaran2008a}.

We turn to produce data under the alternative hypothesis. Let $\bd{\eta}_t:=(\eta_{1t},\cdots,\eta_{Nt})'$ be generated from the above three different distributions under the null hypothesis. Set  $\bd{\epsilon}_{.t}=(\epsilon_{1t},\cdots,\epsilon_{Nt})'=\bms^{1/2}\bd{\eta}_t$. Please differentiate the notation $\bd{\epsilon}_{.t}$ here and $\bd{\epsilon}_{i}$ in \eqref{shadongxi}.
We consider the following two cases of the covariance matrix $\bms=\D^{1/2}\R\D^{1/2}$ with $\D=\diag\{\sigma^2_1,\cdots,\sigma^2_N\}$.
\begin{itemize}
\item[(1)] Non-sparse case. Randomly select a subset $A\subset\{1,\cdots,N\}$ with cardinality $N^{0.5}$. Let $\R=(\rho_{ij})_{N\times N}$ be a symmetric matrix with
$\rho_{ij}=1$ if $i=j$. For $i< j$, define $\rho_{ij}=0$ if $i\not \in A$ or $j\not \in A$, and  $\rho_{ij}$ has the uniform distribution over $(\sqrt{3\log N/T}, \sqrt{5\log N/T})$ if $i\in A$ and $j \in A$.
\item[(2)] Sparse case. Randomly select a subset $A\subset\{1,\cdots,N\}$ with cardinality $N^{0.3}$. Let $\R=(\rho_{ij})_{N\times N}$ be a symmetric matrix with
$\rho_{ij}=1$ if $i=j$.  For $i< j$, define $\rho_{ij}=0$ if $i\not \in A$ or $j\not \in A$, and  $\rho_{ij}$  has the uniform distribution over $(\sqrt{8\log N/T}, \sqrt{10\log N/T})$ if $i\in A$ and $j \in A$.
\end{itemize}
To ensure that the covariance matrix  $\bms=\D^{1/2}\R\D^{1/2}$ is positive definite, we replace the correlation matrix $\R$ with $\R+\lambda \I_N$, where $\lambda:=|\lambda_{\min}(\R)|+0.05$ and
$\lambda_{\min}(\R)$ is the minimum eigenvalue of $\R$.
Then, we consider two choices of the sample size $T=50,100$, and three choices of the dimension $N=50,100,200$.

\subsection{Simulation results}\lbl{shenzhen_Maoming}

\begin{figure}
  \centering
  \includegraphics[width=10cm]{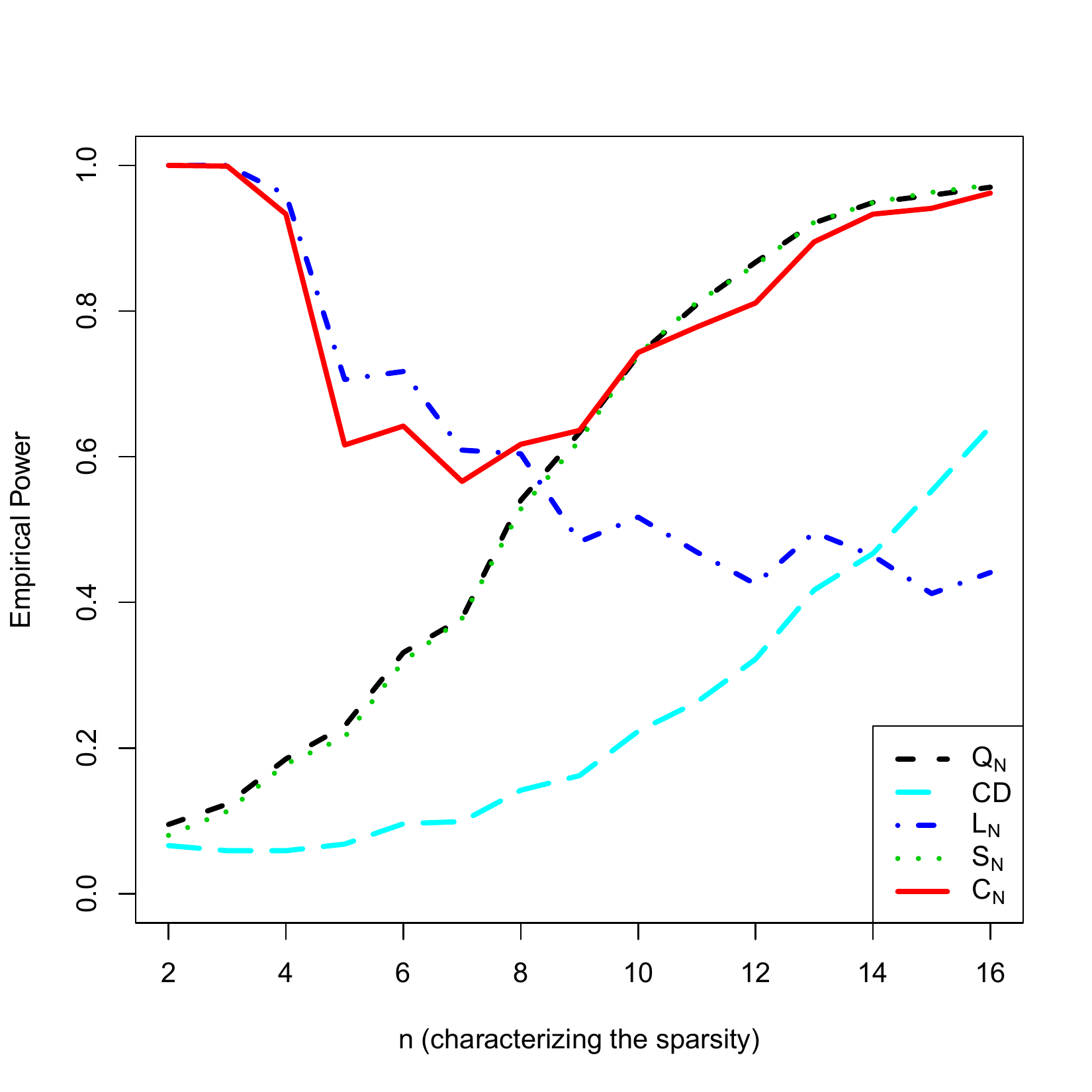}\\
  \caption{Empirical power curves of tests vary with $n$. The number $n$  characterizes
  the degree of sparsity.
  The larger the degree $n$ is, the lower the sparsity is.} \label{Fig1}
\end{figure}

We now present simulation results on the tests of $S_N$, $L_N$,    $C_N$, $Q_N$, $CD$   in    \eqref{sumtest}, \eqref{Write},  \eqref{maxsumtest}, \eqref{R_STY}, \eqref{cotton_air}, respectively.
All the conclusions are based on 1,000 replications.
The empirical sizes and  powers of these tests in non-sparse and sparse cases are summarized in Tables \ref{tab:t1} to \ref{tab:t3}. The power curves are plotted in Figure \ref{Fig1}. We next analyze them in detail.

Table \ref{tab:t1} indicates that all methods have empirical sizes not much larger than 5\%.
Here, the max test $L_N$ and the max-sum test $C_N$  tend to have smaller empirical sizes  than the remaining ones, especially as $T$ is
relatively small. This is not very surprising because it is common for maximum methods designed for raw data models; see, for example, \citep{Liu2008}.

Tables \ref{tab:t2} and \ref{tab:t3} show the information of empirical powers in both non-sparse  and sparse cases. Review the sum test $S_N$ and the sum-based test $Q_N$ are originally proposed in \cite{breusch1980the} and \cite{pesaran2008a}, respectively. The two are well studied in this paper. Tables \ref{tab:t2} and \ref{tab:t3}
 show  that
 %
$S_N$ and $Q_N$ perform best in non-sparse cases in terms of empirical powers, but very poorly in sparse cases.
On the contrary, the proposed max test $L_N$ performs the best in sparse cases, but very poorly in dense cases.
Interestingly, it can be seen from Figure \ref{Fig1} that the empirical power performance of our proposed max-sum test $C_N$ is always very close to the optimal one among
all of the tests, regardless of the local alternative being sparse or not. This shows a very appealing property for the test $C_N$ which compromises the tests for residuals in both sparse and non-sparse cases. In fact, it is hard to tell in reality if residuals are sparse or not.

\begin{table}[htbp]
\footnotesize
\begin{center}
\caption{\label{tab:t1} Empirical sizes (\%) of tests.}
                     \vspace{0.5cm}
                     \renewcommand{\arraystretch}{1}
                     \setlength{\tabcolsep}{7pt}{
\begin{tabular}{cc|ccc|ccc|ccc}
\hline \hline
& \multicolumn{1}{c}{$p$}  & \multicolumn{3}{c}{{2}} & \multicolumn{3}{c}{{3}} & \multicolumn{3}{c}{{ $4$}}\\ \hline
T& $N$ & 50&100&200 & 50&100&200& 50&100&200 \\
\hline
\hline
\multicolumn{11}{c}{Normal distribution}\\\hline
$50$&$Q_N$&6.0&5.8&4.3&4.9&6.3&7.0&5.3&5.5&6.7\\
&CD&6.0&5.5&5.3&4.6&4.6&5.1&4.8&5.5&5.9\\
&$L_N$&1.1&0.7&0.3&1.2&0.5&0.2&1.5&0.1&0.4\\
&$S_N$&5.9&5.7&4.1&4.8&5.7&5.8&5.2&5.5&5.3\\
&$C_N$&2.8&2.2&1.6&2.9&3.2&2.7&2.7&2.7&3.8\\ \hline
$100$&$Q_N$&4.2&3.5&5.6&4.8&5.3&5.8&5.9&4.4&6.9\\
&CD&4.9&4.3&5.3&5.2&5.6&4.8&6.4&4.2&5.0\\
&$L_N$&2.6&2.0&1.3&2.3&1.5&1.0&1.9&1.6&2.0\\
&$S_N$&4.1&3.4&5.6&4.5&5.0&5.2&6.2&3.8&6.5\\
&$C_N$&3.4&2.8&3.0&3.7&3.1&3.5&4.2&2.5&4.5\\ \hline
\multicolumn{11}{c}{$t_6$-distribution}\\\hline
$50$&$Q_N$&5.2&5.9&5.8&5.0&6.9&5.9&6.7&6.1&7.9\\
&CD&4.5&5.0&5.4&4.6&4.3&4.4&4.9&4.2&3.8\\
&$L_N$&1.4&2.5&1.4&1.6&1.3&0.8&2.3&1.4&1.1\\
&$S_N$&4.6&5.4&5.2&4.7&6.8&5.4&5.8&5.6&6.0\\
&$C_N$&3.7&4.8&4.2&3.2&4.1&3.1&4.2&4.0&4.2\\ \hline
$100$&$Q_N$&6.5&6.1&5.3&6.8&5.8&4.6&7.0&6.4&5.0\\
&CD&5.9&5.7&5.0&5.4&6.3&6.4&5.9&6.3&6.0\\
&$L_N$&3.5&4.0&4.3&3.7&4.0&3.9&3.9&4.7&3.2\\
&$S_N$&5.9&5.8&4.8&6.5&5.9&4.8&6.7&6.3&5.0\\
&$C_N$&4.8&5.4&4.4&4.7&4.9&4.1&5.7&5.7&4.1\\\hline
\multicolumn{11}{c}{$\chi^2_5$-distribution}\\\hline
$50$&$Q_N$&5.8&5.5&6.4&6.8&4.7&6.5&7.5&5.1&7.3\\
&CD&3.9&4.8&4.1&6.1&3.9&4.6&5.3&4.4&4.2\\
&$L_N$&2.2&1.9&1.8&2.6&2.1&1.0&3.1&2.2&1.3\\
&$S_N$&5.3&5.3&6.1&6.6&3.7&5.4&6.6&4.4&6.1\\
&$C_N$&3.6&4.2&3.5&4.3&3.1&3.9&3.9&3.4&3.7\\ \hline
$100$&$Q_N$&5.7&5.8&4.5&5.0&5.7&5.8&6.2&5.8&5.1\\
&CD&5.0&5.3&5.0&4.0&4.7&4.7&4.1&5.9&5.0\\
&$L_N$&5.1&3.9&5.8&6.3&4.6&5.4&4.8&6.6&5.7\\
&$S_N$&5.4&5.7&4.2&4.9&5.6&5.7&5.9&5.2&5.3\\
&$C_N$&5.1&5.8&5.2&4.9&5.3&5.1&5.1&5.3&5.5\\
\hline
\hline
\end{tabular}}
\end{center}
\end{table}

\begin{table}[htbp]
\footnotesize
\begin{center}
\caption{\label{tab:t2} Empirical powers (\%) of tests in non-sparse cases.}
                     \vspace{0.5cm}
                     \renewcommand{\arraystretch}{1}
                     \setlength{\tabcolsep}{7pt}{
\begin{tabular}{cc|ccc|ccc|ccc}
\hline\hline
& \multicolumn{1}{c}{$p$}  & \multicolumn{3}{c}{{2}} & \multicolumn{3}{c}{{3}} & \multicolumn{3}{c}{{ $4$}}\\ \hline
T& N &  50&100&200 & 50&100&200& 50&100&200 \\
\hline
\hline
\multicolumn{11}{c}{Normal distribution}\\\hline
$50$&$Q_N$&84.7&96.5&99.9&80.1&97.2&99.9&75.9&95.1&99.3\\
&CD&34.9&50.2&62.8&38.1&49.7&63.2&36.2&50.3&59.9\\
&$L_N$&53.4&77.1&91.9&44.9&74.5&97.8&40.7&65.2&92.6\\
&$S_N$&84.0&96.3&100&80.2&97.0&99.9&75.8&95.2&99.6\\
&$C_N$&80.1&95.5&99.7&75.9&96.3&99.6&72.3&93.6&99.5\\ \hline
$100$&$Q_N$&76.6&92.7&98.5&75.0&88.9&98.1&76.9&89.8&97.2\\
&CD&33.9&43.6&54.3&32.7&42.4&51.8&33.4&38.7&52.1\\
&$L_N$&52.9&72.0&90.9&45.7&66.6&90.2&54.7&65.8&86.4\\
&$S_N$&76.1&92.7&98.3&74.7&88.6&98.2&76.4&89.5&97.3\\
&$C_N$&74.3&92.1&98.3&70.9&88.5&98.5&76.5&90.1&97.2\\ \hline
\multicolumn{11}{c}{$t_6$-distribution}\\\hline
$50$&$Q_N$&83.4&96.9&99.9&82.8&94.8&99.9&83.3&94.2&98.9\\
&CD&39.1&52.8&61.1&40.6&48.9&63.0&37.9&44.9&59.6\\
&$L_N$&57.9&80.6&97.6&58.1&76.1&94.5&56.9&70.3&89.9\\
&$S_N$&82.7&96.6&100&82.3&94.3&99.8&83.2&94.0&98.8\\
&$C_N$&81.0&95.7&99.9&81.2&93.4&99.8&80.5&92.8&98.4\\ \hline
$100$&$Q_N$&75.2&92.8&96.7&76.3&87.4&98.3&72.6&88.4&97.3\\
&CD&31.6&46.3&49.0&32.7&43.6&53.3&30.8&43.1&52.0\\
&$L_N$&55.0&77.6&83.7&51.5&67.0&90.0&50.8&69.9&88.4\\
&$S_N$&74.4&92.6&96.6&75.5&87.0&98.3&72.3&88.0&97.5\\
&$C_N$&74.5&92.5&95.7&73.6&87.4&98.0&71.5&86.6&97.8\\ \hline
\multicolumn{11}{c}{$\chi^2_5$-distribution}\\\hline
$50$&$Q_N$&86.3&97.4&98.9&84.4&95.9&99.9&77.0&94.4&99.5\\
&CD&39.2&51.6&60.7&38.1&53.2&61.5&36.2&47.7&61.2\\
&$L_N$&62.3&84.1&95.2&59.3&82.3&96.7&46.2&71.1&96.5\\
&$S_N$&86.0&97.4&99.3&84.2&95.7&99.9&76.9&94.5&99.5\\
&$C_N$&85.4&97.1&98.9&82.9&96.0&99.7&73.0&93.7&99.7\\ \hline
$100$&$Q_N$&77.8&86.9&97.2&77.1&86.9&97.0&70.9&90.0&95.2\\
&CD&34.8&40.5&50.2&34.9&41.7&52.2&31.8&43.1&49.9\\
&$L_N$&56.5&69.6&91.1&57.7&71.7&88.2&49.2&75.8&86.7\\
&$S_N$&77.2&86.9&97.3&76.4&86.8&97.0&70.3&89.9&95.5\\
&$C_N$&75.6&86.6&97.4&77.2&86.8&98.0&70.2&90.7&96.6\\
\hline
\hline
\end{tabular}}
\end{center}
\end{table}

\begin{table}[htbp]
\footnotesize
\begin{center}
\caption{\label{tab:t3} Empirical powers (\%) of tests in sparse cases.}
                     \vspace{0.5cm}
                     \renewcommand{\arraystretch}{1}
                     \setlength{\tabcolsep}{7pt}{
\begin{tabular}{cc|ccc|ccc|ccc}
\hline \hline
& \multicolumn{1}{c}{$p$}  & \multicolumn{3}{c}{{2}} & \multicolumn{3}{c}{{3}} & \multicolumn{3}{c}{{ $4$}}\\ \hline
T& $N$ & 50&100&200 & 50&100&200& 50&100&200 \\
\hline
\hline
\multicolumn{11}{c}{Normal distribution}\\\hline
$50$&$Q_N$&44.1&20.0&33.4&36.9&26.3&28.3&34.6&24.3&27.7\\
&CD&10.4&6.70&7.00&8.00&7.50&8.00&8.90&6.60&7.50\\
&$L_N$&99.7&99.8&100&98.7&100&100&96.0&100&100\\
&$S_N$&43.0&19.1&31.7&36.0&26.1&28.7&33.8&23.7&25.6\\
&$C_N$&99.6&99.0&100&97.8&100&100&93.2&99.9&100\\ \hline
$100$&$Q_N$&23.1&16.5&19.2&26.2&13.3&20.3&29.6&14.6&20.4\\
&CD&7.70&6.80&5.00&9.40&5.80&6.30&7.10&5.60&7.20\\
&$L_N$&75.5&87.9&98.3&79.3&76.8&99.5&80.5&79.9&98.4\\
&$S_N$&22.8&16.3&19.3&25.6&13.1&19.8&28.9&14.2&19.9\\
&$C_N$&70.4&84.3&97.4&74.3&69.9&99.0&75.9&74.6&98.0\\ \hline
\multicolumn{11}{c}{$t_6$-distribution}\\\hline
$50$&$Q_N$&41.8&19.2&29.9&32.2&26.4&30.1&31.0&23.7&27.6\\
&CD&9.00&8.00&8.00&9.10&6.10&7.70&8.90&9.50&7.20\\
&$L_N$&96.9&98.7&100&93.5&99.9&100&92.7&100&100\\
&$S_N$&39.7&18.5&28.1&32.3&24.9&28.8&30.0&22.9&25.9\\
&$C_N$&95.3&97.4&100&91.3&99.7&100&90.8&100&100\\ \hline
$100$&$Q_N$&26.8&18.9&21.7&23.7&14.7&20.8&27.6&14.8&19.4\\
&CD&8.00&6.10&6.90&7.00&6.30&5.90&7.80&5.50&5.30\\
&$L_N$&80.5&89.9&96.5&74.2&84.8&98.7&80.8&83.8&96.9\\
&$S_N$&25.5&18.3&21.2&23.0&14.6&20.9&27.3&14.7&18.5\\
&$C_N$&76.1&85.1&94.7&68.3&80.4&97.6&77.7&79.7&95.8\\ \hline
\multicolumn{11}{c}{$\chi^2_5$-distribution}\\\hline
$50$&$Q_N$&37.0&25.0&29.8&35.8&24.2&28.7&34.5&21.0&27.1\\
&CD&9.60&6.80&8.50&6.70&7.40&8.70&10.3&6.00&7.60\\
&$L_N$&98.7&99.8&100&95.0&99.8&100&90.7&97.1&100\\
&$S_N$&36.9&24.2&29.3&35.1&24.4&27.3&34.0&20.0&25.5\\
&$C_N$&97.5&99.6&100&93.4&99.8&100&87.7&95.9&100\\ \hline
$100$&$Q_N$&28.3&17.3&21.2&27.3&14.6&22.1&26.1&16.9&20.9\\
&CD&9.50&5.60&6.80&7.50&4.90&5.50&7.00&4.70&7.70\\
&$L_N$&80.3&89.0&98.8&73.6&85.5&98.2&78.8&82.6&96.4\\
&$S_N$&27.5&16.9&20.7&26.4&14.3&21.7&25.8&16.3&20.4\\
&$C_N$&75.4&83.3&97.9&69.5&81.2&97.6&73.3&76.1&93.9\\
\hline
\hline
\end{tabular}}
\end{center}
\end{table}

Figure \ref{Fig1} shows the changes of the powers  of all the tests as the degree of sparsity changes.
Now we explain the procedure to generate
the empirical power curves in Figure \ref{Fig1}. In fact, the horizonal direction  in the plot is $n$, the degree of sparsity
to be defined; the vertical direction represents powers. Specifically, the simulation  is designed as follows. Review the general simulation design in Section \ref{S_design}. We take
$T=50$, $N=200$, $p=2$, $n=2,\cdots,16$; $w_{it}$ are generated from normal distributions;
a subset $A\subset\{1,\cdots,N\}$ is randomly selected with cardinality $n$; $\R=(\rho_{ij})_{1\le i,j\le N}$, where
$\rho_{ij}=1$ if $i=j$; for $i\neq j$, set $\rho_{ij}=0$ if $i\not \in A$ or $j\not \in A$,
and $\rho_{ij}$  has the uniform distribution over $(\sqrt{8(\log n)^{-1}\log N/T}, \sqrt{10(\log n)^{-1}\log N/T})$ if $i\in A$ and $j \in A$. So a larger $n$ means a lower level of sparsity.

Figure \ref{Fig1} indicates that the empirical power of the max-sum test $C_N$ is always very close to the maximum power of
all tests for all $n$. By contrast, the empirical power curves of the remaining methods are monotone,
{\it i.e.},  the empirical powers of both $S_N$ and $Q_N$ generally increase with the decrease of sparsity. On the contrary, the empirical power of the max test increases
with the increase of sparsity. However, every test excluding the max-sum test $C_N$, favors  either the sparse case or the non-sparse case, not both cases simultaneously.

\section{Application}

\label{Sec5}

In this section, we apply the five tests to the securities in the Standard $\&$ Poor (S$\&$P) 500 index of large cap U.S. equity market. As seen earlier, they are $S_N$, $L_N$,    $C_N$, $Q_N$, $CD$   in    \eqref{sumtest}, \eqref{Write},  \eqref{maxsumtest}, \eqref{R_STY}, \eqref{cotton_air}, respectively. This
demonstrates the practical usefulness of the proposed tests. The S$\&$P 500 index is primarily intended as a leading indicator of U.S. equities. The
composition of this index is monitored by Standard and Poor to ensure the widest possible overall market representation
while reducing the index turnover to a minimum. In this section, we consider 374 securities that have been included in the S$\&$P 500 index during
the whole period from January 2005 to November 2018.

In particular, the panel data on the safe rate of return, and the market factors are obtained from Ken French's data library web page. The one-month US treasury bill rate is chosen as the risk-free rate ($r_{ft}$), the value-weighted return on all NYSE, AMEX, and NASDAQ stocks from CRSP is used as a proxy for the market return ($r_{mt}$), the average return on the three small portfolios minus the average return on the three big portfolios ($SMB_t$), and the average return on two value portfolios minus the average return on two growth portfolios ($HML_t$). SMB and HML are based on the stocks listed on the NYSE, AMEX and NASDAQ. All data are measured in percent per month. During January 2005 to November 2018, a total of 163 consecutive observations are obtained.

The Fama-French three-factor model \citep{Fama1993Common} is given as follows:
\[ y_{it}=r_{it}-r_{ft}=\beta_{0i}+\beta_{1i} (r_{mt}-r_{ft})+\beta_{2i}SMB_t+\beta_{3i} HML_t+\epsilon_{it} \]
for each $1\leq i\leq N$ and $1\leq t\leq T$ with $N=374$. We are interested in the following null hypothesis:
\begin{align*}
 H_0: \epsilon_{11},\cdots,\epsilon_{N1}\text{ are independent}.
\end{align*}
That is, we are testing that the  374 variables are independent.

Now we evaluate the performance of the five tests in Section \ref{shenzhen_Maoming}, that is,
$S_N$, $L_N$,    $C_N$, $Q_N$, $CD$   in    \eqref{sumtest}, \eqref{Write},  \eqref{maxsumtest}, \eqref{R_STY}, \eqref{cotton_air}, respectively.
We randomly sample $T=15,25,35$
observations from the 163 monthly returns. At each value of $T$, the experiment is repeated 1000 times. It is trivial to see
\beaa
\binom{163}{15}=\frac{163\cdot 162\cdots 149}{15\cdot 14\cdots 2}>\frac{150^{14}\cdot 149}{15^{14}}>10^{15}.
\eeaa
Similarly,
\beaa
\binom{163}{25}>5\cdot 10^{18}\ \ \mbox{and}\ \ \ \binom{163}{35}>2\cdot 10^{20}.
\eeaa
This says that, although there is a dependency when sample $15$ numbers from a total of $163$ numbers for $1000$ times, comparing to $10^{15}$, the number of repeats $1000$ is still reasonable.  The same also applies to the cases $T=25$ and $T=35$.

The results are summarized in Table \ref{t1}. It
suggests that all tests except the max test always reject the null hypothesis of cross-sectional independence. So this indicates the definite cross-sectional
dependence among stock returns under the  three-factor model by Fama-French. In particular, the max test rejects the null hypothesis when $T$ grows to
35, but never reject it when $T$ reduces to 15. To understand this phenomenon, we point out a well known fact that there may exist a large number of
underlying dependencies between stocks in the same industry or relevant industries. This  leads us to believe that this is indeed a non-sparse case in which
the sum and max-sum tests are more valid.

\section{Concluding remarks}
\label{Sec6}
In this paper we study three tests:  the sum test, the max test and the max-sum test, where the latter two are new ones. Two conjectures on the sum test have been settled. A new method to show asymptotic independence between the maximum  and the sum of squares of a given set of random variables is established. Now we make some comments.

1. Under the Gaussian assumption, we obtain the CLTs for $S_N$ in Theorems \ref{diligent} and \ref{Fall_nice}. However, the Gaussian assumption is  not needed in the study on the maxima of sample correlations in Theorems \ref{linear_case}, \ref{exponential_case}  and \ref{super_exponential_case}. One question is whether the Gaussian assumption can be removed from Theorems \ref{diligent} and \ref{Fall_nice}. Our proofs rely on the framework in Lemma \ref{long_for} where the normal assumption is essential. Another question is about the restriction between $N$ and $T$ in Theorems \ref{diligent} and \ref{Fall_nice}. Can the assumption ``$T/\sqrt{N}\to \infty$" be relaxed? What is the behavior of $S_N$ for other regimes of relationship between $N$ and $T$?

2. The linear regression in  \eqref{haoshichengshuang}  is one of many panel data models; see, for example, the book length treatment in \cite{Baltagi2013}, \cite{Hsiao2014}, \cite{Pesaran2015}, \cite{Wooldridge2010}, among others. Some of other models can be studied similarly for the properties we have pursued in this paper. We leave them as a future work to our authors.

3. A new way is established to show the asymptotic independence between the sum of and the maximum of a set of random variables. The detail of the method is elaborated at the beginning of Section \ref{Key_asym_indept}. We expect this method will also work for other set of random variables of similar feature.

4. For the sum $S_N=\sum_{1\le i<j\le N}T\hat{\rho}_{ij}^2$, Theorem \ref{diligent} states that the central limit theorem of $S_N$ involves with projection matrices $\bd{P}_i$  defined via data; see \eqref{park_1}. However, interestingly enough, as shown in Theorems \ref{linear_case}, \ref{exponential_case} and \ref{super_exponential_case}, the behavior of  the maximum $L_N=\max_{1\le i<j\le N} |\hat{\rho}_{ij}|$ does not depend on $\bd{P}_i$. Only parameters $N$ and $T$ participate in the limiting process.

5. As seen in Section \ref{shenzhen_Maoming}, a simulation study is carried out for tests based on
$S_N$, $L_N$,    $C_N$, $Q_N$, $CD$   in    \eqref{sumtest}, \eqref{Write},  \eqref{maxsumtest}, \eqref{R_STY},  \eqref{cotton_air}, respectively.
It shows that the max-sum test is always
 very close to the maximum power of
all tests for both sparse and non-sparse residuals. By contrast, the empirical power curves of other  methods favor only for one of the two types of residuals. In fact, in practice, it is hard to differentiate if a set of numbers is sparse or not. This implies the max-sum test is also desirable for other statistical models as long as three things are known: the central limit theorem holds for the sum of  a set of random variables; the maximum of the set of random variables is   the Gumbel distribution asymptotically; the sum and the maximum are asymptotically independent.

\begin{table}[htbp]
 \centering
 \caption{The rejection rates of testing cross-sectional independence for the S\&P stock panel data, where $N=374$ and $T=15,25,35$.
For each $T$, we sample 1000 data sets. }
 \vspace{0.2cm}
 \begin{tabular}{cccccc}\hline \hline
    & $T=15$ & $T=25$ & $T=35$ \\
 $Q_N$& 1 & 1 & 1\\
 $CD$    & 1 & 1 & 1\\
 $L_N$ & 0 & 0.08 & 1\\
 $S_N$ & 1 & 1 & 1\\
 $C_n$ & 1 & 1 & 1\\
\hline \hline
\end{tabular}\label{t1}
\end{table}

\section{Proof}
\label{Sec5.5}
%
There are three subsections in this part.
In each subsection we first accumulate some first hand or second hand of understanding before the proofs of main theorems are presented.
 Considering many proofs are involved, we postpone some of them in Appendix. They are interesting in their own right.

In this paper we use the following notation. For a sequence of random variables $\{U_N;\, N \geq 1\}$ and a sequence of constants $\{a_N;\, N\geq 1\}$, the notation $U_N=o_p(a_N)$ means that $U_N/a_N \to 0$ in probability as $N\to\infty;$ we write $U_N=O_p(a_N)$ if $\{U_N/a_N;\, N\geq 1\}$ is stochastically bounded, that is, $\lim_{A\to \infty}\limsup_{N\to \infty}P(|U_N/a_N|\geq A)=0.$ In particular, if $U_N=O_p(a_N)$ then $U_N=o_p(a_Nb_N)$ for any sequence of numbers $\{b_N;\, N\geq 1\}$ with $\lim_{N\to\infty}b_N=\infty.$ We write $a_N\sim b_N$ if $\lim_{N\to\infty}\frac{a_N}{b_N}=1$ for any two sequence of numbers $\{a_N;\, N\geq 1\}$ and $\{b_N;\, N\geq 1\}$.

\subsection{The proofs of Theorems \ref{diligent} and \ref{Fall_nice}}
The proof of Theorem \ref{diligent} is lengthy. The main tool is the Lindeberg-Feller central limit theorem for martingales. Automatically many computations of conditional means and variances as well as higher moments are needed for sample correlation coefficients $\hat{\rho}_{ij}$. They are non-trivial. To make the proof organized, we decide to put key steps in a few of sections. This may best facilitate the understanding of readers.

\subsubsection{Prelude 1: technical lemmas towards proofs of Theorems  \ref{diligent} and \ref{Fall_nice}}\lbl{Pre_123}

The proofs of the results in this section  will be presented in Section \ref{att_123}.

\begin{lemma}\lbl{mac_pencil}Let $\xi$ be a random variable with $E\xi=a$. Let $\tau\geq 2$ be given. The following holds.

(i) If $a=0$, then
\begin{eqnarray*}
 E[ |\xi^2-E\xi^2|^{\tau} ] \leq 32^{\tau}\cdot
E(|\xi|^{2\tau}).
\end{eqnarray*}

(ii) If $a\ne 0$, then
\begin{eqnarray*}
&& E[ |\xi^2-E\xi^2|^{\tau} ]\\
&\leq & 16^{\tau}\cdot \Big[|a|^{-\tau}\cdot {\rm Var}(\xi)^{\tau}+\sqrt{E(|\xi-a|^{2\tau})}\Big]\\
&&\cdot \Big[|a|^{\tau}+|a|^{-\tau}\cdot {\rm Var}(\xi)^{\tau}+\sqrt{E(|\xi-a|^{2\tau})}\Big].
\end{eqnarray*}
\end{lemma}

\smallskip

The following is the Marcinkiewicz-Zygmund inequality; see, e.g., p. 386 and p. 387 from \cite{CT88}.

\begin{lemma} Let $m\geq 1$ and $\{\xi_i;\, 1\leq i \leq m\}$ be independent random variables with $E\xi_i=0$ for each $i$   and $\sup_{1\leq i \leq m}E(|\xi_i|^\tau)<\infty$ for some $\tau\geq 2$. Then there exists  a constant $K_{\tau}>0$  depending on $\tau$ only such that
\begin{align}
E(|\xi_1+\cdots +\xi_m|^{\tau})\leq & K_{\tau}\cdot E\big[\big(\xi_1^2+\cdots +\xi_m^2\big)^{\tau/2}\big] \lbl{hot_sleep1}\\
 \leq & K_{\tau}\cdot m^{(\tau/2)-1}\big(E|\xi_1|^{\tau}+\cdots +E|\xi_m|^{\tau}\big). \lbl{hot_sleep2}
\end{align}
\end{lemma}

\smallskip

\subsubsection{Prelude 2: mixing moments on random variables uniformly distributed on spheres}\lbl{beauty_ugly}

In this subsection we develop some identities and inequalities regarding moments of random vectors with the uniform distribution on high-dimensional unit spheres.
The tools and methods are of independent interest. The proof of Lemma \ref{long_for} is given in this section to show the main idea and starting point.  The remaining proofs of other lemmas will be presented in Section \ref{att_proof}.

Review the setting  above \eqref{park_1} and notation $\bd{P}_i$ and $\bd{\epsilon}_i=(\epsilon_{i1},\cdots,\epsilon_{iT})^{'}\in \mathbb{R}^T$ for each $i$. The notation $\mathbb{S}^{m-1}$ represents the unit sphere in the $m$-dimensional Euclidean space.

\begin{lemma}\lbl{long_for} Set $m=T-p.$  Let $\bd{O}_i$ be  a $T\times T$ orthogonal matrix such that
\begin{align}\lbl{liberation}
\bd{P}_i=\bd{O}_i
\begin{pmatrix}
\bd{I}_m & \bd{0}\\
\bd{0} & \bd{0}
\end{pmatrix}
\bd{O}_i',\ \ \ 1\leq i \leq N.
\end{align}
Write $\bd{O}_i=(\bd{U}_{i}, \bd{V}_{i})$ for each $i$, where $\bd{U}_{i}$ is a $T\times m$ submatrix. Let $\{\epsilon_{ij};\, 1\leq i\leq N, 1\leq j \leq  T\}$ be independent random variables with $\epsilon_{ij}\sim N(0, \sigma_i^2)$, $\sigma_i>0$, for all $i$ and $j$.
 Write $\bd{\epsilon}_i=(\epsilon_{i1},\cdots,\epsilon_{iT})^{'}\in \mathbb{R}^T$ for each $i$. Let $\bd{s}_1, \cdots, \bd{s}_N$ be i.i.d. random vectors uniformly distributed on $\mathbb{S}^{m-1}.$
Then $(\frac{\bd{P}_1\bd{\epsilon}_1}{\|\bd{P}_1\bd{\epsilon}_1\|}, \cdots, \frac{\bd{P}_N\bd{\epsilon}_N}{\|\bd{P}_N\bd{\epsilon}_N\|})$ and $(\bd{U}_1\bd{s}_1, \cdots, \bd{U}_{N}\bd{s}_N)$ have the same  distribution.
\end{lemma}
\noindent\textbf{Proof of Lemma \ref{long_for}}. By the scale-invariance of $\frac{\bd{P}_i\bd{\epsilon}_i}{\|\bd{P}_i\bd{\epsilon}_i\|}$, without loss of generality, assume $\sigma_1=\cdots=\sigma_N=1$.
Evidently, $(\bd{U}_i, \bd{0})\bd{\epsilon}_i=\bd{U}_i\bd{\eta}_i$ for each $i$, where $\bd{\eta}_i=(\epsilon_{i1},\cdots,\epsilon_{im})'$. By the orthogonality and \eqref{liberation},
\begin{align}\lbl{happy_moon}
\bd{U}_i'\bd{U}_i=\bd{I}_m\ \ \mbox{and}\ \ \bd{U}_i\bd{U}_i'=\bd{P}_i.
\end{align}
By the orthogonal invariance of normal distributions and \eqref{liberation} again, $\bd{P}_i\bd{\epsilon}_i=(\bd{U}_i, \bd{0})\bd{O}_i'\bd{\epsilon}_i$ has the same distribution as that of $(\bd{U}_i, \bd{0})\bd{\epsilon}_i=\bd{U}_i\bd{\eta}_i$.  Then $\frac{\bd{P}_i\bd{\epsilon}_i}{\|\bd{P}_i\bd{\epsilon}_i\|}$, as a function of  $\bd{P}_i\bd{\epsilon}_i$, has the same distribution as that of
\begin{eqnarray*}
\frac{\bd{U}_i\bd{\eta}_i}{\|\bd{U}_i\bd{\eta}_i\|}=
\frac{\bd{U}_i\bd{\eta}_i}{(\bd{\eta}_i'\bd{U}_i'\bd{U}_i\bd{\eta}_i)^{1/2}}=\bd{U}_i\frac{\bd{\eta}_i}{\|\bd{\eta}_i\|}
\end{eqnarray*}
for each $i$ by the first identity of \eqref{happy_moon}. The desired conclusion then follows from the independence among $\{\bd{\epsilon}_1, \cdots, \bd{\epsilon}_N\}$. \hfill$\Box$

\begin{lemma}\lbl{trivial12} Let $\bd{U}_i$'s be as in Lemma \ref{long_for}. The following holds.

(i) Set $\bd{M}_{ij}=\bd{U}_i'\bd{U}_j\bd{U}_j'\bd{U}_i$ for any $1\leq i<j \leq N.$ Then both $\bd{M}_{ij}$  and $\bd{I}_m-\bd{M}_{ij}$ are non-negative definite.

(ii) Let $\bd{M}$ be a $m\times m$ non-negative definite matrix satisfying that $\bd{I}_m-\bd{M}$ is  non-negative definite. Then $\bd{I}_m-\bd{M}^2$ is also non-negative definite.
\end{lemma}

\begin{lemma}\lbl{trivial} Let $\{\bd{M}_i, i=1,2\}$ be non-negative definite matrices. Assume $\bd{M}_1$ is idempotent, that is,  $\bd{M}_1^2=\bd{M}_1$. Then $\mbox{tr}\,(\bd{M}_1\bd{M}_2) \leq  \mbox{tr}\,(\bd{M}_2)$.
\end{lemma}

\begin{lemma}\lbl{brother_cat} Let $\bd{M}_1$ and $\bd{M}_2$  be $n\times n$ non-negative definite matrices. Then, $\mbox{tr}(\bd{M}_1\bd{M}_2) \geq 0$ and  $[\mbox{tr}(\bd{M}_1\bd{M}_2)]^2 \leq r\cdot \mbox{tr}((\bd{M}_1\bd{M}_2)^2)$, where $r:=\textrm{rank}(\bd{M}_1\bd{M}_2)$ $\leq$ $n$.
\end{lemma}

Recall notation $(2m-1)!!=1\cdot3\cdots (2m-1)$ for any integer $m\geq 1$. By convention we set  $(-1)!!=1$.

\begin{lemma}\lbl{geese}[Lemma 2.4 from \cite{Jiang09}]. Suppose $m\geq 2$ and $Z_1, \cdots, Z_m$ are i.i.d. $N(0, 1)$-distributed random variables. Define $U_i=Z_i^2/(Z_1^2+\cdots +Z_m^2)$ for $1\leq i \leq m$. Let $a_1, \cdots, a_m$ be nonnegative integers. Set $a=a_1+\cdots +a_m$. Then
\begin{eqnarray*}
E(U_1^{a_1} U_2^{a_2}\cdots U_m^{a_m})=\frac{\prod_{i=1}^m(2a_i-1)!!}{\prod_{i=1}^a(m+2i-2)}.
\end{eqnarray*}
\end{lemma}

\begin{lemma}\lbl{sun_quiet} Let $m\geq 2$ and  $\{Z_i;\, 1\leq i \leq m\}$ be i.i.d. $N(0,1)$-distributed  random variables. Set $\bd{d}=(Z_1, \cdots, Z_m)'/(Z_1^2+ \cdots +Z_m^2)^{1/2}$. Let $\bd{M}$ be a symmetric matrix. Then
\begin{eqnarray*}
&(i)& E(\bd{d}'\bd{M}\bd{d})=\frac{1}{m}\cdot \mbox{tr}(\bd{M});\\
&(ii)&
E[(\bd{d}'\bd{M}\bd{d})^2]=\frac{1}{m(m+2)}\cdot \big\{2\,\mbox{tr}(\bd{M}^2)+[\mbox{tr}(\bd{M})]^2\big\};\\
&(iii)&
\mbox{Var}(\bd{d}'\bd{M}\bd{d})=\frac{2}{m(m+2)}\cdot \mbox{tr}(\bd{M}^2)-\frac{2}{m^2(m+2)}\cdot \big[\mbox{tr}(\bd{M})\big]^2.
\end{eqnarray*}
\end{lemma}

\begin{lemma}\lbl{dust_bone} Let $\{Z_i;\, 1\leq i \leq m\}$ be i.i.d. $N(0,1)$-distributed  random variables. Set $\bd{d}=(Z_1, \cdots, Z_m)'/(Z_1^2+ \cdots +Z_m^2)^{1/2}$ for $1\leq i \leq m$. Let $\bd{M}$ be a $m\times m$ symmetric matrix. Let $\tau\geq 1$ be given. Then,
\begin{eqnarray*}
E\big[|\bd{d}'\bd{M}\bd{d}-E(\bd{d}'\bd{M}\bd{d})|^{\tau}\big]\leq \frac{C_\tau}{m^{\tau}}\cdot \Big\{\mbox{tr}(\bd{M}^2)-\frac{1}{m}[\mbox{tr}(\bd{M})]^2\Big\}^{\tau/2}
\end{eqnarray*}
for all $m\geq 4\tau+1$, where $C_{\tau}>0$ is a constant depending on $\tau$ only.
\end{lemma}

\begin{lemma}\lbl{gan_geming} Let $\{Z_i;\, 1\leq i \leq m\}$ be i.i.d. $N(0,1)$-distributed  random variables. Set $\bd{d}=(Z_1, \cdots, Z_m)'/(Z_1^2+ \cdots +Z_m^2)^{1/2}$ for $1\leq i \leq m$. Let $\bd{a}\in \mathbb{R}^m$ be a vector and  $\bd{M}$ be a $m\times m$ symmetric matrix. Let $\tau\geq 1$ be given. Then, $E(|\bd{a}'\bd{d}|^{2\tau})\leq C_\tau\|\bd{a}\|^{2\tau}/m^{\tau}$ and
\begin{eqnarray*}
E\big(\bd{d}'\bd{M}\bd{d}\big)^{\tau}\leq \frac{C_\tau}{m^{\tau}}\cdot \Big\{|\mbox{tr}(\bd{M})|^{\tau} + \Big[\mbox{tr}(\bd{M}^2)-\frac{1}{m}[\mbox{tr}(\bd{M})]^2\Big]^{\tau/2}\Big\}
\end{eqnarray*}
for all $m\geq 2\tau+1$, where $C_{\tau}>0$ is a constant depending on $\tau$ only.
\end{lemma}

\begin{lemma}\lbl{R_V} Let $\{\bd{h}, \bd{h}_1, \bd{h}_2\}$ be i.i.d. $\mathbb{R}^m$-valued random vectors, where $\bd{h}$ has the same distribution as $\bd{d}$ in Lemma \ref{dust_bone}. Let $\bd{A}$, $\bd{B}$ and $\bd{C}$ be $m\times m$ matrices. Then

(i) $E\big[(\bd{h}'\bd{A}\bd{h})(\bd{h}'\bd{B}\bd{h})\big]
=\frac{1}{m(m+2)}\big[2\,\mbox{tr}(\bd{A}\bd{B})+\mbox{tr}(\bd{A})\cdot\mbox{tr}(\bd{B}) \big]$ if $\bd{A}$ and $\bd{B}$ are symmetric.

(ii) $\mbox{Var}[(\bd{h}_1'\bd{C}\bd{h}_2)^2]\leq \frac{K}{m^{5/2}}\cdot {[\mbox{tr}((\bd{C}\bd{C}')^4)]}^{1/2}$, where $K>0$ is a constant.

(iii) $\mbox{Cov}\big[(\bd{h}'\bd{A}\bd{h}_1)^2, (\bd{h}'\bd{B}\bd{h}_2)^2\big]=\frac{2}{m^3(m+2)}\cdot \mbox{tr}(\bd{A}\bd{A}'\bd{B}\bd{B}')-\frac{2}{m^4(m+2)}\,
\mbox{tr}(\bd{A}\bd{A}')\cdot\mbox{tr}(\bd{B}\bd{B}')$.
\end{lemma}

A quick reminder is that, although we assume that $\bd{A}$ and $\bd{B}$ are symmetric in (i) above, we do no  need that $\bd{A}$, $\bd{B}$ or $\bd{C}$ are symmetric in (ii) and (iii).

\begin{lemma}\lbl{use_tea} Review $\bd{P}_i$ in \eqref{park_1} and $\hat{\rho}_{ij}$ in \eqref{Asia_foods}.  Recall $\bd{U}_i$ and $\bd{s}_i$ from Lemma \ref{long_for} and  $\bd{M}_{ij}=\bd{U}_i'\bd{U}_j\bd{U}_j'\bd{U}_i$ from Lemma \ref{trivial12}.  The  following statements hold for all $i\ne j.$

(i) $E\hat{\rho}_{ij}=0$ and $E(\hat{\rho}_{ij}^2)=\frac{1}{m^2}\cdot\mbox{tr}(\bd{P}_i\bd{P}_j)$.

(ii) $E[\hat{\rho}_{ij}|\bd{s}_i]=0$ and  $E[\hat{\rho}_{ij}^2|\bd{s}_i]=\frac{1}{m}\cdot \bd{s}_i'\bd{M}_{ij}\bd{s}_i$.

\end{lemma}

\smallskip

In the following we will use notation $\mbox{Var}(\xi_2|\xi_1)$ for conditional variance, which is defined by $E(\xi_2^2|\xi_1)-[E(\xi_2|\xi_1)]^2$ for any random variables $\xi_1$ and $\xi_2.$

\begin{lemma}\lbl{good_tie} Review $\bd{P}_i$ in \eqref{park_1} and $\hat{\rho}_{ij}$ in \eqref{Asia_foods}.  Recall $\bd{U}_i$ and $\bd{s}_i$ from Lemma \ref{use_tea} and  $\bd{M}_{ij}=\bd{U}_i'\bd{U}_j\bd{U}_j'\bd{U}_i$ from Lemma \ref{trivial12}. The  following statements are true for all $i\ne j.$

(i) $E\big[(\hat{\rho}_{ij})^4\big|\bd{s}_i\big]
= \frac{3}{m(m+2)}\cdot\big(\bd{s}_i'\bd{M}_{ij}\bd{s}_i\big)^2.$

(ii) $E[(\hat{\rho}_{ij})^4]=\frac{3}{m^2(m+2)^2}\cdot \big\{2\,\mbox{tr}[(\bd{P}_i\bd{P}_j)^2]+   [\mbox{tr}(\bd{P}_i\bd{P}_j)]^2\big\}$.

(iii) $\mbox{Var}(\hat{\rho}_{ij}^2|\bd{s}_i)=\frac{2(m-1)}{m^2(m+2)}\cdot
\big(\bd{s}_i'\bd{M}_{ij}\bd{s}_i\big)^2$.

(iv) $\mbox{Var}(\hat{\rho}_{ij}^2)
 = \frac{6}{m^2(m+2)^2}\cdot \mbox{tr}[(\bd{P}_i\bd{P}_j)^2]+\frac{2(m^2-2m-2)}{m^4(m+2)^2}\cdot [\mbox{tr}(\bd{P}_i\bd{P}_j)]^2.$
\end{lemma}

\subsubsection{Intermezzo 1: calculations of variances of sums related to sample correlation coefficients}\lbl{huangjingou}

In \eqref{haoshichengshuang} and \eqref{shadongxi}, we see parameters  $p, N, T$ and variables $\bd{x}_i$.
In the rest of the paper, we will use or develop many inequalities where a constant $C$ will appear frequently. The constant $C$ does not depend on $p, N, T$ or $\bd{x}_i$'s and it can be different from line to line. The proofs of the lemma in this section  will be given in Section \ref{huangjingou_hei}.

\begin{lemma}\lbl{Love_sound} Review the notations $p, T, N$ and $\bd{P}_i$ in \eqref{park_1}. Let $p$ be fixed and $m=T-p\geq 1$. For any set $S\subset \{1,2,\cdots,N\}$ with $q=|S|\in \{1, \cdots, N-1\}$, define $\bd{P}_S=\sum_{k\in S}\bd{P}_k.$  Let  $j\notin S.$  Then there exists a constant $K>0$ depending on $p$ but not on $N$, $T$ or $\bd{P}_i$ such that the following statements hold uniformly for all $1\leq i< j \leq N$ and $N\geq 4.$
\begin{flalign}
&\ \ \ \ \ (i)\  \frac{1}{T}\cdot \big|[\mbox{tr}(\bd{P}_i\bd{P}_j)]^2-T^2\big|\leq K.& \nonumber\\
&\ \ \ \ \  (ii)\  \big|\mbox{tr}((\bd{P}_i\bd{P}_j)^2)-T\big|\leq K.&\nonumber\\
&\ \ \ \ \  (iii) \ \frac{1}{Tq^2}\cdot \big|[\mbox{tr}(\bd{P}_S\bd{P}_j)]^2-T^2q^2\big|\leq K.&\nonumber\\
&\ \ \ \ \  (iv)\  \frac{1}{q^2}\cdot \big|\mbox{tr}((\bd{P}_S\bd{P}_j)^2)-Tq^2\big|\leq K.&\nonumber\\
&\ \ \ \ \  (v)\  \mbox{Statements (i)-(iv) still hold if symbol}\ ``T"\, \mbox{is replaced by }\, ``m".&\nonumber
\end{flalign}
\end{lemma}

\begin{lemma}\lbl{shadow_1} Recall $\bd{U}_i$ from Lemma \ref{long_for} and  $\bd{M}_{ij}=\bd{U}_i'\bd{U}_j\bd{U}_j'\bd{U}_i$ from Lemma \ref{trivial12}. Let $\bd{e}$ have the uniform distribution on $\mathbb{S}^{m-1}$. Then there is a constant $C>0$ free of $N, T$ and $p$ such that
$\sup_{1\leq i<j \leq N}\mbox{Var}\big((\bd{e}'\bd{M}_{ij}\bd{e})^2\big)\leq Cm^{-2}$
as $N\geq C$.
\end{lemma}

\begin{lemma}\lbl{snow_sound} Review $\bd{P}_i$ in \eqref{park_1} and $\hat{\rho}_{ij}$ in \eqref{Asia_foods}.  Recall $\bd{U}_i$ and $\bd{s}_i$ from Lemma \ref{long_for} and  $\bd{M}_{ij}=\bd{U}_i'\bd{U}_j\bd{U}_j'\bd{U}_i$ from Lemma \ref{trivial12}. Define
$\bd{P}_{j\blacktriangle}=\sum_{i=1}^{j-1}\bd{P}_i$ for $2\leq j\leq N.$ By Lemma \ref{use_tea},
\begin{eqnarray*}
X_j:=\sum_{i=1}^{j-1}[T\hat{\rho}_{ij}^2-E(T\hat{\rho}_{ij}^2|\bd{s}_i)]=\sum_{i=1}^{j-1}T\hat{\rho}_{ij}^2-\frac{T}{m}\sum_{i=1}^{j-1}
\bd{s}_i'\bd{M}_{ij}\bd{s}_i
\end{eqnarray*}
for $2\leq j\leq N.$ Then,
\begin{eqnarray*}
 \frac{1}{T^2}E(X_j^2)&=&\frac{2m-8}{m^3(m+2)^2}\sum_{i=1}^{j-1}\mbox{tr}((\bd{P}_i\bd{P}_j)^2)+ \frac{2m^2+4}{m^4(m+2)^2}\sum_{i=1}^{j-1}[\mbox{tr}(\bd{P}_i\bd{P}_j)]^2+\\
&&\frac{2}{m^3(m+2)}
\Big\{\mbox{tr}\big((\bd{P}_{j\blacktriangle}\bd{P}_{j})^2\big)-
\frac{1}{m}\cdot\big[\mbox{tr}(\bd{P}_{j\blacktriangle}\bd{P}_{j})\big]^2\Big\}.
\end{eqnarray*}
\end{lemma}
A quick comment is that the last term above is non-negative by Lemma \ref{brother_cat}.

\begin{lemma}\lbl{kiss_rain} Let $X_j$ be defined as in Lemma \ref{snow_sound} for $2\leq j\leq N.$
Assume $p$ is fixed and  $N=o(T^2)$ as $N\to\infty$.  Then
\begin{eqnarray*}
\lim_{N\to \infty}\frac{1}{N^2}\sum_{j=2}^NE(X_j^2) = 1.
\end{eqnarray*}
\end{lemma}

\begin{lemma}\lbl{cousin_basket}  Recall $\bd{U}_i$ and $\bd{s}_i$ from Lemma \ref{long_for} and  $\bd{M}_{ij}=\bd{U}_i'\bd{U}_j\bd{U}_j'\bd{U}_i$ from Lemma \ref{trivial12}. Assume $p$ is fixed and $T=T_N\to \infty$.  Then
\begin{eqnarray*}
{\rm Var}\Big[\sum_{j=2}^N\Big(\sum_{i=1}^{j-1}\bd{s}_i'\bd{M}_{ij}\bd{s}_i\Big)^2\Big]
=O\Big(\frac{N^5}{T^2}\Big).
\end{eqnarray*}
In particular, the variance above is of order $o(N^4T^2)$ if $N=o(T^4)$ as $N\to\infty$. 
\end{lemma}

\begin{lemma}\lbl{Lose_Ren}  Recall $\bd{U}_i$ and $\bd{s}_i$ from Lemma \ref{long_for}.  
If $p$ is fixed and $T=T_N\to \infty$, then
\begin{eqnarray*}
{\rm Var}\Big\{\sum_{j=2}^N{\rm tr}\Big[\Big(\sum_{i=1}^{j-1}\bd{U}_j'\bd{U}_i\bd{s}_i
\bd{s}_i'\bd{U}_i'\bd{U}_j\Big)^2\Big]\Big\}=O\Big(\frac{N^4}{T^2}+\frac{N^5}{T^3}\Big)
\end{eqnarray*}
as $N\to\infty$. In particular, the variance is of the order $o(N^4)$ if $N=o(T^3)$. \end{lemma}

\subsubsection{Intermezzo 2: preliminary verifications of the Lindeberg-Feller condition towards proof of  Theorem \ref{diligent}}
\label{PS4}


\begin{lemma}\lbl{leaves_potato} Recall $\bd{U}_i$ and $\bd{s}_i$ from Lemma \ref{long_for} and  $\bd{M}_{ij}=\bd{U}_i'\bd{U}_j\bd{U}_j'\bd{U}_i$ from Lemma \ref{trivial12}. Let $\mu_N$ be as in \eqref{taiyangtaila}. Set $B_j=\frac{T}{m}\sum_{i=1}^{j-1}
\bd{s}_i'\bd{M}_{ij}\bd{s}_i$ for $2\leq j \leq N$. If  $p$ is fixed and $N=o(T^2)$ as $N\to \infty$, then
\begin{eqnarray*}
\frac{1}{N}\Big[ \Big(\sum_{j=2}^NB_j\Big)-\mu_N\Big] \to 0
\end{eqnarray*}
in probability as $N\to\infty$.
\end{lemma}
\noindent\textbf{Proof of Lemma \ref{leaves_potato}}.  By \eqref{happy_moon},  $\bd{U}_{i}\bd{U}_{i}'=\bd{P}_i$. Then
\begin{align}\lbl{thoron_cur}
\mbox{tr}(\bd{M}_{ij})=\mbox{tr}(\bd{U}_i\bd{U}_i'\bd{U}_j\bd{U}_j')
=\mbox{tr}(\bd{P}_i\bd{P}_j).
\end{align}
It follows from Lemma \ref{sun_quiet}(i) that
\begin{eqnarray*}
EB_j =\frac{T}{m}\sum_{i=1}^{j-1}
E(\bd{s}_i'\bd{M}_{ij}\bd{s}_i)
=\frac{T}{m^2}\sum_{i=1}^{j-1}\mbox{tr}(\bd{P}_i\bd{P}_j)
\end{eqnarray*}
for $2\leq j \leq N.$ Therefore,
\begin{eqnarray*}
\mu_N=E\sum_{j=2}^NB_j&=&\frac{T}{m^2}\sum_{j=2}^N\sum_{i=1}^{j-1}\mbox{tr}(\bd{P}_i\bd{P}_j)\\
& = & \frac{T}{2m^2}\cdot \mbox{tr}\Big[\Big(\sum_{i=1}^N\bd{P}_i\Big)^2-\sum_{i=1}^N\bd{P}_i\Big]\\
& = & \frac{T}{2m^2}\cdot \mbox{tr}\Big[\Big(\sum_{i=1}^N\bd{P}_i\Big)^2\Big]-\frac{TN}{2m}
\end{eqnarray*}
by the fact that $\mbox{tr}(\bd{P}_i)=T-p=m$ for each $i.$ On the other hand,
\begin{eqnarray*}
\sum_{j=2}^NB_j
&=&\frac{T}{m}\sum_{j=2}^N\sum_{i=1}^{j-1}
\bd{s}_i'\bd{M}_{ij}\bd{s}_i\\
& = & \frac{T}{m}\sum_{i=1}^{N-1}\sum_{j=i+1}^N
\bd{s}_i'\bd{M}_{ij}\bd{s}_i\\
& = & \frac{T}{m}\sum_{i=1}^{N-1}\bd{s}_i'\bd{M}_{i}\bd{s}_i
\end{eqnarray*}
where $\bd{M}_{i\blacktriangledown}:=\sum_{j=i+1}^N\bd{M}_{ij}.$ By independence among $\{\bd{s}_i\}$'s and Lemma \ref{sun_quiet},
\begin{align}\lbl{nest_fly}
\mbox{Var}\Big(\sum_{j=2}^NB_j\Big)
=&\frac{T^2}{m^2}\sum_{i=1}^{N-1}
\mbox{Var}\big(\bd{s}_i'\bd{M}_{i\blacktriangledown}\bd{s}_i\big)\nonumber\\
=&
\frac{T^2}{m^2}\sum_{i=1}^{N-1}\Big\{\frac{2}{m(m+2)}\cdot \mbox{tr}(\bd{M}_{i\blacktriangledown}^2)-\frac{2}{m^2(m+2)}\cdot \big[\mbox{tr}(\bd{M}_{i\blacktriangledown})\big]^2\Big\}\nonumber\\
=& \frac{2T^2}{m^3(m+2)}\sum_{i=1}^{N-1}\Big\{ \mbox{tr}(\bd{M}_{i\blacktriangledown}^2)-\frac{1}{m}\cdot \big[\mbox{tr}(\bd{M}_{i\blacktriangledown})\big]^2\Big\}.
\end{align}
From Lemma \ref{trivial12}, we know $\bd{I}_m-\bd{M}_{ij}$ is non-negative for each $i \ne j.$ Since the sum of non-negative definite matrices is still non-negative definite, we see that $(N-i)\bd{I}_m-\bd{M}_{i\blacktriangledown}$ is also non-negative definite. By Lemma \ref{trivial12}(ii),  $(N-i)^2\bd{I}_m-\bd{M}_{i\blacktriangledown}^2$ is non-negative definite. In particular,
\begin{align}\lbl{one_flower}
\mbox{tr}(\bd{M}_{i\blacktriangledown}^2)\leq (N-i)^2m.
\end{align}
Moreover,
\begin{align}\lbl{han_pro}
\mbox{tr}(\bd{M}_{i\blacktriangledown})=\sum_{j=i+1}^N\mbox{tr}(\bd{M}_{ij})
=\sum_{j=i+1}^N\mbox{tr}(\bd{P}_i\bd{P}_j)
\end{align}
by \eqref{thoron_cur}.  Now we estimate $\mbox{tr}(\bd{P}_i\bd{P}_j)$.

Recall \eqref{park_1}. Set $\bd{A}_i=\bd{x}_i(\bd{x}_i'\bd{x}_i)^{-1}\bd{x}_i'$ for $1\leq i \leq N.$ Then $\bd{A}_i$ is a $T\times T$ idempotent matrix with rank $p$ and $\mbox{tr}(\bd{A}_i)=p$ for each $i$. Since  $\bd{P}_i=\bd{I}_T-\bd{A}_i$, we see
\begin{align*}
\bd{P}_i\bd{P}_j=\bd{I}_T+\bd{B}_{ij}
\end{align*}
where $\bd{B}_{ij}:=\bd{A}_i\bd{A}_j-\bd{A}_i-\bd{A}_j.$ By Lemma \ref{brother_cat}, %
%
\begin{align*}
\mbox{tr}(\bd{F}_1\bd{F}_2)\geq 0
\end{align*}
for  any non-negative definite matrices $\bd{F}_1$ and $\bd{F}_2$. As a result, $\mbox{tr}(\bd{A}_i\bd{A}_j) \geq 0$.
Easily, $\mbox{tr}(\bd{A}_i\bd{A}_j)\leq p$ by Lemma \ref{trivial}. Thus,
\begin{align*}
-2p\leq \mbox{tr}(\bd{B}_{ij})\leq -p.
\end{align*}
Therefore,  we have $\mbox{tr}(\bd{P}_i\bd{P}_j) \geq  T-2p$. Hence,  $\mbox{tr}(\bd{M}_{i\blacktriangledown})\geq (N-i)(T-2p)$ by \eqref{han_pro}. This and \eqref{one_flower} tell us that
\begin{eqnarray*}
& & \mbox{tr}(\bd{M}_{i\blacktriangledown}^2)-\frac{1}{m}\cdot \big[\mbox{tr}(\bd{M}_{i\blacktriangledown})\big]^2\\
& \leq & (N-i)^2m-\frac{1}{m}(N-i)^2(T-2p)^2\\
& = & (N-i)^2\cdot \frac{m^2-(m-p)^2}{m}\\
& \leq & 2(N-i)^2p
\end{eqnarray*}
by recalling the notation $m=T-p.$ Plugging this into \eqref{nest_fly} we get
\begin{eqnarray*}
\mbox{Var}\Big(\sum_{j=2}^NB_j\Big)
\leq \frac{2T^2}{m^3(m+2)}\sum_{i=1}^{N-1}2N^2p\leq \frac{(4p)T^2N^3}{m^4}.
\end{eqnarray*}
By the Chebyshev inequality, for any $\tau>0$,
\begin{eqnarray*}
 P\Big(\frac{1}{2N}\Big| \Big(\sum_{j=2}^NB_j\Big)-\mu_N\Big| \geq \tau\Big)
& \leq & \frac{1}{4\tau^2N^2}\cdot \mbox{Var}\Big(\sum_{j=2}^NB_j\Big)\\
& \leq & \frac{p}{\tau^2}\cdot \frac{T^2N}{(T-p)^4},
\end{eqnarray*}
which goes to zero provided $N=o(T^2).$ \hfill$\square$

\smallskip

Let $\{\bd{s}_1, \cdots, \bd{s}_j\}$ for $1\leq j \leq N$ be defined in Lemma \ref{long_for}, which
are i.i.d. random vectors uniformly distributed on $\mathbb{S}^{m-1}$. Set
\begin{align}\lbl{cheery_coffee}
\ml{F}_0=\{\emptyset, \Omega\}\ \ \mbox{and}\ \ \ml{F}_j=\sigma(\bd{s}_1, \cdots, \bd{s}_j)
\end{align}
which is the $\sigma$-algebra generated by $\{\bd{s}_1, \cdots, \bd{s}_j\}$ for $1\leq j \leq N$. Here  $\Omega$ is the sample space on which random variables $\{\epsilon_{ij}\}$ are defined on.

\begin{lemma}\lbl{ant_egg} Let $X_j$ be defined as in Lemma \ref{snow_sound} and $\ml{F}_j$ be as in \eqref{cheery_coffee}.
Assume $N=o(T^3)$. Define
\begin{eqnarray*}
 Z_N=\frac{1}{N^2} \sum_{j=2}^NE [X_j^2|\mathcal{F}_{j-1}].
\end{eqnarray*}
Then ${\rm Var}(Z_N)\rightarrow 0$ as $N\to \infty.$
\end{lemma}
\noindent\textbf{Proof of Lemma \ref{ant_egg}}. Set $\bd{H}_{ij}=\bd{U}_i'\bd{U}_j\bd{s}_j
\bd{s}_j'\bd{U}_j'\bd{U}_i$ for $1\leq i, j \leq N$, where $\bd{U}_i$'s and $\bd{s}_j$'s are defined as in Lemma \ref{long_for}. Then $\bd{s}_i'\bd{H}_{ij}\bd{s}_i=\bd{s}_j'\bd{C}_{ij}\bd{s}_j$, where
\begin{eqnarray*}
\bd{C}_{ij}:=\bd{U}_j'\bd{U}_i\bd{s}_i\bd{s}_i'\bd{U}_i'\bd{U}_j.
\end{eqnarray*}
Since $\bd{s}_i'\bd{U}_i'\bd{U}_j\bd{s}_j=(\bd{s}_i'\bd{U}_i'\bd{U}_j\bd{s}_j)'
=\bd{s}_j'\bd{U}_j'\bd{U}_i\bd{s}_i\in \mathbb{R}$, we have
\begin{align}\lbl{army_plane}
\hat{\rho}_{ij}^2=\bd{s}_j'(\bd{U}_j'\bd{U}_i\bd{s}_i
\bd{s}_i'\bd{U}_i'\bd{U}_j)\bd{s}_j=\bd{s}_j'\bd{C}_{ij}\bd{s}_j.
\end{align}
By Lemma \ref{sun_quiet}(i) and the independence between $\bd{s}_i$ and $\bd{s}_j$, we have that
\beaa
E(\hat{\rho}_{ij}^2|\bd{s}_i)=\frac{1}{m}\mbox{tr}(\bd{C}_{ij})
=\frac{1}{m}\bd{s}_i'\bd{U}_i'\bd{U}_j\bd{U}_j'\bd{U}_i\bd{s}_i
=\frac{1}{m}\bd{s}_i'\bd{M}_{ij}\bd{s}_i
\eeaa
for $i<j$, where $\bd{M}_{ij}=\bd{U}_i'\bd{U}_j\bd{U}_j'\bd{U}_i$. Then
\begin{align}
\frac{1}{T}X_j= &\sum_{i=1}^{j-1}[\hat{\rho}_{ij}^2-E(\hat{\rho}_{ij}^2|\bd{s}_i)]\nonumber\\
=& \sum_{i=1}^{j-1}\big[\bd{s}_j'\bd{C}_{ij}\bd{s}_j
-    \frac{1}{m}
\bd{s}_i'\bd{M}_{ij}\bd{s}_i
\big] \nonumber\\
 = & \bd{s}_j'\bd{D}_{j}\bd{s}_j-\bd{W}_{j} \lbl{big_small}
\end{align}
for $2\leq j \leq N$, where
\begin{align}\lbl{BJ_we}
 \bd{D}_{j}:=\sum_{i=1}^{j-1}\bd{C}_{ij}\ \ \ \ \mbox{and}\ \ \ \ \bd{W}_{j}:=\frac{1}{m}\sum_{i=1}^{j-1}
\bd{s}_i'\bd{M}_{ij}\bd{s}_i.
\end{align}
In view of the independence among $\bd{s}_i$'s, it is easy to check from Lemma \ref{sun_quiet} that
\begin{eqnarray*}
E\big(\bd{s}_j'\bd{D}_{j}\bd{s}_j\big|\mathcal{F}_{j-1}\big)
=\frac{1}{m}\,\mbox{tr}(\bd{D}_{j})
=\frac{1}{m}\sum_{i=1}^{j-1}\mbox{tr}(\bd{C}_{ij}).
\end{eqnarray*}
Since $\mbox{tr}(\bd{C}_{ij})=\mbox{tr}(\bd{U}_j'\bd{U}_i\bd{s}_i
\bd{s}_i'\bd{U}_i'\bd{U}_j)=\bd{s}_i'\bd{M}_{ij}\bd{s}_i$, we have  that
\begin{eqnarray*}
\bd{W}_{j}=E\big(\bd{s}_j'\bd{D}_{j}\bd{s}_j\big|\mathcal{F}_{j-1}\big).
\end{eqnarray*}
This, \eqref{big_small} and Lemma \ref{sun_quiet} imply
\begin{align}\lbl{chipoli}
\frac{1}{T^2}E [X_j^2|\mathcal{F}_{j-1}]
=&
\mbox{Var}\big(\bd{s}_j'\bd{D}_{j}\bd{s}_j\Big|\mathcal{F}_{j-1}\big) \nonumber\\
 = & \frac{2}{m(m+2)}\cdot \mbox{tr}(\bd{D}_{j}^2)-\frac{2}{m^2(m+2)}\cdot \big[\mbox{tr}(\bd{D}_{j})\big]^2.
\end{align}
 From \eqref{BJ_we},
 \begin{eqnarray*}
&& \mbox{tr}(\bd{D}_{j})=
\sum_{i=1}^{j-1}\mbox{tr}(\bd{C}_{ij})=
 \sum_{i=1}^{j-1}\bd{s}_i'\bd{M}_{ij}\bd{s}_i;\\
 && \mbox{tr}(\bd{D}_{j}^2)=\mbox{tr}\Big[\big(\sum_{i=1}^{j-1}\bd{C}_{ij}\Big)^2\Big]
 =\mbox{tr}\Big[\Big(\sum_{i=1}^{j-1}\bd{U}_j'\bd{U}_i\bd{s}_i
\bd{s}_i'\bd{U}_i'\bd{U}_j\Big)^2\Big]
 \end{eqnarray*}
 by the definition of $\bd{C}_{ij}$.
Thus, we conclude from \eqref{chipoli} that
\begin{eqnarray*}
\frac{1}{T^2}E [X_j^2|\mathcal{F}_{j-1}]
&=&
\frac{2}{m(m+2)}\cdot\mbox{tr}\Big[\Big(\sum_{i=1}^{j-1}\bd{U}_j'\bd{U}_i\bd{s}_i
\bd{s}_i'\bd{U}_i'\bd{U}_j\Big)^2\Big]-\\
&& \frac{2}{m^2(m+2)}\cdot \Big(\sum_{i=1}^{j-1}\bd{s}_i'\bd{M}_{ij}\bd{s}_i\Big)^2.
\end{eqnarray*}
It follows that
\begin{eqnarray*}
\frac{N^2}{T^2} Z_N &=& \frac{1}{T^2}\sum_{j=2}^NE [X_j^2|\mathcal{F}_{j-1}]\\
&= & \frac{2}{m(m+2)}\cdot\sum_{j=2}^N\mbox{tr}\Big[\Big(\sum_{i=1}^{j-1}\bd{U}_j'\bd{U}_i\bd{s}_i
\bd{s}_i'\bd{U}_i'\bd{U}_j\Big)^2\Big]-\\
&& \frac{2}{m^2(m+2)}\cdot \sum_{j=2}^N\Big(\sum_{i=1}^{j-1}\bd{s}_i'\bd{M}_{ij}\bd{s}_i\Big)^2.
\end{eqnarray*}
Review $T=m+p$. Since ${\rm Var}(\xi_1+\xi_2)\leq 2{\rm Var}(\xi_1)+2{\rm Var}(\xi_2)$ for any random variables $\xi_1$ and $\xi_2$, to show ${\rm Var}(Z_N)\rightarrow 0$, it is enough to prove the following two facts.
\begin{align}
&  {\rm Var}\Big\{\sum_{j=2}^N\mbox{tr}\Big[\Big(\sum_{i=1}^{j-1}\bd{U}_j'\bd{U}_i\bd{s}_i
\bd{s}_i'\bd{U}_i'\bd{U}_j\Big)^2\Big]\Big\}=o(N^4);\lbl{tortoise1}\\
& {\rm Var}\Big[\sum_{j=2}^N\Big(\sum_{i=1}^{j-1}\bd{s}_i'\bd{M}_{ij}\bd{s}_i\Big)^2\Big]
=o(N^4T^2).\lbl{tortoise2}
\end{align}
Under restriction $N=o(T^3)$, the assertion \eqref{tortoise1} is confirmed in Lemma \ref{Lose_Ren} and
\eqref{tortoise2} is proved in Lemma \ref{cousin_basket}.  The proof is completed. \hfill$\square$

\begin{lemma}\lbl{lost_found} Let $X_j$ be defined as in Lemma \ref{snow_sound} and $\ml{F}_j$ be as in \eqref{cheery_coffee}. Assume $N=o(T^4)$.  Then
\begin{eqnarray*}
\frac{1}{N^4} \sum_{j=2}^N E(X_j^4|\ml{F}_{j-1})\to 0
\end{eqnarray*}
in probability  as $N\to \infty.$
\end{lemma}
\noindent\textbf{Proof of Lemma \ref{lost_found}}. It suffices to show
\begin{align}\lbl{day_cough}
\frac{1}{N^4} \sum_{j=2}^N E(X_j^4)\to 0
\end{align}
as $N\to\infty.$ By \eqref{big_small} and \eqref{BJ_we},
\begin{align}\lbl{ice_mountain}
\frac{1}{T}X_j
=  \bd{s}_j'\bd{D}_{j}\bd{s}_j-\frac{1}{m}\sum_{i=1}^{j-1}
\bd{s}_i'\bd{M}_{ij}\bd{s}_i
\end{align}
for $2\leq j \leq N$, where $\bd{M}_{ij}=\bd{U}_i'\bd{U}_j\bd{U}_j'\bd{U}_i$ and
\begin{eqnarray*}
 \bd{D}_{j}=\sum_{i=1}^{j-1}\bd{C}_{ij}\ \ \ \mbox{and}\ \ \
 \bd{C}_{ij}=\bd{U}_j'\bd{U}_i\bd{s}_i\bd{s}_i'\bd{U}_i'\bd{U}_j.
\end{eqnarray*}
 Notice
\begin{eqnarray*}
\bd{s}_j'\bd{D}_{j}\bd{s}_j
&=&\sum_{i=1}^{j-1}\bd{s}_j'\bd{U}_j'\bd{U}_i\bd{s}_i\bd{s}_i'\bd{U}_i'\bd{U}_j\bd{s}_j\\
&=& \sum_{i=1}^{j-1}\bd{s}_i'\bd{H}_{ij}\bd{s}_i
\end{eqnarray*}
where
\begin{eqnarray*}
\bd{H}_{ij}:=\bd{U}_i'\bd{U}_j\bd{s}_j\bd{s}_j'\bd{U}_j'\bd{U}_i
\end{eqnarray*}
for $1\leq i<j\leq N.$ By Lemma \ref{sun_quiet},
\begin{align}
& \mu_{ij}:=E(\bd{s}_i'\bd{H}_{ij}\bd{s}_i|\bd{s}_j)=\frac{1}{m}\mbox{tr}(\bd{H}_{ij})
=\frac{1}{m}\bd{s}_j'\bd{M}_{ji}\bd{s}_j; \nonumber\\
& \nu_{ij}:=E(\bd{s}_j'\bd{M}_{ji}\bd{s}_j)=\frac{1}{m}\mbox{tr}(\bd{M}_{ji})
=\frac{1}{m}\mbox{tr}(\bd{P}_i\bd{P}_j) \lbl{shuo_luan}
\end{align}
for any $1\leq i<j\leq N.$ We rewrite \eqref{ice_mountain} to have
\begin{eqnarray*}
\frac{1}{T}X_j
&=&  \sum_{i=1}^{j-1}(\bd{s}_i'\bd{H}_{ij}\bd{s}_i-\mu_{ij})-
\frac{1}{m}\sum_{i=1}^{j-1}(\bd{s}_i'\bd{M}_{ij}\bd{s}_i-\nu_{ij})+ \sum_{i=1}^{j-1}(\mu_{ij}-\frac{1}{m}\nu_{ij})\\
&:=& A_j +B_j +C_j.
\end{eqnarray*}
Therefore,
\begin{align}\lbl{liquor_nice}
\frac{1}{T^4}E(X_j^4) \leq 3^3\cdot \big[E(|A_j|^4) +E(|B_j|^4) +E(|C_j|^4)\big].
\end{align}
Note that $A_j$ is the sum of independent random variables. By \eqref{hot_sleep2} with $\tau=4$,
\begin{align}
E\big(|A_j|^4|\bd{s}_j\big) \leq & C\cdot (j-1)\cdot \sum_{i=1}^{j-1}E\big[(\bd{s}_i'\bd{H}_{ij}\bd{s}_i-\mu_{ij})^4|\bd{s}_j\big] \nonumber\\
 \leq & C\cdot \frac{j}{m^4}\cdot \sum_{i=1}^{j-1}\Big\{\mbox{tr}(\bd{H}_{ij}^2)-\frac{1}{m}[\mbox{tr}(\bd{H}_{ij})]^2\Big\}^{2} \lbl{who_water}\\
 \leq & C\cdot \frac{j}{m^4}\cdot \sum_{i=1}^{j-1}\big[\mbox{tr}(\bd{H}_{ij}^2)\big]^2, \nonumber
\end{align}
where the second inequality follows from Lemma \ref{dust_bone}, and where the fact that $\mbox{tr}(\bd{H}_{ij}^2)\geq \frac{1}{m}[\mbox{tr}(\bd{H}_{ij})]^2$ from Lemma \ref{brother_cat} is used in the last step. Easily, $\mbox{tr}(\bd{H}_{ij}^2)=(\bd{s}_j'\bd{M}_{ji}\bd{s}_j)^2.$
Take another expectation to see
\begin{align}\lbl{roll_over}
E(|A_j|^4) \leq  C\cdot \frac{j}{m^4}\cdot \sum_{i=1}^{j-1}E[(\bd{s}_j'\bd{M}_{ji}\bd{s}_j)^4].
\end{align}
By Lemma \ref{gan_geming} with $\tau=4$ and the fact that $\mbox{tr}(\bd{M}^2)\geq \frac{1}{m}[\mbox{tr}(\bd{M})]^2$ for any $m\times m$ symmetric matrix $\bd{M}$ from Lemma \ref{brother_cat}, we obtain
\begin{eqnarray*}
E[(\bd{s}_j'\bd{M}_{ji}\bd{s}_j)^4]
& \leq & \frac{C}{m^{4}}\cdot \Big\{\mbox{tr}(\bd{M}_{ji}^2)-\frac{1}{m}\big[\mbox{tr}(\bd{M}_{ji})]^2\Big\}^2+ \frac{C}{m^4}\cdot\big[\mbox{tr}(\bd{M}_{ji})\big]^{4}\\
& \leq & \frac{C}{m^{4}}\cdot \big[\mbox{tr}(\bd{M}_{ji}^2)\big]^{2}+ \frac{C}{m^4}\cdot\big[\mbox{tr}(\bd{M}_{ji})\big]^{4}.
\end{eqnarray*}
It is used before that $\mbox{tr}(\bd{M}_{ji})=\mbox{tr}(\bd{P}_i\bd{P}_j)$ and $\mbox{tr}[\bd{M}_{ji}^2]=\mbox{tr}[(\bd{P}_i\bd{P}_j)^2]$. By Lemma \ref{trivial},  both quantities are bounded by $m.$ Hence,  $E\big[(\bd{s}'\bd{M}_{ji}\bd{s})^4\big]\leq C$ uniformly for all  $1\leq i<j\leq N$. We conclude from \eqref{roll_over} that
\begin{align}\lbl{JJJ}
E(|A_j|^4) \leq C\cdot \frac{j^2}{m^4}
\end{align}
uniformly for all $2\leq j \leq N$.

Now we estimate $B_j$. Replace ``$\bd{H}_{ij}$" in \eqref{who_water} with  ``$\bd{M}_{ij}$" to see
\begin{align}
E(|B_j|^4) \leq &  C\cdot\frac{1}{m^4}\cdot \frac{j}{m^4}\cdot \sum_{i=1}^{j-1}\Big\{\mbox{tr}(\bd{M}_{ij}^2)-
\frac{1}{m}[\mbox{tr}(\bd{M}_{ij})]^2\Big\}^{2}\nonumber\\
 = & C\cdot \frac{j}{m^8}\cdot \sum_{i=1}^{j-1}\Big[\mbox{tr}((\bd{P}_j\bd{P}_i)^2)
-\frac{1}{m}(\mbox{tr}(\bd{P}_j\bd{P}_i))^2\Big]^{2}\nonumber\\
 \leq & C\cdot \frac{j^2}{m^8}, \lbl{wind_cycle}
\end{align}
where the last step holds by (i) and (ii) from Lemma \ref{Love_sound}.

Finally, by \eqref{shuo_luan},
\begin{eqnarray*}
C_j&=&\frac{1}{m}\sum_{i=1}^{j-1}\big[\bd{s}_j'\bd{M}_{ji}\bd{s}_j
-\frac{1}{m}\mbox{tr}(\bd{P}_i\bd{P}_j)\big]\\
&=& \frac{1}{m}\big[\bd{s}_j'\bd{M}_{j\blacktriangle}\bd{s}_j
-E(\bd{s}_j'\bd{M}_{j\blacktriangle}\bd{s}_j)\big],
\end{eqnarray*}
where
\begin{eqnarray*}
\bd{M}_{j\blacktriangle}:=\sum_{i=1}^{j-1}\bd{M}_{ji}.
\end{eqnarray*}
 Since $\mbox{tr}(\bd{M}_{ji})=\mbox{tr}(\bd{P}_i\bd{P}_j)$, by defining
\begin{eqnarray*}
\bd{P}_{j\blacktriangle}:=\sum_{i=1}^{j-1}\bd{P}_i,
\end{eqnarray*}
we have
$\mbox{tr}(\bd{M}_{j\blacktriangle})=\mbox{tr}(\bd{P}_{j\blacktriangle}\bd{P}_j)$. Recall  \eqref{happy_moon},  $\bd{U}_{i}\bd{U}_{i}'=\bd{P}_i$. Easily,
\begin{align*}
\mbox{tr}(\bd{M}_{ji}\bd{M}_{jk})=\mbox{tr}(\bd{P}_{i}\bd{P}_{j}\bd{P}_{k}\bd{P}_{j})
\end{align*}
 for any $1\leq i, j, k\leq N$. It follows that
 \begin{align*}
\mbox{tr}(\bd{M}_{j\blacktriangle}^2)=\mbox{tr}\Big[\Big(\sum_{i=1}^{j-1}\bd{M}_{ji}\Big)^2\Big]
=\sum_{1\leq i, k \leq j-1}\mbox{tr}(\bd{P}_{i}\bd{P}_{j}\bd{P}_{k}\bd{P}_{j})
=  \mbox{tr}\big(\big(\bd{P}_{j\blacktriangle}\bd{P}_{j}\big)^2\big).
\end{align*}
On the other hand, recall $m=T-p.$ By Lemma \ref{Love_sound}(v), there exists a constant $K$ not depending on $T$ or $N$ such that
\begin{eqnarray*}
&& \frac{1}{mj^2}\cdot \big|[\mbox{tr}(\bd{P}_{j\blacktriangle}\bd{P}_j)]^2-m^2(j-1)^2\big|\leq K,\\
&& \frac{1}{j^2}\cdot \big|\mbox{tr}((\bd{P}_{j\blacktriangle}\bd{P}_j)^2)-m(j-1)^2\big|\leq K
\end{eqnarray*}
for every $2\leq j \leq N.$ It follows from the triangle inequality that
\begin{align*}
\Big|\mbox{tr}\big((\bd{P}_{j\blacktriangle}\bd{P}_{j})^2\big)-
\frac{1}{m}\cdot\big[\mbox{tr}(\bd{P}_{j\blacktriangle}\bd{P}_{j})\big]^2\Big|\leq 2Kj^2
\end{align*}
for $2\leq j \leq N.$ Consequently, by taking $\tau=4$ in  Lemma \ref{dust_bone} we have that
\begin{eqnarray*}
E(|C_j|^4) &\leq & \frac{1}{m^4}\cdot \frac{C}{m^{4}}\cdot \Big\{\mbox{tr}[\bd{M}_{j\blacktriangle}^2]-
\frac{1}{m}(\mbox{tr}(\bd{M}_{j\blacktriangle}))^2\Big\}^2\\
& = &\frac{C}{m^{8}}\cdot \Big\{\mbox{tr}[\big(\bd{P}_{j\blacktriangle}\bd{P}_{j}\big)^2]-
\frac{1}{m}[\mbox{tr}(\bd{P}_{j\blacktriangle}\bd{P}_j)]^2\Big\}^2\\
& \leq & C\cdot \frac{j^4}{m^8}
\end{eqnarray*}
uniformly for all $2\leq j \leq N$.
Combining this with \eqref{liquor_nice}, \eqref{JJJ} and \eqref{wind_cycle}, we arrive at
\begin{eqnarray*}
E(X_j^4) &\leq & CT^4\Big(\frac{j^2}{m^4}+ \frac{j^2}{m^8}+\frac{j^4}{m^8}\Big)\\
&\leq & C\cdot \Big(j^2+\frac{j^4}{m^4}\Big)
\end{eqnarray*}
uniformly for all $2\leq j \leq N$ as $N$ is large (reviewing $m=T-p$ and $T=T_N\to \infty).$ As a result,
\begin{eqnarray*}
\frac{1}{N^4} \sum_{j=2}^N E(X_j^4)=O\Big(\frac{1}{N}+\frac{N}{m^4}\Big)\to 0
\end{eqnarray*}
as $N\to \infty$ as long as $N=o(T^4).$ We obtain \eqref{day_cough}. \hfill$\Box$

\begin{lemma}\lbl{CLT_step} Let $X_j$ be defined as in Lemma \ref{snow_sound}.  Assume $N=o(T^2)$ as $N\to \infty.$ Then $\frac{1}{N}\sum_{j=2}^NX_j\to N(0, 1)$ in distribution as $N\to \infty.$
\end{lemma}
\noindent\textbf{The Proof of Lemma \ref{CLT_step}}.
Reviewing  Lemma \ref{snow_sound}, we know
\begin{eqnarray*}
X_j=\sum_{i=1}^{j-1}[T\hat{\rho}_{ij}^2-E(T\hat{\rho}_{ij}^2|\bd{s}_i)]=\sum_{i=1}^{j-1}T\hat{\rho}_{ij}^2-\frac{T}{m}\sum_{i=1}^{j-1}
\bd{s}_i'\bd{M}_{ij}\bd{s}_i
\end{eqnarray*}
for $2\leq j\leq N.$ Let $\ml{F}_j$ be as in \eqref{cheery_coffee}. Next we will verify  that, for each $N\geq 2$,  $\{X_j;\, 2\leq j \leq N\}$ forms a sequence of martingale differences with respect to the $\sigma$-algebras  $\{\mathcal{F}_{j};\, 1\le j\le N-1\}.$
Define  $J_1=0$ and
\begin{align*}
J_j=\sum_{i=1}^{j-1}T\hat{\rho}_{ij}^2
\end{align*}
for $2\leq j \leq N-1$.
By Lemma \ref{long_for}, $\hat{\rho}_{ij}$ depends on $\bd{s}_i$ and $\bd{s}_j$ only. From independence of $\{\bd{s}_1, \cdots, \bd{s}_N\}$ and Lemma \ref{use_tea},
\begin{eqnarray*}
E(J_j|\ml{F}_{j-1})=\sum_{i=1}^{j-1}E(T\hat{\rho}_{ij}^2|\bd{s}_i)=\frac{T}{m}\sum_{i=1}^{j-1}
\bd{s}_i'\bd{M}_{ij}\bd{s}_i
\end{eqnarray*}
for $2\leq j \leq N,$ where $\bd{M}_{ij}=\bd{U}_i'\bd{U}_j\bd{U}_j'\bd{U}_i$.  Therefore,
\begin{align}\lbl{chat_2}
X_j=J_j-E(J_j|\ml{F}_{j-1}),\ 2\leq j \leq N,
\end{align}
forms a martingale difference with respect to the $\sigma$-algebras $\{\mathcal{F}_{j}; 2\le j\le N\}.$

Now, in order to prove
\begin{eqnarray*}
\frac{1}{N}\sum_{j=2}^NX_j\to N(0, 1)
\end{eqnarray*}
in distribution as $N\to\infty$, we will employ the Lindeberg-Feller central limit theorem
(see, for example, p. 476 from \cite{Billingsley95} or p. 344 from \cite{Durrett}). To achieve so,  it is enough to verify that
 \be\lbl{squareX21}
 Z_N:=\frac{1}{N^2} \sum_{j=2}^NE [X_j^2|\mathcal{F}_{j-1}]\to  1
 \ee
 in probability and
\begin{align}\lbl{fourX13}
\frac{1}{N^4} \sum_{j=2}^N E\big(X_j^4|\mathcal{F}_{j-1}\big)\to 0
\end{align}
in probability as $N\to\infty.$ Lemma \ref{lost_found} has showed \eqref{fourX13}. Now, to prove \eqref{squareX21}, it suffices to show
\begin{align}\lbl{Weight_loss7}
E (Z_N)\rightarrow 1
\end{align}
and
\begin{align}\lbl{Weight_loss8}
 {\rm Var}(Z_N)\rightarrow 0
\end{align}
as $N\to\infty.$ Lemma \ref{kiss_rain} proves \eqref{Weight_loss7} under the assumption $N=o(T^2)$. The assertion \eqref{Weight_loss8} is confirmed in Lemma \ref{ant_egg} by assuming $N=o(T^3)$. Inspect all restrictions between $N$ and $T$ in the lemmas used earlier, the condition  $N=o(T^2)$ meets all requirement. The proof is then completed. \hfill$\square$

\subsubsection{Finale: proofs of Theorems \ref{diligent} and \ref{Fall_nice}}
\label{PS3}

With the preparations in Sections in \ref{Pre_123}-\ref{PS4}, we now are ready to prove the central limit theorem stated in Theorem \ref{diligent}. The main idea is to write the sum of squares of sample correlation coefficients as sums of martingale differences. Then the Lindeberg-Feller martingale CLT is applied.

\smallskip

\noindent\textbf{Proof of Theorem \ref{diligent}}. 
Review  $J_1=0$ and
\begin{align*}
J_j=\sum_{i=1}^{j-1}T\hat{\rho}_{ij}^2
\end{align*}
for $2\leq j \leq N-1$. Then $S_N=\sum_{j=2}^NJ_j.$ Review $\ml{F}_0$ and $\ml{F}_j$ in \eqref{cheery_coffee}.
By Lemma \ref{use_tea},  the conditional expectation,
\begin{eqnarray*}
B_j:=E(J_j|\ml{F}_{j-1})=\frac{T}{m}\sum_{i=1}^{j-1}
\bd{s}_i'\bd{M}_{ij}\bd{s}_i
\end{eqnarray*}
for $2\leq j \leq N,$ where $\bd{M}_{ij}=\bd{U}_i'\bd{U}_j\bd{U}_j'\bd{U}_i$.  As in \eqref{chat_2},
\begin{align*}
X_j=J_j-E(J_j|\ml{F}_{j-1}),\ 2\leq j \leq N,
\end{align*}
forms a martingale difference with respect to the $\sigma$-algebras $\{\mathcal{F}_{j}; 2\le j\le N\}.$ Therefore $\frac{1}{N}(S_N-\mu_N)$ can be further written by
\begin{eqnarray*}
\frac{1}{N}(S_N-\mu_N)=\frac{1}{N}\Big(\sum_{j=2}^NX_j\Big) + \frac{1}{N}\Big[ \Big(\sum_{j=2}^NB_j\Big)-\mu_N\Big].
\end{eqnarray*}
From Lemma \ref{leaves_potato},
\begin{eqnarray*}
\frac{1}{N}\Big[ \Big(\sum_{j=2}^NB_j\Big)-\mu_N\Big] \to 0
\end{eqnarray*}
in probability as $N\to\infty$. By Lemma \ref{CLT_step},
\begin{eqnarray*}
\frac{1}{N}\sum_{j=2}^NX_j\to N(0, 1)
\end{eqnarray*}
in distribution as $N\to\infty.$ The proof then follows from the Slutsky lemma. \hfill$\square$

\smallskip

\noindent\textbf{Proof of Theorem \ref{Fall_nice}}. Set $m=T-p$. First,
\begin{align}\lbl{sun_tai}
\sqrt{\frac{N-1}{N}}\cdot Q_N=\frac{\sqrt{2}}{N}\cdot\sum_{i=1}^{N-1}\sum_{j=i+1}^N
\frac{m\hat{\rho}_{ij}^2-\mu_{Nij}}{v_{Nij}}.
\end{align}
It is easy to see
\begin{eqnarray*}
a_{2N}=\frac{3}{T^2}\Big[1+O\Big(\frac{1}{T}\Big)\Big]
\end{eqnarray*}
as $N\to\infty$. It follows that
\begin{eqnarray*}
a_{1N}&=& \frac{3}{T^2}\Big[1+O\Big(\frac{1}{T}\Big)\Big]-\frac{1}{m^2}\\
&= & \frac{3}{T^2}\Big[1+O\Big(\frac{1}{T}\Big)\Big]-\frac{1}{T^2}\cdot\Big[1+O\Big(\frac{1}{T}\Big)\Big]\\
&= & \frac{2}{T^2}\Big[1+O\Big(\frac{1}{T}\Big)\Big].
\end{eqnarray*}
By Lemma \ref{Love_sound}(i) and (ii), there exists a constant  $K>0$ depending on $p$ but not on $N$, $T$ or $\bd{P}_i$'s such that
\begin{eqnarray*}
\frac{1}{T}\cdot \big|[\mbox{tr}(\bd{P}_i\bd{P}_j)]^2-T^2\big|\leq K \ \ \mbox{and}\ \ \   \big|\mbox{tr}[(\bd{P}_i\bd{P}_j)^2]-T\big|\leq K
\end{eqnarray*}
uniformly for all $1\leq i< j \leq N$ and $N\geq 4.$ Therefore, by the definition of $v_{Nij}$, we have
\begin{eqnarray*}
v_{Nij}^2&=&\frac{2}{T^2}\Big[1+O\Big(\frac{1}{T}\Big)\Big]\cdot [T^2+O(T)]+\frac{6}{T^2}\Big[1+O\Big(\frac{1}{T}\Big)\Big]\cdot [T+O(1)]\\
&=& 2+O\Big(\frac{1}{T}\Big)
\end{eqnarray*}
uniformly for all $1\leq i< j \leq N$ as $N\to\infty.$ Immediately,
\begin{align}\lbl{nice_edu}
\frac{1}{v_{Nij}}=\frac{1}{\sqrt{2}}+O\Big(\frac{1}{T}\Big)
\end{align}
uniformly for all $1\leq i< j \leq N$ as $N\to\infty.$ Now write $\frac{1}{v_{Nij}}=\frac{1}{\sqrt{2}}(1+\omega_{Nij})$. Then
\begin{align}\lbl{Moon_Cake}
\sup_{1\leq i< j \leq N}|\omega_{Nij}|\leq \frac{C}{T}
\end{align}
as $N\geq 4.$ By Lemma \ref{use_tea}, $E\hat{\rho}_{ij}^2=m^{-1}\mu_{Nij}.$ It follows that
\begin{eqnarray*}
&& \frac{\sqrt{2}}{N}\sum_{i=1}^{N-1}\sum_{j=i+1}^N
\frac{m\hat{\rho}_{ij}^2-\mu_{Nij}}{v_{Nij}}\\
&=  &  \frac{m}{NT}\sum_{i=1}^{N-1}\sum_{j=i+1}^N
T[\hat{\rho}_{ij}^2-m^{-1}\mu_{Nij}](1+\omega_{Nij})\\
& = & \frac{m}{T}\cdot\frac{1}{N} (S_N-\mu_N) + \frac{m}{N}\sum_{i=1}^{N-1}\sum_{j=i+1}^N \omega_{Nij}
(\hat{\rho}_{ij}^2-E\hat{\rho}_{ij}^2),
\end{eqnarray*}
where $S_N$ is defined as in Theorem \ref{diligent}. By the Slutsky lemma and Theorem \ref{diligent},
\begin{eqnarray*}
\frac{m}{T}\cdot\frac{1}{N} (S_N-\mu_N)\to N(0, 1)
\end{eqnarray*}
in distribution as $N\to \infty.$ Recall \eqref{sun_tai}. To prove $Q_N\to N(0, 1)$ in distribution, by the Slutsky lemma again, it is enough to show
\begin{align}\lbl{flower_beer}
\Delta_n:=\frac{m}{N}\sum_{i=1}^{N-1}\sum_{j=i+1}^N \omega_{Nij}
(\hat{\rho}_{ij}^2-E\hat{\rho}_{ij}^2) \to 0
\end{align}
in probability as $N\to \infty$. Since $\hat{\rho}_{ij}^2$ and $\hat{\rho}_{kl}^2$ are independent if $\{i, j\}\cap \{k, l\}=\emptyset$, then
\begin{eqnarray*}
\mbox{Var}(\Delta_n) =
\Big(\frac{m}{N}\Big)^2\sum_{1\leq i<j \leq N}\sum\omega_{Nij}^2\mbox{Cov}(\hat{\rho}_{ij}^2, \hat{\rho}_{kl}^2)
\end{eqnarray*}
where the last sum runs over all $(k, l)$ with $1\leq k<l\leq N$ and $\{i, j\}\cap \{k, l\}\ne \emptyset$. The total number of such $(k, l)$'s is no more than $2N+ 2N=4N.$ Since $|\mbox{Cov}(U, V)|\leq [\mbox{Var}(U)]^{1/2}\cdot [\mbox{Var}(V)]^{1/2}$ for any random variables $U$ and $V$, we have from  \eqref{Moon_Cake} that
\begin{align}\lbl{Huier}
\mbox{Var}(\Delta_n) \leq & \Big(\frac{m}{N}\Big)^2\cdot \frac{C^2}{T^2}\cdot \frac{1}{2}N(N-1)\cdot (4N)\cdot \max_{1\leq i< j \leq N}\mbox{Var}(\hat{\rho}_{ij}^2)\nonumber\\
 \leq & (2C^2N)\cdot \max_{1\leq i< j \leq N}\mbox{Var}(\hat{\rho}_{ij}^2).
\end{align}
By Lemma \ref{trivial}, $\mbox{tr}(\bd{P}_i\bd{P}_j)\leq \mbox{tr}(\bd{P}_i)= m$ and $\mbox{tr}[(\bd{P}_i\bd{P}_j)^2]\leq \mbox{tr}(\bd{P}_i)= m$. By Lemma \ref{good_tie}(iv),
\begin{eqnarray*}
\mbox{Var}(\hat{\rho}_{ij}^2)
& = & \frac{6}{m^2(m+2)^2}\cdot \mbox{tr}[(\bd{P}_i\bd{P}_j)^2]+\frac{2(m^2-2m-2)}{m^4(m+2)^2}\cdot [\mbox{tr}(\bd{P}_i\bd{P}_j)]^2\\
& \leq & \frac{6}{m^3}+ \frac{2}{m^2}.
\end{eqnarray*}
Thus, $\mbox{Var}(\hat{\rho}_{ij}^2)\leq \frac{8}{m^2}$ for all $1\leq i<j \leq N.$ Combing this with \eqref{Huier}, we get
\begin{eqnarray*}
\mbox{Var}(\Delta_n)\leq C\cdot \frac{N}{T^2}\to 0
\end{eqnarray*}
by the assumption $N=o(T^2).$ This implies \eqref{flower_beer}. The proof is completed.\hfill$\square$

\subsection{The proofs of Theorems \ref{linear_case}, \ref{exponential_case} and \ref{super_exponential_case} }

Theorems \ref{linear_case}-\ref{super_exponential_case} will be proved via approximating $L_N=\max_{1\le i<j\le N} |\hat{\rho}_{ij}|$ for any $p\geq 0$ by $L_N$ for the case $p=0$. The latter one has the asymptotic result known in \cite{CJF13}. The main job is reduced to show the difference between the two versions of $L_N$ is small enough.

\subsubsection{Prelude: auxilary results towards proofs of Theorems \ref{linear_case}, \ref{exponential_case} and \ref{super_exponential_case}}\lbl{Pre_linear}

The results stated in this section  will be proved in Section \ref{Pre_linear_app}.

\begin{lemma}\lbl{shock}  Let $m\geq 1$ and $\{\xi_i;\, 1\leq i \leq m\}$ be i.i.d. random variables with $E\xi_1=0$,\, $E\xi_1^2=1$  and $E(|\xi_1|^\tau)<\infty$ for some $\tau\geq 2$. Let $\{a_i;\, 1\leq i \leq m\}$ be constants such that $a_1^2+\cdots + a_m^2=1.$ Then, there exists a constant $K>0$ satisfying
\begin{eqnarray*}
P(|a_1\xi_1+\cdots+ a_m\xi_m|\geq x) \leq \frac{K}{x^\tau}
\end{eqnarray*}
for all $x\geq 3.$
\end{lemma}

It is easy to see that the bound in the lemma is tight by simply taking $a_1=1$ and $a_2=\cdots=a_m=0.$

\begin{lemma}\lbl{Root}  Let $m\geq 1$ and $\{\xi_i;\, 1\leq i \leq m\}$ be i.i.d. random variables with $E\xi_1=0$,\, $E\xi_1^2=1$  and $Ee^{\omega|\xi_1|}<\infty$ for some $\omega>0$. Let $\{a_i;\, 1\leq i \leq m\}$ be constants satisfying $a_1^2+\cdots + a_m^2=1.$ Then, there exists $K>0$ such that
\begin{eqnarray*}
P(|a_1\xi_1+\cdots+ a_m\xi_m|\geq x) \leq K\cdot e^{-x/K}
\end{eqnarray*}
for all $x\geq 0.$
\end{lemma}

The above inequality is tight, which can be seen by taking $a_1=1$ and  $a_i=0$ for $2\leq i \leq m.$

Recall the definition of subgaussian random variables defined before the statement of Theorem \ref{super_exponential_case}.

\begin{lemma}\lbl{sleepless} Let $m\geq 1$ and $\{\xi_i;\, 1\leq i \leq m\}$ be i.i.d. subgaussian random variables. Let $\{a_i;\, 1\leq i \leq m\}$ be constants such that $a_1^2+\cdots + a_m^2=1.$ Then, there exists a positive constant $K$ not depending on $m$ or $\{a_i;\, 1\leq i \leq m\}$ such that
\begin{eqnarray*}
P(|a_1\xi_1+\cdots+ a_m\xi_m|\geq x) \leq 2\cdot e^{-K x^2}
\end{eqnarray*}
for all $x>0.$
\end{lemma}

The upper bound in the lemma is optimized, which can be seen evidently by choosing $a_1=1$ and $a_2=\cdots=a_m=0$.

\subsubsection{Intermezzo: approximation of sample correlation coefficients by simple versions}

Recall the setting in \eqref{shadongxi}, \eqref{park_1} and \eqref{Asia_foods}.
Let $p$ be fixed.  Let  $\bd{e}_i=\frac{\bd{\epsilon}_i}{\|\bd{\epsilon}_i\|}$ and $\tilde{\rho}_{ij}=\bd{e}_i'\bd{e}_j$ for $1\leq i, j \leq N$. In this section, we always assume  that $\{\epsilon_{ij};\, i\geq 1, j \geq 1 \}$ are i.i.d. continuous random variables. The ``continuous" requirement guarantees that $\hat{\rho}_{ij}$ in \eqref{Asia_foods} is well-defined. See the comment below \eqref{assa}.

\begin{lemma}\lbl{Seattle} 
Assume $\{\epsilon_{ij};\, i\geq 1, j \geq 1 \}$ are i.i.d. continuous random variables. Let $\hat{\rho}_{ij}$ be defined as in \eqref{Asia_foods}.
Set $\bd{A}_i=\bd{x}_i(\bd{x}_i'\bd{x}_i)^{-1}\bd{x}_i'$ for $1\leq i \leq N$. Then,
\begin{eqnarray*}
\max_{1\leq i <j \leq N}\big|\hat{\rho}_{ij} - \tilde{\rho}_{ij}\big| \leq
14\cdot \big(\max_{1\leq i \leq j \leq N}\bd{e}_j'\bd{A}_i\bd{e}_j\big).
\end{eqnarray*}
\end{lemma}
\noindent\textbf{Proof of Lemma \ref{Seattle}}. Notice
 $\bd{P}_i\bd{\epsilon}_i=\bd{\epsilon}_i-\bd{A}_i\bd{\epsilon}_i.$ It follows that
\begin{eqnarray*}
 \bd{\epsilon}_i'\bd{P}_i\bd{P}_j\bd{\epsilon}_j&=&(\bd{\epsilon}_i-\bd{A}_i\bd{\epsilon}_i)'
(\bd{\epsilon}_j-\bd{A}_j\bd{\epsilon}_j)\\
&= & \bd{\epsilon}_i'\bd{\epsilon}_j-\bd{\epsilon}_i'\bd{A}_i\bd{\epsilon}_j-\bd{\epsilon}_i'\bd{A}_j\bd{\epsilon}_j
+\bd{\epsilon}_i'\bd{A}_i\bd{A}_j\bd{\epsilon}_j.
\end{eqnarray*}
Take $i=j$ to see that
\begin{eqnarray*}
\|\bd{P}_i\bd{\epsilon}_i\|^2=\|\bd{\epsilon}_i\|^2- \bd{\epsilon}_i'\bd{A}_i\bd{\epsilon}_i,
\end{eqnarray*}
since $\bd{A}_i^2=\bd{A}_i$. In particular, $\|\bd{\epsilon}_i\|^2- \bd{\epsilon}_i'\bd{A}_i\bd{\epsilon}_i\geq 0$, hence
\begin{align}\lbl{sucaknvt}
0 \leq  \bd{e}_i' \bd{A}_i\bd{e}_i \leq 1
\end{align}
for each $i$. Combining the last two identities, we have
\begin{eqnarray*}
\hat{\rho}_{ij}=\frac{\bd{\epsilon}_i'\bd{\epsilon}_j-\bd{\epsilon}_i'\bd{A}_i\bd{\epsilon}_j
-\bd{\epsilon}_i'\bd{A}_j\bd{\epsilon}_j
+\bd{\epsilon}_i'\bd{A}_i\bd{A}_j\bd{\epsilon}_j}{\sqrt{\|\bd{\epsilon}_i\|^2- \bd{\epsilon}_i'\bd{A}_i\bd{\epsilon}_i}\cdot\sqrt{\|\bd{\epsilon}_j\|^2- \bd{\epsilon}_j'\bd{A}_j\bd{\epsilon}_j}}.
\end{eqnarray*}
Dividing the numerator and denominator by $\|\bd{\epsilon}_i\|\cdot \|\bd{\epsilon}_j\|$, we have
\begin{align}
\hat{\rho}_{ij}=&\big[\tilde{\rho}_{ij}-\bd{e}_i'\bd{A}_i\bd{e}_j
-\bd{e}_i'\bd{A}_j\bd{e}_j
+\bd{e}_i'\bd{A}_i\bd{A}_j\bd{e}_j\big]\nonumber\\
&\cdot \big(1-\bd{e}_i'\bd{A}_i\bd{e}_i\big)^{-1/2}
\big(1-\bd{e}_j'\bd{A}_j\bd{e}_j\big)^{-1/2}.\lbl{39671}
\end{align}
Write
\begin{eqnarray*}
1-\frac{1}{\sqrt{1-x}}=\frac{\sqrt{1-x}-1}{\sqrt{1-x}}=\frac{-x}{(\sqrt{1-x}+1)\sqrt{1-x}}.
\end{eqnarray*}
It is easy to see $|1-(1-x)^{-1/2}|\leq 2|x|$ if $|x|\leq \frac{1}{2}.$ Then  $(1-x)^{-1/2}\leq 1+2|x|$ as $|x|\leq \frac{1}{2}.$ For brevity of notation, set $h_{ij}=\big(1-\bd{e}_i'\bd{A}_i\bd{e}_i\big)^{-1/2}
\big(1-\bd{e}_j'\bd{A}_j\bd{e}_j\big)^{-1/2}.$ Then
\begin{eqnarray*}
0\leq h_{ij}-1 &\leq & (1+2\bd{e}_i'\bd{A}_i\bd{e}_i)
(1+2\bd{e}_j'\bd{A}_j\bd{e}_j)-1\\
& = & 2\bd{e}_i'{A}_i\bd{e}_i+
2\bd{e}_j'\bd{A}_j\bd{e}_j + 4(\bd{e}_i'\bd{A}_i\bd{e}_i)\cdot (\bd{e}_j'\bd{A}_j\bd{e}_j)\\
& \leq & 4(\bd{e}_i'{A}_i\bd{e}_i+
\bd{e}_j'\bd{A}_j\bd{e}_j)
\end{eqnarray*}
provided $\max_{1\leq i  \leq N}\bd{e}_i'\bd{A}_i\bd{e}_i\leq \frac{1}{2}$, and at the same time  $h_{ij}\leq 2$ by definition. From \eqref{39671},
\begin{align}\lbl{ye_menglong}
\hat{\rho}_{ij}=\tilde{\rho}_{ij} + \tilde{\rho}_{ij}(h_{ij}-1)+(-\bd{e}_i'\bd{A}_i\bd{e}_j
-\bd{e}_i'\bd{A}_j\bd{e}_j
+\bd{e}_i'\bd{A}_i\bd{A}_j\bd{e}_j)\cdot h_{ij}.
\end{align}
By the Cauchy-Schwartz inequality and the fact $\bd{A}_i^2=\bd{A}_i$,
\begin{eqnarray*}
|\bd{e}_i'\bd{A}_i\bd{A}_j\bd{e}_j|\leq \|\bd{A}_i\bd{e}_i\|\cdot\|\bd{A}_j\bd{e}_j\|.
\end{eqnarray*}
Similarly, $|\bd{e}_i'\bd{A}_i\bd{e}_j|\leq \|\bd{A}_i\bd{e}_i\|\cdot\|\bd{A}_i\bd{e}_j\|$ and $|\bd{e}_i'\bd{A}_j\bd{e}_j| \leq \|\bd{A}_j\bd{e}_i\|\cdot\|\bd{A}_j\bd{e}_j\|$ since $\bd{A}_i^2=\bd{A}_i$.
Consequently
\begin{eqnarray*}
& & \big| \tilde{\rho}_{ij}(h_{ij}-1)+\big(-\bd{e}_i'\bd{A}_i\bd{e}_j
-\bd{e}_i'\bd{A}_j\bd{e}_j
+\bd{e}_i'\bd{A}_i\bd{A}_j\bd{e}_j\big)\cdot h_{ij}\big|\\
& \leq & 4(\bd{e}_i'{A}_i\bd{e}_i+
\bd{e}_j'\bd{A}_j\bd{e}_j)\\
& &+ 2\big(\|\bd{A}_i\bd{e}_i\|\cdot\|\bd{A}_i\bd{e}_j\|+ \|\bd{A}_j\bd{e}_i\|\cdot\|\bd{A}_j\bd{e}_j\|+
\|\bd{A}_i\bd{e}_i\|\cdot\|\bd{A}_j\bd{e}_j\|\big)
\end{eqnarray*}
by the fact $|\tilde{\rho}_{ij}| \leq 1$. 
Use the trivial fact that $2xy\leq x^2+y^2$ to see
\begin{eqnarray*}
& & \max_{1\leq i < j \leq N}\big| \tilde{\rho}_{ij}(h_{ij}-1)+\big(-\bd{e}_i'\bd{A}_i\bd{e}_j
-\bd{e}_i'\bd{A}_j\bd{e}_j
+\bd{e}_i'\bd{A}_i\bd{A}_j\bd{e}_j\big)\cdot h_{ij}\big|\\
&\leq& \max_{1\leq i < j \leq N}\{6 (\bd{e}_i'{A}_i\bd{e}_i+
\bd{e}_j'\bd{A}_j\bd{e}_j)+(\bd{e}_j'{A}_i\bd{e}_j+
\bd{e}_i'\bd{A}_j\bd{e}_i)\}\\
&\leq & 14\cdot \big(\max_{1\leq i \leq j \leq N}\bd{e}_j'\bd{A}_i\bd{e}_j\big)
\end{eqnarray*}
provided $\max_{1\leq i  \leq N}\bd{e}_i'\bd{A}_i\bd{e}_i\leq \frac{1}{2}$, where in the last inequality we use the fact that each term is bounded by $\max_{1\leq i \leq j \leq N}\bd{e}_j'\bd{A}_i\bd{e}_j$. We then have from \eqref{ye_menglong} that
\begin{align}\lbl{ai_qin}
\max_{1\leq i <j \leq N}\big|\hat{\rho}_{ij} - \tilde{\rho}_{ij}\big| \leq 14\cdot \big(\max_{1\leq i \leq j \leq N}\bd{e}_j'\bd{A}_i\bd{e}_j\big)
\end{align}
provided $\max_{1\leq i  \leq N}\bd{e}_i'\bd{A}_i\bd{e}_i\leq \frac{1}{2}$. If $\max_{1\leq i  \leq N}\bd{e}_i'\bd{A}_i\bd{e}_i> \frac{1}{2}$, \eqref{ai_qin} holds automatically due to the facts that $|\hat{\rho}_{ij}|\leq 1$ and $|\tilde{\rho}_{ij}|\leq 1$ for all $i,j$. The proof is completed. \hfill$\Box$

\begin{prop}\lbl{Fall} Assume $\{\epsilon_{ij};\, i\geq 1, j \geq 1 \}$ are i.i.d. continuous random variables with $E\epsilon_{11}=0$ and  $E(|\epsilon_{11}|^\tau)<\infty$ for some $\tau\geq 4$.  Suppose $T/(N^{8/\tau}\log N)\to \infty$ as $N\to\infty$. Then
\begin{equation*}
\sqrt{T\log N}\cdot\max_{1\leq i <j \leq N}\big|\hat{\rho}_{ij} - \tilde{\rho}_{ij}\big|\to 0
\end{equation*}
in probability as $N\to\infty.$
\end{prop}
\noindent\textbf{Proof of Proposition \ref{Fall}}. To prove the result,  by the homogeneity of $\hat{\rho}_{ij}$ from \eqref{Asia_foods}, without loss of generality, we assume $E(\bd{\epsilon}_{11}^2)=1$. Set $\alpha_N=1/\sqrt{T\log N}$. Then $\alpha_N\to 0$ as $N\to\infty$ by assumption. From Lemma \ref{Seattle}, for any $h\in (0, 14)$,
\begin{equation}
 P\Big(\max_{1\leq i <j \leq N}\big|\hat{\rho}_{ij} - \tilde{\rho}_{ij}\big|\geq \alpha_N h\Big)
\leq   P\Big(\max_{1\leq i \leq j \leq N}\bd{e}_j'\bd{A}_i\bd{e}_j> \frac{h}{14}\alpha_N\Big).\lbl{subway1}
\end{equation}
Next we estimate the last probability.

For any $v>0$,
\begin{align}
P\Big(\max_{1\leq i \leq j \leq N}\bd{e}_j'\bd{A}_i\bd{e}_j> 2v\alpha_N\Big)  \leq  & N^2\cdot \max_{1\leq i \leq j \leq N}P\Big(\bd{e}_j'\bd{A}_i\bd{e}_j> 2\alpha_Nv\Big)\nonumber\\
  = &  N^2\cdot \max_{1\leq i  \leq N}P\Big(\bd{e}_1'\bd{A}_i\bd{e}_1> 2\alpha_Nv\Big).\label{good_eat}
\end{align}
Now
\begin{align}
P\big(\bd{e}_1'\bd{A}_i\bd{e}_1> 2\alpha_N v\big)
=& P\big(\bd{\epsilon}_1'\bd{A}_i\bd{\epsilon}_1> 2\alpha_N v\|\bd{\epsilon}_1\|^2\big)\nonumber\\
 \leq & P\big(\bd{\epsilon}_1'\bd{A}_i\bd{\epsilon}_1> v\sqrt{T/\log N}\big)+ P\Big(\|\bd{\epsilon}_1\|^2\leq \frac{1}{2}T\Big).\lbl{love1}
\end{align}
Note that $\|\bd{\epsilon}_1\|^2=\sum_{j=1}^T\bd{\epsilon}_{1j}^2$. Since $E(\bd{\epsilon}_{11}^2)=1$ and $E(|\bd{\epsilon}_{11}|^\tau)<\infty$, we see
\begin{align}\lbl{hahaha}
P\Big(\|\bd{\epsilon}_1\|^2\leq \frac{1}{2}T\Big)
 \leq & P\Big(\Big|\sum_{j=1}^T(\epsilon_{1j}^2-1)\Big|\geq  \frac{1}{2}T\Big)\nonumber\\
 \leq & \Big(\frac{2}{T}\Big)^{\tau/2}\cdot E\Big(\Big|\sum_{j=1}^T(\epsilon_{1j}^2-1)\Big|^{\tau/2}\Big) \nonumber\\
 = & O\Big(\frac{1}{T^{\tau/4}}\Big)
\end{align}
by \eqref{hot_sleep2}, where the Markov inequality is applied in the second inequality. Write $\bd{Q}=\Gamma'\bd{D}\Gamma$ where $\Gamma=(\gamma_{ij})_{T\times T}$ is an orthogonal matrix and
\begin{eqnarray*}
\bd{D}=\begin{pmatrix}
\bd{I}_p & \bd{0}\\
\bd{0} & \bd{0}
\end{pmatrix}
.
\end{eqnarray*}
Then $\bd{\epsilon}_1'\bd{Q}\bd{\epsilon}_1=\sum_{k=1}^p\big(\sum_{j=1}^T\gamma_{kj}\epsilon_{1j}\big)^2$. It follows that
\begin{align}
 & P\big(\bd{\epsilon}_1'\bd{Q}\bd{\epsilon}_1> v\sqrt{T/\log N}\big) \nonumber\\
\leq  &  p\cdot \max_{1\leq k \leq p}P\Big(\big(\sum_{j=1}^T\gamma_{kj}\epsilon_{1j}\big)^2> \frac{v}{p}\sqrt{T/\log N}\Big) \nonumber\\
 = & p\cdot \max_{1\leq k \leq p}P\Big(\big|\sum_{j=1}^T\gamma_{kj}\epsilon_{1j}\big|> v'(T/\log N)^{1/4}\Big), \lbl{bird1}
\end{align}
where $v':=(v/p)^{1/2}.$ Note that $\sum_{j=1}^T\gamma_{kj}^2=1$ for each $1\leq k \leq p$ by orthogonality. Thus, from Lemma \ref{shock} we have that there exists some $K>0$, such that
\begin{eqnarray*}
 P\big(\bd{\epsilon}_1'\bd{Q}\bd{\epsilon}_1> v\sqrt{T/\log N}\big)
 \leq \frac{pK}{v'^\tau}\cdot (T/\log N)^{-\tau/4}.
\end{eqnarray*}
Join this with \eqref{good_eat}, \eqref{love1} and \eqref{hahaha} to get
 \begin{eqnarray*}
 P\big(\max_{1\leq i \leq j \leq N}\bd{e}_j'\bd{Q}\bd{e}_j> 2\alpha_N v\big)= O\Big(\frac{N^2}{(T/\log N)^{\tau/4}}\Big)+ O\Big(\frac{N^2}{T^{\tau/4}}\Big).
 \end{eqnarray*}
By taking $v=\frac{h}{28}$, we have from the above  and \eqref{subway1} that
\begin{eqnarray*}
P\Big(\max_{1\leq i <j \leq N}\big|\hat{\rho}_{ij} - \tilde{\rho}_{ij}\big|\geq \alpha_N h\Big)
 &\leq & O\Big(\frac{N^2}{(T/\log N)^{\tau/4}}\Big)+ O\Big(\frac{N^2}{T^{\tau/4}}\Big)\\
 & = &  O\Big(\frac{N^2(\log N)^{\tau/4}}{T^{\tau/4}}\Big)+o\Big(\frac{1}{(\log N)^{\tau/4}}\Big),
\end{eqnarray*}
which goes to zero by the assumption that $T/(N^{8/\tau}\log N)\to \infty$ as $N\to\infty$. \hfill$\Box$

\begin{prop}\lbl{Winter} Assume $\{\epsilon_{ij};\, i\geq 1, j \geq 1 \}$ are i.i.d. continuous random variables with $E\epsilon_{11}=0$ and  $Ee^{\omega|\epsilon_{11}|}<\infty$ for some $\omega>0$. Suppose $\log N=o(T^{1/5})$ as $N\to\infty$. Then
\begin{eqnarray*}
\sqrt{T\log N}\cdot\max_{1\leq i <j \leq N}\big|\hat{\rho}_{ij} - \tilde{\rho}_{ij}\big|\to 0
\end{eqnarray*}
in probability as $N\to\infty.$
\end{prop}
\noindent\textbf{Proof of Proposition \ref{Winter}}. To prove the result,  by the homogeneity of $\hat{\rho}_{ij}$ from \eqref{Asia_foods}, without loss of generality, assume $E(\bd{\epsilon}_{11}^2)=1$. Set $\alpha_N=1/\sqrt{T\log N}$. Then $\alpha_N\to 0$ as $N\to\infty$ by assumption. From \eqref{subway1}, we have that, for any $h\in (0, 14)$,
\begin{align}
 P\Big(\max_{1\leq i <j \leq N}\big|\hat{\rho}_{ij} - \tilde{\rho}_{ij}\big|\geq \alpha_N h\Big)
\leq   P\Big(\max_{1\leq i \leq j \leq N}\bd{e}_j'\bd{A}_i\bd{e}_j> \frac{h}{14}\alpha_N\Big),\lbl{subway}
\end{align}
as $N$ is sufficiently large. So to finish the proof it suffices to show the second probability goes to zero.

For any $v>0$, by \eqref{good_eat} and \eqref{love1},
\begin{align}\lbl{rest_lu}
 &  P\big(\max_{1\leq i \leq j \leq N}\bd{e}_j'\bd{A}_i\bd{e}_j> 2\alpha_N v\big)\nonumber\\
\leq & N^2\cdot \max_{1\leq i  \leq N} P\big(\bd{\epsilon}_1'\bd{A}_i\bd{\epsilon}_1> v\sqrt{T/\log N}\big)+ N^2\cdot P\Big(\|\bd{\epsilon}_1\|^2\leq \frac{1}{2}T\Big).
\end{align}
Since $E(\bd{\epsilon}_{11}^2)=1$, by large deviations, there exists a constant $\eta_0>0$ such that
\begin{align}\lbl{control}
P\Big(\|\bd{\epsilon}_1\|^2\leq \frac{1}{2}T\Big) =P\Big(\frac{1}{T}\sum_{j=1}^T\epsilon_{1j}^2<\frac{1}{2}\Big)\leq  e^{-\eta_0 T}
\end{align}
for large enough $T$; see, for example, \cite{DZ98}. By \eqref{bird1},
\begin{eqnarray*}
&& P\big(\bd{\epsilon}_1'\bd{A}_i\bd{\epsilon}_1> v\sqrt{T/\log N}\big) \\
& \leq & p\cdot \max_{1\leq k \leq p}P\Big(\big|\sum_{j=1}^T\gamma_{kj}\epsilon_{1j}\big|> v'(T/\log N)^{1/4}\Big),
\end{eqnarray*}
where $v':=(v/p)^{1/2}.$ Note that $\sum_{j=1}^T\gamma_{kj}^2=1$ for each $1\leq k \leq p$ by orthogonality. From Lemma \ref{Root}, there exists $K>0$ such that
\begin{eqnarray*}
& & P\big(\max_{1\leq i \leq j \leq N}\bd{e}_j'\bd{A}_i\bd{e}_j> 2\alpha_N v\big)\\
& \leq &   N^2\cdot \Big[p\cdot \max_{1\leq k \leq p}P\Big(\big|\sum_{j=1}^T\gamma_{kj}\epsilon_{1j}\big|> v'(T/\log N)^{1/4}\Big)+ e^{-\eta_0 T}\Big]\\
& \leq & N^2\cdot\Big[(pK)e^{-(v'/K)(T/\log N)^{1/4}}+e^{-\eta_0 T}\Big]
\end{eqnarray*}
 as $T$ is large enough, where $v':=(v/p)^{1/2}$ and $K>0$ is a constant not depending on $p$, $N$, $T$ or $\gamma_{kj}$'s. It is easy to see the above goes to zero if $T/(\log N)^5\to \infty.$ It follows that
\begin{eqnarray*}
P\big(\max_{1\leq i \leq j \leq N}\bd{e}_j'\bd{A}_i\bd{e}_j> 2\alpha_N v\big)\to 0
\end{eqnarray*}
provided $T/(\log N)^5\to \infty.$ The proof is completed. \hfill$\Box$

\begin{prop}\lbl{huns}  Assume $\{\epsilon_{ij};\, i\geq 1, j \geq 1 \}$ are i.i.d. continuous and subgaussian random variables.
If $\log N=o(T^{1/3})$, then
\begin{eqnarray*}
\sqrt{T\log N}\cdot\max_{1\leq i <j \leq N}\big|\hat{\rho}_{ij} - \tilde{\rho}_{ij}\big|\to 0
\end{eqnarray*}
in probability as $N\to\infty.$
\end{prop}
\noindent\textbf{Proof of Proposition \ref{huns}}. First, the subgaussian assumption implies that $Ee^{\omega|\epsilon_{11}|^2}<\infty$ for some $\omega>0$. Hence $Ee^{t|\epsilon_{11}|}<\infty$ for all $t>0$. We will use the same notation as in the proof of Proposition \ref{Winter}. Reviewing \eqref{rest_lu} and \eqref{control}, to get our desired result,  it suffices to show that $N^2\cdot \max_{1\leq i  \leq N}P\big(\bd{e}_1'\bd{A}_i\bd{e}_1> 2\alpha_N v\big)\to 0$ for any $v>0$.

By \eqref{bird1}, for each $i$,
\begin{eqnarray*}
&& P\big(\epsilon_1'\bd{A}_i\epsilon_1> v\sqrt{T/\log N}\big) \\
& \leq & p\cdot \max_{1\leq k \leq p}P\Big(\big|\sum_{j=1}^T\gamma_{kj}\epsilon_{1j}\big|> v'(T/\log N)^{1/4}\Big),
\end{eqnarray*}
where $v':=(v/p)^{1/2}.$ Note that $\sum_{j=1}^T\gamma_{kj}^2=1$ for each $1\leq k \leq p$ by orthogonality. From Lemma \ref{sleepless},  there exists $K>0$ such that
\begin{eqnarray*}
P\big(\epsilon_1'\bd{A}_i\epsilon_1> v\sqrt{T/\log N}\big) \leq (2p)\cdot \exp\big(-Kv'^2\sqrt{T/\log N}\big).
\end{eqnarray*}
Therefore, by \eqref{rest_lu} and \eqref{control}, there exists a constant $\beta_0>0$ such that
\begin{eqnarray*}
 N^2\cdot \max_{1\leq i  \leq N}P\big(\bd{e}_1'\bd{A}_i\bd{e}_1> 2\alpha_N v\big)
  \leq (2pN^2)\cdot \exp(-Kv'^2\sqrt{T/\log N}\,)+N^2e^{-\beta_0 T}.
\end{eqnarray*}
It is easy to see the above goes to zero if $\log N=o(T^{1/3}).$   \hfill$\Box$

\smallskip

Assume $\{\epsilon_{ij};\, i\geq 1, j \geq 1 \}$ are i.i.d. continuous random variables. Set $\bar{\bd{\epsilon}}_i=(1/T)\sum_{j=1}^T\epsilon_{ij}$ for all $i$.
Define $\bd{1}=(1, \dots, 1)'\in \mathbb{R}^T.$ The Pearson correlation coefficient $\rho_{ij}$ is then defined by
\begin{align}\lbl{Pearson}
\rho_{ij}=\frac{(\bd{\epsilon}_i-\bar{\bd{\epsilon}}_i\bd{1})'(\epsilon_j-\bar{\epsilon}_j\bd{1})}
{\|\bd{\epsilon}_i-\bar{\bd{\epsilon}}_i\bd{1}\|\cdot \|\epsilon_j-\bar{\epsilon}_j\bd{1}\|}
\end{align}
for $1\leq i, j\leq N.$ Similar to the clarification below \eqref{assa}, the ``i.i.d. continuous" assumption  justifies that $\rho_{ij}$ is well-defined.

\begin{prop}\lbl{BJ} Let $\tilde{\rho}_{ij}$ be as in Lemma \ref{Seattle}. Assume $E\epsilon_{11}=0$ and  $E(|\epsilon_{11}|^\tau)<\infty$ for some $\tau\geq 2$. If $\frac{N^{4/\alpha}\log N}{T}\to 0$, then
\begin{eqnarray*}
\sqrt{T\log N}\cdot \max_{1\leq i <j \leq N}\big|\tilde{\rho}_{ij}-\rho_{ij}\big|\to 0
\end{eqnarray*}
in probability as $N\to\infty$.
\end{prop}
\noindent\textbf{Proof of Proposition \ref{BJ}}. The proof consists of two steps. In the first step we will show
\begin{align}\lbl{honest}
\max_{1\leq i <j \leq N}\big|\tilde{\rho}_{ij}-\rho_{ij}\big| \leq (12T)\cdot \max_{1\leq i  \leq N}\frac{\bar{\bd{\epsilon}}_i^2}{\|\bd{\epsilon}_i\|^2}.
\end{align}
By using this we will prove the desired conclusion in the second step.

{\it Step 1}. Set $\delta_i=\sqrt{T}\cdot\frac{\bar{\bd{\epsilon}}_i}{\|\bd{\epsilon}_i\|}$ for each $i.$ Recall $\tilde{\rho}_{ij}=\frac{\bd{\epsilon}_i'\bd{\epsilon}_j}
{\|\bd{\epsilon}_i\|\cdot\|\bd{\epsilon}_j\|}$. Write
\begin{eqnarray*}
\rho_{ij} &= & \frac{\bd{\epsilon}_i'\bd{\epsilon}_j-T\bar{\bd{\epsilon}}_i\bar{\epsilon}_j}
{\sqrt{\|\bd{\epsilon}_i\|^2-T\bar{\bd{\epsilon}}_i^2}\cdot \sqrt{\|\bd{\epsilon}_j\|^2-T\bar{\epsilon}_j^2}}\\
& = & (\tilde{\rho}_{ij}-\delta_i\delta_j)(1-\delta_i^2)^{-1/2}(1-\delta_j^2)^{-1/2}.
\end{eqnarray*}
It follows that
\begin{align}\lbl{maple}
\rho_{ij}-\tilde{\rho}_{ij}=& \tilde{\rho}_{ij}\cdot \big[(1-\delta_i^2)^{-1/2}(1-\delta_j^2)^{-1/2}-1\big]- \nonumber\\
 & \delta_i\delta_j(1-\delta_i^2)^{-1/2}(1-\delta_j^2)^{-1/2}.
\end{align}
By the inequality $(1-x)^{-1/2}\leq 1+2|x|$ for $|x|\leq \frac{1}{2}$ appeared in the proof of Lemma \ref{Seattle}, we have
\begin{eqnarray*}
0 &\leq & (1-\delta_i^2)^{-1/2}(1-\delta_j^2)^{-1/2}-1\\
&\leq & (1+2\delta_i^2)(1+2\delta_j^2)-1\\
& \leq & 4(\delta_i^2+\delta_j^2)
\end{eqnarray*}
provided $\max_{1\leq i \leq N} |\delta_i|<\frac{1}{2}.$ Under this restriction, $(1-\delta_i^2)^{-1/2}\leq \frac{2}{\sqrt{3}} \leq 2$ for each $i$. Therefore,
\begin{eqnarray*}
\Big|\delta_i\delta_j(1-\delta_i^2)^{-1/2}(1-\delta_j^2)^{-1/2}\Big|\leq 4|\delta_i\delta_j|.
\end{eqnarray*}
Since $|\tilde{\rho}_{ij}|\leq 1$, the above two estimates joining with \eqref{maple} implies that
\begin{eqnarray*}
\max_{1\leq i< j \leq N}|\rho_{ij}-\tilde{\rho}_{ij}|
& \leq & \max_{1\leq i< j \leq N}\big[4(\delta_i^2+\delta_j^2)+4|\delta_i\delta_j|\big]\\
& \leq & 12\cdot \max_{1\leq i \leq N}\delta_i^2
\end{eqnarray*}
as $\max_{1\leq i \leq N} |\delta_i| \leq \frac{1}{2}.$ Moreover, the above naturally holds if $\max_{1\leq i \leq N} |\delta_i| > \frac{1}{2}.$ This leads to \eqref{honest}.

{\it Step 2}. Set $\alpha_N=1/\sqrt{T\log N}$. Then $\lim_{N\to\infty}\alpha_N= 0$. From {\it Step 1}, for any $t>0$,
\begin{eqnarray*}
 P\Big(\max_{1\leq i< j \leq N}|\rho_{ij}-\tilde{\rho}_{ij}|\geq \alpha_N t\Big)
& \leq & P\Big(12\cdot \max_{1\leq i \leq N}\delta_i^2\geq \alpha_N t
\Big)
\nonumber\\
& \leq & P\Big(\max_{1\leq i \leq N}\delta_i\geq \Big(\frac{\alpha_N t}{12}\Big)^{1/2}\Big).
\end{eqnarray*}
Therefore, to show $\sqrt{T\log N}\cdot\max_{1\leq i< j \leq N}|\rho_{ij}-\tilde{\rho}_{ij}| \to 0$, it is enough to prove that
\begin{eqnarray*}
P\Big(\max_{1\leq i \leq N}\delta_i> s\sqrt{\alpha_N}\Big)\to 0
\end{eqnarray*}
for any $s>0.$ In fact,
\begin{eqnarray*}
 & & P\Big(\max_{1\leq i \leq N}\delta_i> s\sqrt{\alpha_N}\Big)\\
& \leq & N\cdot P\Big(\frac{|\xi_1+\cdots +\xi_T|}{\sqrt{\xi_1^2+\cdots +\xi_T^2}}> s\sqrt{T\alpha_N}\Big)\\
& \leq & N\cdot P\Big(\xi_1^2+\cdots +\xi_T^2\leq \frac{1}{2}T\Big) + N\cdot  P\Big(|\xi_1+\cdots +\xi_T|> \frac{1}{\sqrt{2}}sT\sqrt{\alpha_N}\Big),
\end{eqnarray*}
where $\{\xi_j;\, 1\leq j \leq T\}$ are i.i.d. random variables with the same distribution of $\bd{\epsilon}_{11}$. The reason we switch the notations from $\{\bd{\epsilon}_{ij}\}$'s to $\{\xi_j;\, 1\leq j \leq T\}$ is for the brevity of symbols. By \eqref{hahaha},
\begin{eqnarray*}
 P\Big(\xi_1^2+\cdots +\xi_T^2\leq \frac{1}{2}T\Big)=O\Big(\frac{1}{T^{\tau/4}}\Big).
\end{eqnarray*}
By the Markov inequality and \eqref{hot_sleep2} as used in  \eqref{hahaha},
\begin{eqnarray*}
P\Big(|\xi_1+\cdots+ \xi_T|> \frac{1}{\sqrt{2}}sT\sqrt{\alpha_N}\Big)=O\Big( \frac{T^{\alpha/2}}{(T\sqrt{\alpha_N})^{\alpha}}\Big)=O\Big(\frac{(\log N)^{\alpha/4}}{T^{\alpha/4}}\Big)
\end{eqnarray*}
since $\alpha_N=1/\sqrt{T\log N}$. Combing the above assertions, we arrive at
\begin{eqnarray*}
P\Big(\max_{1\leq i \leq N}\delta_i> s\sqrt{\alpha_N}\Big)=O\Big(\frac{N}{T^{\alpha/4}}\Big) + O\Big(\frac{N(\log N)^{\alpha/4}}{T^{\alpha/4}}\Big),
\end{eqnarray*}
which converges to zero provided $\frac{N(\log N)^{\alpha/4}}{T^{\alpha/4}}\to 0$, or equivalently,  $\frac{N^{4/\alpha}\log N}{T}\to 0$     \hfill$\Box$

\subsubsection{Finale: proofs of Theorems \ref{linear_case}, \ref{exponential_case} and \ref{super_exponential_case}}
\label{PS2}

With preparations earlier, we are now ready to prove the main theorems on the maximum statistics of sample correlation coefficients.

\smallskip

\noindent\textbf{Proof of Theorem \ref{linear_case}}. Under the condition $E|\epsilon_{11}|^{6}<\infty$, \cite{Jiang2004} and \cite{Zhou2007} show that
\begin{align}\lbl{0407}
TL_N'^2-4\log N +\log\log N
\end{align}
converges weakly to a distribution with distribution function $F(y)$, where $L_N'=\max_{1\leq i< j \leq N}|\rho_{ij}|$ and $\rho_{ij}$ is as in \eqref{Pearson}. Set $L_N''=\max_{1\leq i< j \leq N}|\tilde{\rho}_{ij}|$ and $\tilde{\rho}_{ij}$ is as in Lemma \ref{Seattle}. Observe that
\begin{align}\lbl{WU}
|L_N-L_N'|
 \leq & |L_N-L_N''| + |L_N''-L_N'|\nonumber\\
 \leq & \max_{1\leq i< j \leq N}|\hat{\rho}_{ij}-\tilde{\rho}_{ij}| + \max_{1\leq i< j \leq N}|\tilde{\rho}_{ij}-\rho_{ij}|.
\end{align}
Since $E|\epsilon_{11}|^{\tau}<\infty$ with $\tau>8$, by using the assumption $T/N\to c\in (0, \infty)$, we see that  $\lim_{N\to\infty}T/(N^{8/\tau}\log N)=\infty$. Hence, by the Proposition \ref{Fall}
\begin{eqnarray*}
\sqrt{T\log N}\cdot\max_{1\leq i <j \leq N}\big|\hat{\rho}_{ij} - \tilde{\rho}_{ij}\big|\to 0
\end{eqnarray*}
in probability as $N\to\infty$.
On the other hand,  since $T/N\to c\in (0, \infty)$ and $\tau > 8$ we have $\frac{N^{4/\alpha}\log N}{T}\to 0$. Then, by  Proposition \ref{BJ},
\begin{eqnarray*}
\sqrt{T\log N}\cdot \max_{1\leq i <j \leq N}\big|\tilde{\rho}_{ij}-\rho_{ij}\big|\to 0
\end{eqnarray*}
in probability as $N\to\infty$. From \eqref{WU} we see that
\begin{align}\lbl{minimum}
\sqrt{T\log N}\cdot(L_N-L_N')\to 0
\end{align}
in probability. Set $\Delta=L_N-L_N'$. Then
\begin{align}\lbl{copy}
TL_N^2=T(L_N'+\Delta)^2=TL_N'^2 + 2TL_N'\Delta +T\Delta^2.
\end{align}
The Slutsky lemma and \eqref{0407} say that $(T/\log N)^{1/2}L_N'\to 2$ in probability. Consequently,
\begin{eqnarray*}
&& TL_N'\Delta = \Big(\frac{T}{\log N}\Big)^{1/2}L_N'\cdot \big(\sqrt{T\log N}\,\Delta\big)\to 0\\
&& T\Delta^2=\big[\sqrt{T\log N}\,\Delta\big]^2\cdot \frac{1}{\log N} \to 0
\end{eqnarray*}
in probability by \eqref{minimum}. These together with \eqref{copy} conclude that
\begin{align}\lbl{bioStat}
TL_N^2=TL_N'^2+o_p(1).
\end{align}
By the Slutsky lemma again, this fact and \eqref{0407} imply the desired result. \hfill$\Box$

\smallskip

\noindent\textbf{Proof of Theorem \ref{exponential_case}}. By assumption,  $Ee^{\omega|\epsilon_{11}|}<\infty$ and $\log N=o(T^{1/5})$. Using Theorem 3 and Remark 2.1 from \cite{Cai_J2011} with  $``\mu=0"$ and $``\alpha=1"$, we get
\begin{align}\lbl{2011}
T(L_N'')^2-4\log N +\log\log N
\end{align}
converges weakly to a distribution with distribution function $F(y)$, where $L_N''=\max_{1\leq i< j \leq N}|\tilde{\rho}_{ij}|$ and $\tilde{\rho}_{ij}$ is as in Lemma \ref{Seattle}. By Proposition \ref{Winter},
\begin{eqnarray*}
\sqrt{T\log N}\cdot\max_{1\leq i <j \leq N}\big|\hat{\rho}_{ij} - \tilde{\rho}_{ij}\big|\to 0
\end{eqnarray*}
in probability as $N\to\infty.$ Recall $L_N=\max_{1\leq i< j \leq N}|\hat{\rho}_{ij}|$. By the triangle inequality, the above says that
\begin{eqnarray*}
\sqrt{T\log N}\cdot (L_N-L_N'') \to 0
\end{eqnarray*}
in probability as $N\to\infty.$ Repeating the argument from \eqref{minimum} to \eqref{bioStat}, we obtain
\begin{eqnarray*}
TL_N^2=T(L_N'')^2+o_p(1)
\end{eqnarray*}
as $N\to\infty.$ The conclusion follows from \eqref{2011}.  \hfill$\Box$

\smallskip

\noindent\textbf{Proof of Theorem \ref{super_exponential_case}}. By assumption, $\log N=o(T^{1/3})$. Taking $``\mu=0"$ and $``\alpha=2"$ in  Theorem 3 and Remark 2.1 from
\cite{Cai_J2011}, we have
\begin{align}\lbl{ccpr}
T(L_N'')^2-4\log N +\log\log N
\end{align}
converges weakly to distribution function $F(y)$ for $y \in \mathbb{R}$, where $L_N''=\max_{1\leq i< j \leq N}|\tilde{\rho}_{ij}|$ and $\tilde{\rho}_{ij}$ is as in Lemma \ref{Seattle}. Recall $L_N=\max_{1\leq i< j \leq N}|\hat{\rho}_{ij}|$. By Proposition \ref{huns}, under the restriction $\log N=o(T^{1/3})$,
\begin{eqnarray*}
\sqrt{T\log N}\cdot\max_{1\leq i <j \leq N}\big|\hat{\rho}_{ij} - \tilde{\rho}_{ij}\big|\to 0
\end{eqnarray*}
in probability as $N\to\infty.$  By the triangle inequality, the above says that
\begin{eqnarray*}
\sqrt{T\log N}\cdot (L_N-L_N'') \to 0
\end{eqnarray*}
in probability as $N\to\infty.$ From the argument between \eqref{minimum} and  \eqref{bioStat}, we have
\begin{align}\label{order_1}
TL_N^2=T(L_N'')^2+o_p(1)
\end{align}
as $N\to\infty.$ This and \eqref{ccpr} yield the conclusion. \hfill$\Box$

\subsection{The proof of Theorem \ref{Asym_indept}}
\label{PS5}

We create a new method to prove Theorem \ref{Asym_indept} which gives the asymptotic independence between the sum $S_N$ and the maximum $L_N$. The idea is employing the inclusion-exclusion formula twice. We expect this method to work for other problems regarding asymptotic independence between sums of and maxima of weakly dependent random variables.

\subsubsection{Prelude: auxiliary results towards proof of Theorem \ref{Asym_indept}}\lbl{Pre_Asym}

The results stated in this section are about the estimates of probabilities of events related to Gaussian random variables. They are useful in their own right. Their proofs will be presented in Section \ref{Pre_Asym_last}.

\begin{lemma}\lbl{Lili} For each $N\geq 1$, let $T=T_N\geq 2$ be an integer. Suppose $\bd{s}_1$ and $\bd{s}_2$ are i.i.d. random vectors uniformly distributed on $\mathbb{S}^{T-1}$. Given $y \in \mathbb{R}$, set $l_N=T^{-1/2}\cdot (4\log N -\log\log N+y)^{1/2}$ which makes sense for large $N$. Assume $\log N=o(\sqrt{T})$ as $N\to\infty.$ Then
\begin{eqnarray*}
\lim_{N\to\infty}N^2\cdot P(\bd{s}_1'\bd{s}_2\geq l_N)=  \frac{1}{2\sqrt{2\pi}}e^{-y/2}.
\end{eqnarray*}
\end{lemma}

\begin{lemma}\lbl{Gaussian_inter} Suppose $\bd{s}_1$ and $\bd{s}_2$ are two i.i.d. random vectors uniformly distributed on $\mathbb{S}^{T-1}$ with $T\geq 2.$ Let $\{\xi_1, \cdots, \xi_k\}$ be random variables (not necessarily independent), each of which has the same distribution as that of  $\bd{s}_1'\bd{s}_2$.  Then
\begin{eqnarray*}
P\big(\max_{1\leq i \leq k}|\xi_i|\geq t\big) \leq k\cdot e^{-Tt^2/4}+(2k)\cdot e^{-cT}
\end{eqnarray*}
for all $t>\frac{2}{\sqrt{T}}$, where $c>0$ is a constant free of $k$, $t$ and $T$.
\end{lemma}

\begin{lemma}\lbl{hua_dan} Let  $\{Z, Z_1, \cdots, Z_k\}$ be i.i.d. standard normals. Let  $\delta \in (0, 1)$ be given. Set $v_i=\sqrt{\delta}Z+\sqrt{1-\delta}Z_i$ for $1\leq i \leq k$. Then
\begin{eqnarray*}
P\big(\min_{1\leq i \leq k}v_i>x\big)\leq \frac{1}{y}\exp\Big(-\frac{y^2}{2\delta}\Big)+ \frac{1}{(x-y)^k}\cdot \exp\Big[-\frac{k(x-y)^2}{2(1-\delta)}\Big]
\end{eqnarray*}
for all $x>y>0$.
\end{lemma}

\smallskip


\begin{lemma}\lbl{Kahane_ineq}(Slepian's lemma from \cite{Slepian}) Suppose $(U_1, \cdots, U_k)'$ and $(V_1, \cdots, V_k)'$ are two $\mathbb{R}^k$-valued centered Gaussian random vectors such that $EU_i^2=EV_i^2$ and $E(U_iU_j)\leq E(V_iV_j)$
for all $1\leq i, j \leq k.$ Then, for any real numbers $t_1, \cdots, t_k$,
\begin{eqnarray*}
P(U_i\leq t_i\ \mbox{for all}\ 1\leq i \leq k) \leq P(V_i\leq t_i\,\mbox{for all}\ 1\leq i \leq k).
\end{eqnarray*}
\end{lemma}

\begin{lemma}\lbl{Zheng_zhou} Suppose $\bd{a}_1, \cdots, \bd{a}_k$ are constant unit vectors on $\mathbb{S}^{T-1}$ for some $T\geq 2$. Let $\bd{s}$ be a vector with the uniform distribution on $\mathbb{S}^{T-1}.$ Assume $\max_{1\leq i< j \leq k}|\bd{a}_i'\bd{a}_j|\leq \delta$ for some $\delta\in [0, 1)$. Then
\begin{eqnarray*}
P\big(\min_{1\leq i \leq k}|\bd{a}_i'\bd{s}|>z\big)  &\leq& \frac{2^k}{y}\cdot\exp\Big(-\frac{y^2}{2\delta}\Big)+2\exp(-cT)\\
& &+ \frac{2^k}{(z\sqrt{rT}-y)^k}\cdot \exp\Big[-\frac{k\big(z\sqrt{rT}-y\big)^2}{2(1-\delta)}\Big]
\end{eqnarray*}
for all $z>0$,  $y \in (0, z\sqrt{rT})$, $r\in (0, 1)$ and $c$ is a constant depending on $r$  only.
\end{lemma}

\begin{lemma}\lbl{wen_fei} Let $\hat{\rho}_{ij}$ be as in \eqref{Asia_foods}.  Suppose assumption \eqref{assa} holds with $\{\epsilon_{ij};\, 1\leq i\leq N, 1\leq j \leq  T\}$ being Gaussian random variables. Recall $\bd{U}_i$ and $\bd{s}_i$ from Lemma \ref{long_for}. Let $\tau\geq 2$ be given. Then
\begin{eqnarray*}
E\big(\big|\hat{\rho}_{ij}^2-E(\hat{\rho}_{ij}^2|\bd{s}_i)\big|^{\tau}\big) \leq  \frac{K}{m^{\tau}},
\end{eqnarray*}
\begin{eqnarray*}
E\Big[\,\Big|\sum_{j=i+1}^N (\hat{\rho}_{ij}^2-E\hat{\rho}_{ij}^2)\Big|^{\tau}\Big] \leq
K\cdot \Big[\frac{(N-i)^{\tau/2}}{m^{\tau}}+ \frac{(N-i)^{\tau}}{m^{2\tau}}\Big]
\end{eqnarray*}
and
\begin{eqnarray*}
E\Big[\,\Big|\sum_{i=1}^{j-1} (\hat{\rho}_{ij}^2-E\hat{\rho}_{ij}^2)\Big|^{\tau}\Big] \leq
K\cdot \Big[\frac{(j-1)^{\tau/2}}{m^{\tau}} + \frac{(j-1)^{\tau}}{m^{2\tau}}\Big]
\end{eqnarray*}
for all $1\leq i < j\leq N$  and $N\geq 3$, where $K$ is a constant depending on $\tau$ only.
\end{lemma}

\subsubsection{Intermezzo: key steps in the proof of  Theorem \ref{Asym_indept}}\lbl{Key_asym_indept}

After collecting some useful facts in Section \ref{Pre_Asym}, we are now ready to prove Theorem \ref{Asym_indept}. To make the discussion easier to follow, we give the outline first.

First, Let $S_N$,  $L_N$ and $\mu_N$ be as in \eqref{sumtest}, \eqref{Write} and \eqref{taiyangtaila}, respectively. Review the framework between \eqref{park_1} and \eqref{Write}. In particular,
\begin{align}
\hat{\rho}_{ij}
=\frac{\bd{\epsilon}_i'\bd{P}_i\bd{P}_j\bd{\epsilon}_j}
{\|\bd{P}_i\bd{\epsilon}_i\|\cdot\|\bd{P}_j\bd{\epsilon}_j\|}\label{Asia_foods_66}
\end{align}
for $1\leq i, j \leq N.$ Assume \eqref{assa} holds with $\bd{\epsilon}_{i} \sim N_T(\bd{0}, \sigma_i^2\bd{I})$ for each $i$. Then
\begin{align}\lbl{rain_dear_lost}
\bd{e}_i:=\frac{\bd{\epsilon}_i}{\|\bd{\epsilon}_i\|},\ 1\leq i \leq N,
\end{align}
are i.i.d. uniformly distributed over the $T$-dimensional unit sphere $\mathbb{S}^{T-1}.$
For fixed $y\in \mathbb{R}$, set
\begin{align}\lbl{wet_soil}
l_N=T^{-1/2}\cdot (4\log N -\log\log N+y)^{1/2}.
\end{align}

Here is the structure of the proof of  Theorem \ref{Asym_indept}.

1. Let $\bd{e}_i$ be as in \eqref{rain_dear_lost}.
Define $\tilde{L}_N=\max_{1\leq i<j \leq N}|\bd{e}_i'\bd{e}_j|$.
To show that $TL_N^2-4\log N +\log\log N$ and $(S_N-\mu_N)/N$ are asymptotically independent, it is enough to prove that $T\tilde{L}_N^2-4\log N +\log\log N$ and $(S_N-\mu_N)/N$ are asymptotically independent (Lemma \ref{change_to_and}). The benefit of this step is that $\{\bd{e}_i'\bd{e}_j;\ 1\leq i<j\leq  N\}$ are identically distributed. This is not true for $\{\hat{\rho}_{ij};\, 1\leq i<j \leq N\}$ appeared in definition of $L_N$.

2. Review \eqref{wet_soil}. To show the asymptotic independence, it suffices to prove
\begin{align}\lbl{cuba}
\lim_{N\to \infty}P\Big(\frac{1}{N}(S_N-\mu_N)\leq x,\ \tilde{L}_N>l_N\Big)= \Phi(x)\cdot [1-F(y)]
\end{align}
for any real numbers $x$ and $y$,
where $\Phi(x)$ is the cdf of $N(0, 1)$ and is also the limiting distribution function  of $\frac{1}{N}(S_N-\mu_N)$; $F(y)$ is the Gumbel distribution and is also the  limiting distribution function  of $\tilde{L}_N$. Recall the definition of $\tilde{L}_N$, we are able to write the event in \eqref{cuba} as the union of $\binom{N}{2}$ many events which are exchangeable. Then, by using the inclusion-exclusion formula, the probability in  \eqref{cuba} is sandwiched between two bounds  [\eqref{Upper_bound} and \eqref{Lower_bound}]. The advantage is that we reduce the probability on the global maximum ``$\tilde{L}_N$" to sums of probabilities on ``local maxima".

3. In dealing with the ``local maxima", each probability in the sum is of the form $P(\frac{1}{N}(S_N-\mu_N)\leq x, |\bd{e}_{i_1}'\bd{e}_{j_1}|>l_N, \cdots, |\bd{e}_{i_n}'\bd{e}_{j_n}|>l_N )$, where $n$ is a fixed number free of $N$ and $T$, and where the indices $\{(i_l, j_l);\, 1\leq l \leq n\}$ are different. Review $S_N$  is  the sum of $(\bd{e}_i'\bd{e}_j)^2$ over all $1\leq i<j \leq N$. Remove the terms related to $\{\bd{e}_{i_l}'\bd{e}_{j_l};\, 1\leq l \leq n\}$ from $S_N$, in other words, eliminate the terms $(\bd{e}_i'\bd{e}_j)^2$ for all $(i, j)$ with $\{i, j\}\cap \{i_l, j_l\} \ne \emptyset$ for some $1\leq l\leq n$. Then the resulting sum is independent of $\{\bd{e}_{i_l}'\bd{e}_{j_l};\, 1\leq l \leq n\}$, and hence $P(\frac{1}{N}(S_N-\mu_N)\leq x, |\bd{e}_{i_1}'\bd{e}_{j_1}|>l_N, \cdots, |\bd{e}_{i_n}'\bd{e}_{j_n}|>l_N )$ is asymptotically the product of  $P(\frac{1}{N}(S_N-\mu_N)\leq x)$ and  $P(|\bd{e}_{i_1}'\bd{e}_{j_1}|>l_N, \cdots, |\bd{e}_{i_n}'\bd{e}_{j_n}|>l_N )$. Of course we have to handle the ``loss" after removing the terms. It turns out that the removed terms are very concentrated at their mean values by the second and third conclusions from Lemma \ref{wen_fei}. So the probability $P(\frac{1}{N}(S_N-\mu_N)\leq x)$ and the modified version $P(\frac{1}{N}(\tilde{S}_N-\tilde{\mu}_N)\leq x)$ are asymptotically equal. The total errors in the above approximations is negligible (Lemma \ref{akxbcvgf}).

4. In step 3, we have showed that
\begin{center}
$P(\frac{1}{N}(S_N-\mu_N)\leq x, |\bd{e}_{i_1}'\bd{e}_{j_1}|>l_N, \cdots, |\bd{e}_{i_n}'\bd{e}_{j_n}|>l_N )$
\end{center}
is asymptotically the product of  $P(\frac{1}{N}(S_N-\mu_N)\leq x)$ and  $P(|\bd{e}_{i_1}'\bd{e}_{j_1}|>l_N, \cdots, |\bd{e}_{i_n}'\bd{e}_{j_n}|>l_N )$ in \eqref{Upper_bound} and \eqref{Lower_bound}, where $A_N=\{\frac{1}{N}(S_N-\mu_N)\leq x\}$ and  $B_{I}=\{|\bd{e}_i'\bd{e}_j|>l_N\}$ for $I=(i, j)$. We will use one more time the inclusion-exclusion formula to regroup the sum of probabilities  $P(|\bd{e}_{i_1}'\bd{e}_{j_1}|>l_N, \cdots, |\bd{e}_{i_n}'\bd{e}_{j_n}|>l_N )$ and change it to $P(\max_{1\leq i<j \leq N}|\bd{e}_{i}'\bd{e}_{j}|>l_N)$. Note that the original upper bound becomes the lower bound  of $P(\max_{1\leq i<j \leq N}|\bd{e}_{i}'\bd{e}_{j}|>l_N)$ and similarly the original lower bound becomes the new upper-bound. There are some ``middle" terms in between the bounds, we have to show they are negligible. This is guaranteed by Lemma \ref{one_key}.

Now let us execute the steps streamlined above.

\begin{lemma}\lbl{change_to_and} Let $S_N$,  $L_N$ and $\mu_N$ be as in \eqref{sumtest}, \eqref{Write} and \eqref{taiyangtaila}, respectively. Let $\{\bd{e}_i;\ 1\leq i \leq N\}$ be defined in \eqref{rain_dear_lost}.  Set $\tilde{L}_N=\max_{1\leq i<j \leq N}|\bd{e}_i'\bd{e}_j|$. Assume $N=o(T^2)$ and \eqref{assa} holds with $\bd{\epsilon}_{i} \sim N_T(\bd{0}, \sigma_i^2\bd{I})$ for each $i$. If $T\tilde{L}_N^2-4\log N +\log\log N$ and $\frac{1}{N}(S_N-\mu_N)$ are asymptotically independent, then $TL_N^2-4\log N +\log\log N$ and $\frac{1}{N}(S_N-\mu_N)$ are also asymptotically independent.
\end{lemma}
\noindent\textbf{Proof of Lemma \ref{change_to_and}}. Let $m=T-p.$ Under assumption \eqref{assa} with $\bd{\epsilon}_{i} \sim N_T(\bd{0}, \sigma_i^2\bd{I})$ for each $i$, we know  $\{\bd{e}_i;\ 1\leq i \leq N\}$ are i.i.d. uniformly distributed over $\mathbb{S}^{T-1}.$
Define $\tilde{\rho}_{ij}=\bd{e}_i'\bd{e}_j$ for $1\leq i<j \leq N.$ To organize the proof clearly, we list the relevant quantities as follows.
\begin{eqnarray*}
&& L_N=\max_{1\leq i< j \leq N}|\hat{\rho}_{ij}|\ \ \ \ \mbox{and}\ \ \ \ \tilde{L}_N=\max_{1\leq i<j \leq N}|\tilde{\rho}_{ij}|;\\
&& S_N=\sum_{1\leq i< j \leq N}T\hat{\rho}_{ij}^2\ \ \ \mbox{and}\ \ \  \mu_N=\frac{T}{m^2}\sum_{j=2}^N\sum_{i=1}^{j-1}\mbox{tr}(\bd{P}_i\bd{P}_j)
\end{eqnarray*}
where $\bd{P}_i$ is defined  in \eqref{park_1}.
By Theorems \ref{super_exponential_case} and \ref{diligent}, the following hold.
\begin{align}
&TL_N^2-4\log N +\log\log N \to F(y)=\exp(-(1/\sqrt{8\pi})e^{-y/2})\ \mbox{weakly}; \lbl{jian_jiaoooo}\\
&\frac{1}{N}(S_N-\mu_N)\to N(0, 1)\ \mbox{weakly}.\ \ \ \ \ \ \lbl{gan_doufuuuu}
\end{align}
By Theorem 6 from \cite{CJF13},  the assertion \eqref{jian_jiaoooo} is also true if ``$L_N$" is replaced by ``$\tilde{L}_N$".
To show asymptotic independence, it is enough to show
\begin{align}\lbl{qwkucb}
\lim_{N\to \infty}P\Big(\frac{1}{N}(S_N-\mu_N)\leq x,\ TL_N^2-4\log N +\log\log N\leq y\Big)= \Phi(x)\cdot F(y)
\end{align}
for any $x\in \mathbb{R}$ and $y \in \mathbb{R}$, where $\Phi(x)=(2\pi)^{-1/2}\int_{-\infty}^xe^{-t^2/2}\,dt.$ Let $l_N$ be as in \eqref{wet_soil}.
Due to \eqref{jian_jiaoooo} and \eqref{gan_doufuuuu} we know  \eqref{qwkucb} is equivalent to that
\begin{align}\lbl{wealth_noo}
\lim_{N\to \infty}P\Big(\frac{1}{N}(S_N-\mu_N)\leq x,\ L_N>l_N\Big)= \Phi(x)\cdot [1-F(y)]
\end{align}
for any $x\in \mathbb{R}$ and $y \in \mathbb{R}$. By assumption, we know that
\begin{align}\lbl{dslkfnuec}
\lim_{N\to \infty}P\Big(\frac{1}{N}(S_N-\mu_N)\leq x,\ \tilde{L}_N>l_N\Big)= \Phi(x)\cdot [1-F(y)]
\end{align}
for any $x\in \mathbb{R}$ and $y \in \mathbb{R}$. We show next that \eqref{dslkfnuec} implies \eqref{wealth_noo}.

By Proposition \ref{huns},
\begin{eqnarray*}
\sqrt{T\log N}\cdot\max_{1\leq i <j \leq N}\big|\hat{\rho}_{ij} - \tilde{\rho}_{ij}\big|\to 0
\end{eqnarray*}
in probability as $N\to\infty$ provided $\log N=o(T^{1/3})$ as $N\to\infty$. By the triangle inequality,
\begin{eqnarray*}
\sqrt{T\log N}\cdot\big|L_N - \tilde{L}_N\big|\to 0
\end{eqnarray*}
in probability.

Given $\epsilon \in (0,1)$. Set
\begin{eqnarray*}
\Omega_N=\big\{\sqrt{T\log N}\cdot\big|L_N - \tilde{L}_N\big|<\epsilon\big\}
\end{eqnarray*}
for $N\geq 3.$ Then
\begin{align}\lbl{ba_ma}
\lim_{N\to \infty}P(\Omega_N)=1.
\end{align}
Now,
\begin{align}
& P\Big(\frac{1}{N}(S_N-\mu_N)\leq x,\ L_N>l_N\Big)\nonumber\\
\leq & P\Big(\frac{1}{N}(S_N-\mu_N)\leq x,\ L_N>l_N,\ \Omega_N\Big)+ P(\Omega_N^c).\lbl{2jhr3}
\end{align}
On $\Omega_N$, if $L_N>l_N$ then
\begin{align}\lbl{jsalfs}
\tilde{L}_N\geq L_N-\big|L_N - \tilde{L}_N\big|>l_N- \frac{\epsilon}{\sqrt{T\log N}}.
\end{align}
Define
\begin{eqnarray*}
\tilde{l}_N=T^{-1/2}\cdot (4\log N -\log\log N+y-5\epsilon)^{1/2},
\end{eqnarray*}
which makes sense for large $N.$
Use the formula $\sqrt{x}-\sqrt{y}=(x-y)/(\sqrt{x}+\sqrt{y})$ for any $x\geq 0$ and $y\geq 0$ to see
\begin{align}
& T^{1/2}\cdot(l_N-\tilde{l}_N) \nonumber\\
=& (4\log N -\log\log N+y)^{1/2}-(4\log N -\log\log N+y-5\epsilon)^{1/2} \nonumber\\
 = & \frac{5\epsilon}{(4\log N -\log\log N+y)^{1/2}+(4\log N -\log\log N+y-5\epsilon)^{1/2}} \nonumber\\
 \sim & \frac{5\epsilon}{4\sqrt{\log N}} \lbl{sdusicb}
\end{align}
as $N \to\infty.$
Thus,
\begin{eqnarray*}
l_N-\tilde{l}_N > \frac{\epsilon}{\sqrt{T\log N}}.
\end{eqnarray*}
as $N$ is sufficiently large. This and \eqref{jsalfs} conclude that
\begin{eqnarray*}
\tilde{L}_N\geq \tilde{l}_N
\end{eqnarray*}
as $N$ is sufficiently large. Review \eqref{2jhr3}. We have
\begin{eqnarray*}
&& P\Big(\frac{1}{N}(S_N-\mu_N)\leq x,\ L_N>l_N\Big)\nonumber\\
&\leq & P\Big(\frac{1}{N}(S_N-\mu_N)\leq x,\ \tilde{L}_N\geq \tilde{l}_N\Big)+ P(\Omega_N^c).
\end{eqnarray*}
Immediately from \eqref{dslkfnuec} and \eqref{ba_ma} we get
\begin{eqnarray*}
\limsup_{N\to\infty}P\Big(\frac{1}{N}(S_N-\mu_N)\leq x,\ L_N>l_N\Big)\leq \Phi(x)\cdot [1-F(y-5\epsilon)]
\end{eqnarray*}
for any $\epsilon \in (0,1).$ Inspect that the left-hand side of the above  does not depend on $\epsilon$. Letting $\epsilon \downarrow 0$, we obtain
\begin{align}\lbl{IXBCATRS}
\limsup_{N\to\infty}P\Big(\frac{1}{N}(S_N-\mu_N)\leq x,\ L_N>l_N\Big)\leq \Phi(x)\cdot [1-F(y)]
\end{align}
for any $x\in \mathbb{R}$ and $y \in \mathbb{R}.$ In the following we will show the lower limit.

Evidently,
\begin{align}\lbl{84659}
& P\Big(\frac{1}{N}(S_N-\mu_N)\leq x,\ L_N>l_N\Big)\nonumber\\
\geq & P\Big(\frac{1}{N}(S_N-\mu_N)\leq x,\ L_N>l_N, \Omega_N\Big).
\end{align}
Set
\begin{eqnarray*}
\tilde{l}_N'=T^{-1/2}\cdot (4\log N -\log\log N+y+5\epsilon)^{1/2}.
\end{eqnarray*}
Similar to \eqref{sdusicb}, it is checked that
\begin{eqnarray*}
T^{1/2}\cdot (\tilde{l}_N'-l_N) \sim \frac{5\epsilon}{4\sqrt{\log N}}
\end{eqnarray*}
as $N\to\infty$. Therefore,
\begin{eqnarray*}
\tilde{l}_N'>l_N+\frac{\epsilon}{\sqrt{T\log N}}
\end{eqnarray*}
as $N$ is sufficiently large. It is straightforward to verify that
\begin{eqnarray*}
\big\{\tilde{L}_N>\tilde{l}_N',\ \Omega_N\big\} \subset \big\{\tilde{L}_N>l_N+\frac{\epsilon}{\sqrt{T\log N}},\ \Omega_N\big\} \subset   \big\{L_N>l_N,\ \Omega_N\big\}
\end{eqnarray*}
as $N$ is sufficiently large, where the last inclusion follows from the definition of $\Omega_N$. By \eqref{84659},
\begin{eqnarray*}
&& P\Big(\frac{1}{N}(S_N-\mu_N)\leq x,\ L_N>l_N\Big)\\
&\geq & P\Big(\frac{1}{N}(S_N-\mu_N)\leq x,\, \tilde{L}_N>\tilde{l}_N',\, \Omega_N \Big).
\end{eqnarray*}
Thus, from \eqref{dslkfnuec} and \eqref{ba_ma} we get
\begin{eqnarray*}
\liminf_{N\to\infty}P\Big(\frac{1}{N}(S_N-\mu_N)\leq x,\ L_N>l_N\Big)\geq \Phi(x)\cdot [1-F(y+5\epsilon)]
\end{eqnarray*}
for any $\epsilon \in (0,1).$ Sending $\epsilon \downarrow 0$ we see
\begin{eqnarray*}
\liminf_{N\to\infty}P\Big(\frac{1}{N}(S_N-\mu_N)\leq x,\ L_N>l_N\Big)\geq \Phi(x)\cdot [1-F(y)]
\end{eqnarray*}
for any $x\in \mathbb{R}$ and $y \in \mathbb{R}.$ This together with \eqref{IXBCATRS} concludes \eqref{wealth_noo}. \hfill$\Box$

\smallskip

We need some notations now. Let $S_N$,  $L_N$ and $\mu_N$ be as in \eqref{sumtest}, \eqref{Write} and \eqref{taiyangtaila}, respectively.
Let $\{\bd{e}_i;\ 1\leq i \leq N\}$ be as in \eqref{rain_dear_lost}. Define
\begin{align}
& \Lambda_N=\{(i, j);\, 1\leq i<j\leq N\}; \nonumber\\
& A_N=\Big\{\frac{1}{N}(S_N-\mu_N)\leq x\Big\}\ \ \ \mbox{and}\ \ \ B_{I}=\{|\bd{e}_i'\bd{e}_j|>l_N\} \lbl{haokanne}
\end{align}
for any $I=(i,j)\in \Lambda_N$. To make a clear presentation, we impose a trivial ordering for elements in  $\Lambda_N$. For any $I_1=(i_1, j_1)\in \Lambda_N$ and $I_2=(i_2, j_2) \in \Lambda_N$, we say $I_1<I_2$ if $i_1<i_2$ or $i_1=i_2$ but $j_1<j_2$.

\begin{lemma}\lbl{one_key} Recall the notations from \eqref{wet_soil} and \eqref{haokanne}. Assume $\log N=o(\sqrt{T})$ as $N\to\infty$. Assume  $\{\bd{e}_i;\ 1\leq i \leq N\}$ are i.i.d. uniformly distributed over $\mathbb{S}^{T-1},$ which is particularly true if \eqref{assa} holds with $\bd{\epsilon}_{i} \sim N_T(\bd{0}, \sigma_i^2\bd{I})$ for each $i$. Set
\begin{eqnarray*}
H(N, n)=\sum_{I_1< I_2< \cdots < I_{n}\in \Lambda_N}P(B_{I_1}B_{I_2}\cdots B_{I_{n}}).
\end{eqnarray*}
Then $\lim_{n\to\infty}\limsup_{N\to\infty}H(N, n)=0.$
\end{lemma}
\noindent\textbf{Proof of Lemma \ref{one_key}}. For $I_l$ appeared in $H(N, n)$, write $I_l=(i_l, j_l)$ for $l=1, \cdots, n.$ Now we classify the indices $I_1< I_2< \cdots < I_{n}\in \Lambda_N$ in the definition of $H(N, n)$ into three cases. Let $\Gamma_{N,1}$ be the set of indices $(I_1, \cdots, I_n)$ such that no two of the $2n$ indices $\{i_l, j_l\,; l=1, \cdots, n\}$ are identical. Let $\Gamma_{N,2}$ be the set of indices $(I_1, \cdots, I_n)$ such that either $i_1=\cdots =i_n$ or $j_1=\cdots = j_n.$ Let $\Gamma_{N,3}$ be the set of indices $I_1< I_2< \cdots < I_{n}\in \Lambda_N$ excluding $\Gamma_{N,1}\cup \Gamma_{N,2}$. In the following we will estimate
\begin{eqnarray*}
F_j:=\sum_{I_1< I_2< \cdots < I_{n}\in \Gamma_{N,j}}P(B_{I_1}B_{I_2}\cdots B_{I_{n}})
\end{eqnarray*}
for $j=1, 2, 3$ one by one. We will see $F_1$ contributes essentially the sum in the expression of $H(N, n)$ by an easy argument; the term $F_2$ is negligible and its computation is trivial; the term $F_3$ is also negligible but its estimate is most involved.

{\it Step 1: the estimate of $F_1$}. Recall $B_{I}=\{|\bd{e}_i'\bd{e}_j|>l_N\}$ if $I=(i,j)\in \Lambda_N$, where $l_N$ is defined in \eqref{wet_soil}. By the definition of $\Gamma_{N,1}$,  we know that $B_{I_1}, B_{I_2} \cdots, B_{I_n}$ are independent. By Lemma \ref{Lili} and the symmetry of $\bd{e}_1'\bd{e}_2$,
\begin{align}\lbl{kshjg}
\max_{I \in \Lambda_N}P(B_{I})=P(|\bd{e}_1'\bd{e}_2|\geq l_N)=2P(\bd{e}_1'\bd{e}_2\geq l_N)\leq \frac{C}{N^2}
\end{align}
for all $N\geq 3$. Then, by the elementary fact $\binom{k}{n}=\frac{1}{n!}k(k-1)\cdots (k-n+1)\leq \frac{k^n}{n!}$ for all $k> n\geq 1.$
\begin{align}\lbl{thankful}
 F_1\leq \frac{C}{N^{2n}}\cdot \binom{\frac{N(N-1)}{2}}{n} \leq \frac{C}{n!}.
\end{align}

{\it Step 2: the estimate of $F_2$}. Evidently, the size of $\Gamma_{N,2}$ is no more than $\binom{N}{1}\cdot\binom{N}{n}\cdot 2\leq 2N^{n+1}$. We first claim that $\{\bd{e}_1'\bd{e}_2, \bd{e}_1'\bd{e}_3, \cdots, \bd{e}_1'\bd{e}_n\}$ are independent. In fact, let $\bd{e}$ be uniformly distributed on $\mathbb{S}^{T-1}$. Then, $\bd{a}'\bd{e}$ has the same distribution as that of $(1, 0, \cdots, 0)'\bd{e}$ for any $a\in \mathbb{S}^{T-1}$ (see, e.g., Theorem 1.5.7(i) and the argument for (5) on p.147 from \cite{mui}). Since $\bd{e}_1, \cdots, \bd{e}_n$ are i.i.d. random vectors, we know that, conditioning on $\bd{e}_1$, the random variables $\{\bd{e}_1'\bd{e}_2, \bd{e}_1'\bd{e}_3, \cdots, \bd{e}_1'\bd{e}_n\}$ are i.i.d. with a common distribution of $(1, 0, \cdots, 0)'\bd{e}$. In particular, their conditional distributions do not depend on $\bd{e}_1$. This proves the claim. Consequently,
\begin{align}\lbl{nice_hao}
F_2
\leq & 2N^{n+1}\cdot P(|\bd{e}_1'\bd{e}_2|>l_N, \cdots, |\bd{e}_1'\bd{e}_{n+1}|>l_N) \nonumber\\
 = & 2N^{n+1}\cdot \big[P(|\bd{e}_1'\bd{e}_2|>l_N)\big]^n  \leq  \frac{2C^n}{N^{n-1}}
\end{align}
by \eqref{kshjg}.

{\it Step 3: the estimate of $F_3$}. Fix a tuple $(I_1, I_2, \cdots , I_{n})\in \Gamma_{N,3}$. By the ordering imposed on $\Lambda_N$, we see that $i_1\leq i_2\leq \cdots \leq i_n$. There are two different  cases: (1) $i_1<i_2$;
(2) there exists $2\leq k\leq n-1$ such that $i_1=\cdots =i_k<i_{k+1}$.

Under case (1), let $\ml{F}_1$ be the set of random vectors  $\{\bd{e}_{j_1}, \bd{e}_{i_l}, \bd{e}_{j_l};\, 2\leq l\leq n\}$ (the first index is ``$j_1$" which is different from the third one ``$j_l$"). Then, by independence and the property ``take out what is known" for the conditional probability,
\begin{eqnarray*}
P(B_{I_1}B_{I_2}\cdots B_{I_{n}})
&=&E[P(B_{I_1}B_{I_2}\cdots B_{I_{n}}|\ml{F}_1)]\nonumber\\
& = & E\Big[ {P\Big(|\bd{e}_{i_1}'\bd{e}_{j_1}|\geq l_N|\bd{e}_{j_1}\Big)}\cdot \prod_{l=2}^n{I(B_{I_{l}})}\Big].
\end{eqnarray*}
As a fact used earlier, the conditional distribution of $\bd{e}_{i_1}'\bd{e}_{j_1}$ given $\bd{e}_{j_1}$ and the unconditional distribution of $\bd{e}_{i_1}'\bd{e}_{j_1}$ are identical. Therefore, by \eqref{kshjg},
\begin{align}\lbl{boat_pink}
P(B_{I_1}B_{I_2}\cdots B_{I_{n}}) \leq \frac{C}{N^2}\cdot P(B_{I_2}\cdots B_{I_{n}}).
\end{align}

Let us study case (2).
Without loss of generality, for notational clarity, we assume $i_1=\cdots =i_k=1$ and $i_{k+1}=2$. Denote by $\ml{F}_2$ the set of random vectors  $\{\bd{e}_{i_l}, \bd{e}_{j_l};\, 1\leq l\leq n\}$ excluding $\bd{e}_1$. Then use conditional probability and independence to see
\begin{align}\lbl{196720}
P(B_{I_1}B_{I_2}\cdots B_{I_{n}})
=&E[P(B_{I_1}B_{I_2}\cdots B_{I_{n}}|\ml{F}_2)]\nonumber\\
 = & E\Big[ P_1\Big(\min_{1\leq l \leq k}|\bd{e}_1'\bd{e}_{j_l}|\geq {l_N}\Big)\cdot \prod_{l: i_l\ne 1}{I(B_{I_{l}})}\Big],
\end{align}
where $P_1$ stands for the condition probability given $\ml{F}_2$. By independence, the last probability in \eqref{196720} is computed by treating $\bd{e}_1$ as a random variable  while fixing the values of $\bd{e}_{j_1}, \cdots, \bd{e}_{j_{k}}$. To study the $P_1$, we need to understand  the relationship among $\{\bd{e}_{j_1}, \cdots, \bd{e}_{j_{k}}\}.$ To do so, set
\begin{eqnarray*}
\Omega_N=\Big\{\max_{j_1\leq l_1<l_2\leq j_{k}}|\bd{e}_{j_{l_1}}'\bd{e}_{j_{l_2}}|< \delta_N\Big\}
\ \ \ \mbox{and}\ \ \ \delta_N=4\sqrt{\frac{n\log N}{T}}.
\end{eqnarray*}
By Lemma \ref{Gaussian_inter} and the fact $k\leq n$,
\begin{align}\lbl{83838}
P(\Omega_N^c) \leq & 2k\cdot \Big[\exp\Big(-\frac{1}{4}T\delta_N^2\Big)+ e^{-cT}\Big]\nonumber\\
 \leq  &  \frac{3n}{N^{3n}}
\end{align}
as $N$ is sufficiently large provided $\log N=o(T)$.  Notice that
\begin{align}\lbl{0458282}
 P_1\Big(\min_{1\leq l \leq k}|\bd{e}_1'\bd{e}_{j_l}|\geq {l_N}\Big) \leq {I(\Omega_N^c)}+ P_1\Big(\min_{1\leq l \leq k}|\bd{e}_1'\bd{e}_{j_l}|\geq {l_N} \Big)\cdot {I(\Omega_N)}.
\end{align}
We claim that, for any $\epsilon\in (0,1)$, there exists an integer $N_{\epsilon}\geq 1$ such that
\begin{align}\lbl{tu_wu}
P_1\Big(\min_{1\leq l \leq k}|\bd{e}_1'\bd{e}_{j_l}|\geq {l_N} \Big)\cdot {I(\Omega_N)}
\leq \frac{1}{N^{2k-\epsilon}}
\end{align}
as $N\geq N_{\epsilon}$. On $\Omega_N$, we know $\max_{j_1\leq l_1<l_2\leq j_{k}}|\bd{e}_{j_{l_1}}'\bd{e}_{j_{l_2}}|< \delta_N.$ Take  $r=1-\frac{\epsilon}{4k}$, $y=(\log N)^{1/4}$, $z=l_N$, $\delta=\delta_N$ in Lemma \ref{Zheng_zhou}. Observe that $2^k/y\to 0$. By \eqref{wet_soil}, $l_N\sim  2\sqrt{(\log N)/T}$ as $N\to\infty$, and hence $y=o(z\sqrt{rT})$. Also, $\delta_N\to 0$ since $\log N=o(T)$.
\begin{eqnarray*}
\frac{k\big(z\sqrt{rT}-y\big)^2}{2(1-\delta)} \sim (2rk)\cdot\log N.
\end{eqnarray*}
Thus, by the lemma, use the facts that $2rk>2k-\frac{3}{4}\epsilon$ and that $z\sqrt{rT}-y\to \infty$ to get
\begin{eqnarray*}
&& P_1\Big(\min_{1\leq l \leq k}|\bd{e}_1'\bd{e}_{j_l}|\geq l_N \Big)\cdot I(\Omega_N)\\
&\leq & \exp\Big(-\frac{1}{2\delta_N}(\log N)^{1/2}\Big)+  \exp\Big[-\frac{k\big(z\sqrt{rT}-y\big)^2}{2(1-\delta)}\Big]+2\exp(-cT)\\
& \leq &\exp\Big(-\frac{1}{8\sqrt{n}}\sqrt{T}\Big) +\exp\Big[-\Big(2k-\frac{3}{4}\epsilon\Big)\cdot\log N\Big]+2\exp(-cT)\\
&\leq & \frac{1}{N^{2k-\epsilon}}+2\exp(-cT)
\end{eqnarray*}
as $N\geq N_{\epsilon}$ thanks to the assumption $\log N=o(\sqrt{T})$, where $N_{\epsilon}\geq 1$ is an integer depending on $\epsilon$ only. This leads to \eqref{tu_wu}.

Now, combining \eqref{0458282} and \eqref{tu_wu}, we arrive at
\begin{eqnarray*}
P_1\Big(\min_{2\leq l \leq k+1}|\bd{e}_1'\bd{e}_{j_l}|\geq l_N\Big) \leq I(\Omega_N^c)+ \frac{1}{N^{2k-\epsilon}}+2\exp(-cT)
\end{eqnarray*}
as $N\geq N_{\epsilon}$. This together with \eqref{196720} and \eqref{83838} implies
\begin{eqnarray*}
P(B_{I_1}B_{I_2}\cdots B_{I_{n}}) \leq \frac{1}{N^{2k-\epsilon}} \cdot  P\Big(\bigcap_{l: i_l\ne 1}B_l\Big) + \frac{3n}{N^{3n}}+2\exp(-cT)
\end{eqnarray*}
as $N$ is sufficiently large. In summary, by using the above conclusion and \eqref{boat_pink}, for any $\epsilon\in (0, 1)$ and any $(I_1, \cdots, I_n)\in \Gamma_{n,3}$,
\begin{eqnarray*}
P(B_{I_1}B_{I_2}\cdots B_{I_{n}}) \leq \frac{1}{N^{2k_1-\epsilon}} \cdot  P\Big(\bigcap_{l: i_l> i_1}B_{I_l}\Big) + \frac{3n}{N^{3n}}+2\exp(-cT)
\end{eqnarray*}
as $N\geq N_{\epsilon}$, where $k_1$ is the number of elements on the $i_1$-th row of $\{I_1, \cdots, I_n\}$. In words, when we consider $P(B_{I_1}B_{I_2}\cdots B_{I_{n}})$ based on the positions of $I_j$'s appeared in the upper triangular matrix $\Lambda_N=\{(i, j);\, 1\leq i<j \leq N\}$, after reducing the first row we see the connection between the old and new probabilities. Similarly, let $k_j$ be the number of elements from $\{I_1, \cdots, I_n\}$ on the $j$-th row for $j\geq 1.$ Then
\begin{eqnarray*}
P(B_{I_1}B_{I_2}\cdots B_{I_{n}}) &\leq & \frac{1}{N^{2k_1-\epsilon}} \cdot  \Big[\frac{1}{N^{2k_2-\epsilon}} \cdot  P\Big(\bigcap B_l\Big)\\
 & &+ \frac{3n}{N^{3n}}+2\exp(-cT)\Big] + \frac{3n}{N^{3n}}+2\exp(-cT)\\
& \leq & \frac{1}{N^{2k_1+2k_2-2\epsilon}} \cdot   P\Big(\bigcap B_l\Big)+ 2\cdot \frac{3n}{N^{3n}}+4\exp(-cT)
\end{eqnarray*}
where the two intersections above run over all elements from $\{I_1, \cdots, I_n\}$ excluding the first two rows. Continue the process recursively to see
\begin{eqnarray*}
P(B_{I_1}B_{I_2}\cdots B_{I_{n}}) \leq  \frac{1}{N^{2k_1+\cdots+2k_b-b\epsilon}} + b\cdot \frac{3n}{N^{3n}}+2b\exp(-cT)
\end{eqnarray*}
where $b$ is the total number of rows of $\{I_1, \cdots, I_n\}$ in the upper triangular matrix $\Lambda_N=\{(i, j);\, 1\leq i<j \leq N\}$. Obviously, $k_1+\cdots +k_b=n$ and $b\leq n.$ Therefore, for each $\epsilon \in (0, 1)$,
\begin{eqnarray*}
P(B_{I_1}B_{I_2}\cdots B_{I_{n}}) \leq  \frac{1}{N^{2n-n\epsilon}} +  \frac{3n^2}{N^{3n}}+2n\exp(-cT)
\end{eqnarray*}
for $N\geq N_{\epsilon}$. This gives that
\begin{align}\lbl{comarade}
F_3  = & \sum_{I_1< I_2< \cdots < I_{n}\in \Gamma_{N,3}}P(B_{I_1}B_{I_2}\cdots B_{I_{n}}) \nonumber\\
\leq & \ |\Gamma_{N,3}| \cdot \Big(\frac{1}{N^{2n-n\epsilon}} +  \frac{3n^2}{N^{3n}}+2ne^{-cT}\Big).
\end{align}
Recall $I_l=(i_l, j_l)$ for each $1\leq l \leq n.$ In view of the definition of $\Gamma_{N,3}$, there are at least two of the $2n$ indices from $\{(i_l, j_l);\, 1\leq l \leq n\}$ are identical for any $(I_1, \cdots, I_n)\in \Gamma_{N,3}$.  Let $\kappa=|\{i_l, j_l; 1\leq l \leq n\}|$ for such $(I_1, I_2, \cdots , I_{n})$. Easily, $n+1\leq \kappa\leq 2n-1$. To see how many such $(I_1, \cdots, I_n)$ with $|\{i_l, j_l; 1\leq l \leq N\}|=\kappa$, first pick $\kappa$ many indices from $\{1, 2, \cdots, N\}$, which has the total number of ways $\binom{N}{\kappa}\leq N^{\kappa}$, then use the $\kappa$ many indices to make a $(I_1, \cdots, I_n)\in \Gamma_{N,3}$. The total number of ways to do so is no more than $\kappa^{2n}$. Therefore,
\begin{eqnarray*}
|\Gamma_{N,3}|\leq \sum_{\kappa=n+1}^{2n-1}N^{\kappa}\cdot\kappa^{2n}\leq  (2n)^{2n}\cdot N^{2n-1}.
\end{eqnarray*}
As a consequence, for each $\epsilon \in (0, 1)$, from \eqref{comarade} we have
\begin{eqnarray*}
F_3\leq (2n)^{2n}\cdot \Big(\frac{1}{N^{1-n\epsilon}}+\frac{3n^2}{N^{n+1}}+2n\exp(-cT)\Big)
\end{eqnarray*}
as $N\geq N_{\epsilon}$. Take $\epsilon=\frac{1}{2n}$ to see $\lim_{N\to\infty}F_3= 0.$ Joining this with \eqref{thankful} and \eqref{nice_hao}, we eventually arrive at
\begin{align}\lbl{lncsivlbt}
\limsup_{N\to\infty}H(N, n)\leq \frac{C}{n!}
\end{align}
for each $n\geq 3.$ The desired conclusion then follows by sending $n\to \infty.$  \hfill$\Box$



\begin{lemma}\lbl{akxbcvgf} Recall the notations from \eqref{wet_soil} and \eqref{haokanne}.  Assume  \eqref{assa} holds with $\bd{\epsilon}_{i} \sim N_T(\bd{0}, \sigma_i^2\bd{I})$ for each $i$.  If $N=o(T^2)$ as $N\to \infty$, then
\begin{eqnarray*}
\sum_{I_1< I_2< \cdots < I_{n}\in \Lambda_N}\big[P(A_NB_{I_1}B_{I_2}\cdots B_{I_{n}}) - P(A_N)\cdot P(B_{I_1}B_{I_2}\cdots B_{I_{n}})\big]\to 0
\end{eqnarray*}
as $N\to\infty$ for each $n\geq 1.$
\end{lemma}
\noindent\textbf{Proof of Lemma \ref{akxbcvgf}}. From assumption that \eqref{assa} holds with $\bd{\epsilon}_{i} \sim N_T(\bd{0}, \sigma_i^2\bd{I})$ for each $i$, we know from \eqref{rain_dear_lost} that $\{\bd{e}_i;\ 1\leq i \leq N\}$ are i.i.d. uniformly distributed over $\mathbb{S}^{T-1}$. For $I_1< I_2< \cdots < I_{n}\in \Lambda_N$, write $I_l=(i_l, j_l)$ for $l=1,2,\cdots, n$.  Set
\begin{eqnarray*}
\Lambda_{n,N}=\big\{(i_l, j);\, i_l<j \leq N, 1\leq l \leq n\big\}
\bigcup\big\{(i, j_l);\, 1\leq i<j_l, 1\leq l \leq n\big\}
\end{eqnarray*}
for $n\geq 1.$ It is easy to check that $|\Lambda_{n,N}| =\sum_{l=1}^n(N-i_l+j_l-2)$. Since $i_l<j_l$ for each $l$, we see that
\begin{eqnarray*}
n(N-1)\leq |\Lambda_{n,N}| \leq \sum_{l=1}^n(N+j_l)\leq 2nN.
\end{eqnarray*}
 Recall $m=T-p$ and
\begin{eqnarray*}
S_N=\sum_{1\leq i< j \leq N}T\hat{\rho}_{ij}^2\ \ \ \mbox{and}\ \ \ \mu_N=ES_N=\frac{T}{m^2}\sum_{j=2}^N\sum_{i=1}^{j-1}\mbox{tr}(\bd{P}_i\bd{P}_j)
\end{eqnarray*}
where $\bd{P}_i$ is defined as in \eqref{park_1}. Define
\begin{eqnarray*}
A_N(x)=\Big\{\frac{1}{N}(S_N-\mu_N)\leq x\Big\},\ \ x \in \mathbb{R},
\end{eqnarray*}
for $N\geq 3$ and
\begin{eqnarray*}
S_{N,n}=\sum_{(i, j)\in \Lambda_{n,N}} T\hat{\rho}_{ij}^2
\end{eqnarray*}
for $N\geq n\geq 1.$ Observe that $B_{I_1}B_{I_2}\cdots B_{I_{n}}$ is an event generated by random vectors $\{\bd{e}_i, \bd{e}_j ;\, (i, j) \in \Lambda_{n,N}\}$. A crucial observation is that $S_N-S_{N,n}$ is independent of $B_{I_1}B_{I_2}\cdots B_{I_{n}}$.
It is easy to see that
\begin{eqnarray*}
S_{N,n} &=& \sum_{l=1}^n\sum_{j=i_l+1}^N T\hat{\rho}_{i_lj}^2+\sum_{l=1}^{N}\sum_{i=1}^{j_l-1} T\hat{\rho}_{ij_l}^2 -\sum_{s=1}^n \sum_{l=1}^nT\hat{\rho}_{i_lj_s}^2\\
& := & Q_{N,1} +Q_{N,2}-Q_{N,3}.
\end{eqnarray*}
For any integer $\tau\geq 2$, from a convex inequality we have
\begin{eqnarray*}
E\big(|Q_{N,1}-EQ_{N,1}|^{\tau}\big) &\leq & n^{\tau-1}\cdot \sum_{l=1}^n E\Big(\Big|\sum_{j=i_l+1}^N T(\hat{\rho}_{i_lj}^2-E\hat{\rho}_{i_lj}^2)\Big|^{\tau}\Big)\\
& \leq & C n^{\tau} T^{\tau}\cdot\Big(\frac{N^{\tau/2}}{m^{\tau}}+ \frac{N^{\tau}}{m^{2\tau}}\Big)\\
& \leq & C\cdot n^{\tau}N^{\tau/2}
\end{eqnarray*}
by Lemma \ref{wen_fei}, where the constant $C$ is free of $N$ and $T$, and where the last step follows from the assumption $N=o(T^2)$. Similarly,
\begin{eqnarray*}
E\big(|Q_{N,2}-EQ_{N,2}|^{\tau}\big) \leq  C\cdot n^{\tau}N^{\tau/2}.
\end{eqnarray*}
Lastly, by Lemma \ref{wen_fei} again,
\begin{eqnarray*}
E\big(|Q_{N,3}-EQ_{N,3}|^{\tau}\big)
&\leq & T^{\tau}\cdot n^{2(\tau-1)}\cdot \sum_{s=1}^n \sum_{l=1}^nE\big[|\hat{\rho}_{i_lj_s}^2-E\hat{\rho}_{i_lj_s}^2|^{\tau}\big]\\
& \leq & C\cdot n^{2\tau}.
\end{eqnarray*}
Therefore,
\begin{eqnarray*}
E|S_{N,n}-ES_{N,n}|^{\tau}\leq C\big(n^{\tau}N^{\tau/2}+n^{2\tau}\big).
\end{eqnarray*}
Fix $\epsilon\in (0,1).$ By the Markov inequality,
\begin{align}\lbl{skfenhy}
P\Big(\frac{1}{N}|S_{N,n}-ES_{N,n}|
\geq & \epsilon\Big) \leq \frac{C}{\epsilon^{\tau}}\cdot\frac{n^{\tau}}{N^{\tau/2}}\nonumber\\
= &C'\cdot\frac{n^{\tau}}{N^{\tau/2}}
\end{align}
for all $N\geq n^2$, where $C'$ is a constant depending on $\epsilon$ but free of $N$, $T$ or indices $\{I_1, \cdots, I_n\}.$

Fix $I_1< I_2< \cdots < I_{n}\in \Lambda_N$. By \eqref{skfenhy} and the definition of $A_N(x)$,
\begin{eqnarray*}
&& P(A_N(x)B_{I_1}B_{I_2}\cdots B_{I_{n}})\\
&\leq & P\Big(A_N(x)B_{I_1}B_{I_2}\cdots B_{I_{n}},\ \frac{1}{N}|S_{N,n}-ES_{N,n}|< \epsilon\Big) + C'\cdot\frac{n^{\tau}}{N^{\tau/2}}\\
& \leq & P\Big(\frac{1}{N}[(S_N-S_{N,n})-E(S_N-S_{N,n})]\leq x+\epsilon
,\ B_{I_1}B_{I_2}\cdots B_{I_{n}}\Big) +C'\cdot\frac{n^{\tau}}{N^{\tau/2}}\\
& = & P\Big(\frac{1}{N}[(S_N-S_{N,n})-E(S_N-S_{N,n})]\leq x+\epsilon\Big)\cdot P\big(
B_{I_1}B_{I_2}\cdots B_{I_{n}}\big) +C'\cdot\frac{n^{\tau}}{N^{\tau/2}}
\end{eqnarray*}
by the independence between $S_N-S_{N,n}$ and $B_{I_1}B_{I_2}\cdots B_{I_{n}}$. Now
\begin{eqnarray*}
&& P\Big(\frac{1}{N}[(S_N-S_{N,n})-E(S_N-S_{N,n})]\leq x+\epsilon\Big)\\
& \leq & P\Big(\frac{1}{N}[(S_N-S_{N,n})-E(S_N-S_{N,n})]\leq x+\epsilon,\ \\
& &  \frac{1}{N}|S_{N,n}-ES_{N,n}|< \epsilon\Big) + C'\cdot\frac{n^{\tau}}{N^{\tau/2}}\\
& \leq & P\Big(\frac{1}{N}(S_N-ES_N)\leq x+2\epsilon\Big) + C'\cdot\frac{n^{\tau}}{N^{\tau/2}}\\
&\leq & P\big(A_N(x+2\epsilon)\big) + C'\cdot\frac{n^{\tau}}{N^{\tau/2}}.
\end{eqnarray*}
Combing the two inequalities to get
\begin{align}\lbl{shakings}
& P(A_N(x)B_{I_1}B_{I_2}\cdots B_{I_{n}})\nonumber\\
\leq &P\big(A_N(x+2\epsilon)\big)\cdot P\big(
B_{I_1}B_{I_2}\cdots B_{I_{n}}\big)  +2C'\cdot\frac{n^{\tau}}{N^{\tau/2}}.
\end{align}
Similarly,
\begin{eqnarray*}
&&P\Big(\frac{1}{N}[(S_N-S_{N,n})-E(S_N-S_{N,n})]\leq x-\epsilon,\
B_{I_1}B_{I_2}\cdots B_{I_{n}}\Big)\\
&\leq & P\Big(\frac{1}{N}[(S_N-S_{N,n})-E(S_N-S_{N,n})]\leq x-\epsilon,
B_{I_1}B_{I_2}\cdots B_{I_{n}},\\
 & &\frac{1}{N}|S_{N,n}-ES_{N,n}|< \epsilon\Big)+C'\cdot\frac{n^{\tau}}{N^{\tau/2}}\\
& \leq & P\Big(\frac{1}{N}(S_N-ES_N)\leq x,\ B_{I_1}B_{I_2}\cdots B_{I_{n}}\Big) +C'\cdot\frac{n^{\tau}}{N^{\tau/2}}.
\end{eqnarray*}
In other words, by independence,
\begin{eqnarray*}
&& P(A_N(x)B_{I_1}B_{I_2}\cdots B_{I_{n}}) \\
&\geq &
P\Big(\frac{1}{N}[(S_N-S_{N,n})-E(S_N-S_{N,n})]\leq x-\epsilon\Big)\cdot \\
& &P\Big(B_{I_1}B_{I_2}\cdots B_{I_{n}}\Big)-C'\cdot\frac{n^{\tau}}{N^{\tau/2}}.
\end{eqnarray*}
Furthermore,
\begin{eqnarray*}
&& P\Big(\frac{1}{N}(S_N-ES_N)\leq x-2\epsilon\Big)\\
& \leq & P\Big(\frac{1}{N}(S_N-ES_N)\leq x-2\epsilon,\ \frac{1}{N}|S_{N,n}-ES_{N,n}|< \epsilon\Big) + C'\cdot\frac{n^{\tau}}{N^{\tau/2}}\\
& \leq & P\Big(\frac{1}{N}[(S_N-S_{N,n})-E(S_N-S_{N,n})]\leq x-\epsilon\Big) +C'\cdot\frac{n^{\tau}}{N^{\tau/2}}.
\end{eqnarray*}
The above two strings of inequalities imply
\begin{eqnarray*}
&& P(A_N(x)B_{I_1}B_{I_2}\cdots B_{I_{n}}) \\
&\geq & P\Big(\frac{1}{N}(S_N-ES_N)\leq x-2\epsilon\Big)\cdot P\Big(
B_{I_1}B_{I_2}\cdots B_{I_{n}}\Big)-2C'\cdot\frac{n^{\tau}}{N^{\tau/2}},
\end{eqnarray*}
which joining with \eqref{shakings} yields
\begin{eqnarray*}
&&\big|P(A_N(x)B_{I_1}B_{I_2}\cdots B_{I_{n}})-P(A_N(x))\cdot P(B_{I_1}B_{I_2}\cdots B_{I_{n}})\big|\\
& \leq & \Delta_{N, \epsilon}\cdot  P(B_{I_1}B_{I_2}\cdots B_{I_{n}})+4C'\cdot\frac{n^{\tau}}{N^{\tau/2}}
\end{eqnarray*}
where
\begin{eqnarray*}
\Delta_{N, \epsilon}:&=&|P(A_N(x))-P(A_N(x+2\epsilon)| + |P(A_N(x))-P(A_N(x-2\epsilon)|. \end{eqnarray*}
In particular,
\begin{align}\lbl{864123}
\Delta_{N, \epsilon} \to  |\Phi(x+2\epsilon)-\Phi(x)|+|\Phi(x-2\epsilon)-\Phi(x)|
\end{align}
as $N\to\infty$ by Theorem \ref{diligent}.  As a consequence,
\begin{eqnarray*}
\zeta(N,n):&=& \sum_{I_1< I_2< \cdots < I_{n}\in \Lambda_N}\big[P(A_N(x) B_{I_1}B_{I_2}\cdots B_{I_{n}}) - \\
 & &P(A_N(x))\cdot P(B_{I_1}B_{I_2}\cdots B_{I_{n}})\big]\\
& \leq & \sum_{I_1< I_2< \cdots < I_{n}\in \Lambda_N}\Big[\Delta_{N, \epsilon}\cdot  P(B_{I_1}B_{I_2}\cdots B_{I_{n}})+4C'\cdot\frac{n^{\tau}}{N^{\tau/2}}\Big]\\
& \leq & \Delta_{N, \epsilon}\cdot H(N, n)+ (4C')\cdot \binom{\frac{1}{2}N(N-1)}{n} \cdot\frac{n^{\tau}}{N^{\tau/2}},
\end{eqnarray*}
where
\begin{eqnarray*}
H(N, n)=\sum_{I_1< I_2< \cdots < I_{n}\in \Lambda_N}P(B_{I_1}B_{I_2}\cdots B_{I_{n}})
\end{eqnarray*}
as defined in Lemma \ref{one_key}. From \eqref{lncsivlbt}, we know $\limsup_{N\to\infty}H(N, n)\leq C/n!$, where $C$ is a universal constant. Picking  $\tau=6n$, and using the trivial fact $\binom{r}{s}\leq r^s$ for any integers $1\leq s \leq r$, we have that
\begin{eqnarray*}
\binom{\frac{1}{2}N(N-1)}{n} \cdot\frac{n^{\tau}}{N^{\tau/2}}
\leq N^{2n}\cdot \frac{n^{\tau}}{N^{\tau/2}}\leq \frac{n^{\tau}}{N^n}.
\end{eqnarray*}
Hence, from \eqref{864123}
\begin{eqnarray*}
\limsup_{N\to \infty}\zeta(N,n) &\leq & \frac{C}{n!}\cdot
\limsup_{N\to\infty}\Delta_{N, \epsilon}\\
& = & \frac{C}{n!}\cdot\big[|\Phi(x+2\epsilon)-\Phi(x)|+|\Phi(x-2\epsilon)-\Phi(x)|\big]
\end{eqnarray*}
for any $\epsilon>0$. The desired result follows by sending $\epsilon \downarrow 0.$ \hfill $\square$

\subsubsection{Finale: proof of Theorem \ref{Asym_indept}} We now are ready to assemble everything together.

\smallskip

\noindent\textbf{Proof of Theorem \ref{Asym_indept}}. Recall $\{\bd{e}_i;\ 1\leq i \leq N\}$ in \eqref{rain_dear_lost}. By assumption  \eqref{assa}, we see that $\{\bd{e}_i;\ 1\leq i \leq N\}$ are i.i.d. uniformly distributed over $\mathbb{S}^{T-1}$.  As in Lemma \ref{change_to_and}, define
\begin{eqnarray*}
\tilde{L}_N=\max_{1\leq i<j \leq N}|\bd{e}_i'\bd{e}_j|.
\end{eqnarray*}
Let $m=T-p.$ Recall
\begin{eqnarray*}
&& L_N=\max_{1\leq i< j \leq N}|\hat{\rho}_{ij}|;\\
&& S_N=\sum_{1\leq i< j \leq N}T\hat{\rho}_{ij}^2,\ \ \  \mu_N=\frac{T}{m^2}\sum_{j=2}^N\sum_{i=1}^{j-1}\mbox{tr}(\bd{P}_i\bd{P}_j).
\end{eqnarray*}
By Theorem 3 and Remark 2.1 from \cite{Cai_J2011}  and  Theorem \ref{diligent}, the following hold.
\begin{align}
&T\tilde{L}_N^2-4\log N +\log\log N \to F(y)=\exp(-(1/\sqrt{8\pi})e^{-y/2})\ \mbox{weakly}; \lbl{jian_jiao}\\
&\frac{1}{N}(S_N-\mu_N)\to N(0, 1)\ \mbox{weakly}.\ \ \ \ \ \ \lbl{gan_doufu}
\end{align}
To show asymptotic independence, by Lemma \ref{change_to_and}, it is enough to show
\begin{eqnarray*}
\lim_{N\to \infty}P\Big(\frac{1}{N}(S_N-\mu_N)\leq x,\ T\tilde{L}_N^2-4\log N +\log\log N\leq y\Big)= \Phi(x)\cdot F(y)
\end{eqnarray*}
for any $x\in \mathbb{R}$ and $y \in \mathbb{R}$, where $\Phi(x)=(2\pi)^{-1/2}\int_{-\infty}^xe^{-t^2/2}\,dt.$ Review \eqref{wet_soil} to see
\begin{align}\lbl{whale}
l_N=T^{-1/2}\cdot (4\log N -\log\log N+y)^{1/2},
\end{align}
which makes sense for large $N$. Because of \eqref{jian_jiao} and \eqref{gan_doufu}, the above is equivalent to that
\begin{align}\lbl{wealth_no}
\lim_{N\to \infty}P\Big(\frac{1}{N}(S_N-\mu_N)\leq x,\ \tilde{L}_N>l_N\Big)= \Phi(x)\cdot [1-F(y)]
\end{align}
for any $x\in \mathbb{R}$ and $y \in \mathbb{R}$. Review notations $ \Lambda_N$, $A_N$ and $B_{I}$ for any $I=(i,j)\in \Lambda_N$ in \eqref{haokanne}.
Write
\begin{align}\lbl{abci}
P\Big(\frac{1}{N}(S_N-\mu_N)\leq x,\ \tilde{L}_N>l_N\Big)=P\Big(\bigcup_{I\in \Lambda_N}A_NB_{I}\Big).
\end{align}
Here the notation $A_NB_I$ stands for $A_N\cap B_I$.
From the inclusion-exclusion principle,
\begin{align}
P\Big(\bigcup_{I\in \Lambda_N}A_NB_{I}\Big)  \leq&  \sum_{I_1\in \Lambda_N}P(A_NB_{I_1})-\sum_{I_1< I_2\in \Lambda_N}P(A_NB_{I_1}B_{I_2})+\cdots+\nonumber\\
 & \sum_{I_1< I_2< \cdots < I_{2k+1}\in \Lambda_N}P(A_NB_{I_1}B_{I_2}\cdots B_{I_{2k+1}})\nonumber\\
&  \lbl{Upper_bound}
\end{align}
and
\begin{align}
P\Big(\bigcup_{I\in \Lambda_N}A_NB_{I}\Big)  \geq  \sum_{I_1\in \Lambda_N}P(A_NB_{I_1})-&\sum_{I_1< I_2\in \Lambda_N}P(A_NB_{I_1}B_{I_2})+\cdots- \nonumber\\
 & \sum_{I_1< I_2< \cdots < I_{2k}\in \Lambda_N}P(A_NB_{I_1}B_{I_2}\cdots B_{I_{2k}}) \nonumber\\
&  \lbl{Lower_bound}
\end{align}
for any integer $k\geq 1$. Reviewing the definition
\begin{eqnarray*}
H(N, n)=\sum_{I_1< I_2< \cdots < I_{n}\in \Lambda_N}P(B_{I_1}B_{I_2}\cdots B_{I_{n}})
\end{eqnarray*}
for $n\geq 1$ in Lemma \ref{one_key}, we have from the lemma that
\begin{align}\lbl{Maya}
\lim_{n\to\infty}\limsup_{N\to\infty}H(N, n)=0.
\end{align}
Set
\begin{eqnarray*}
&&\zeta(N,n)\\
&=&\sum_{I_1< I_2< \cdots < I_{n}\in \Lambda_N}\big[P(A_NB_{I_1}B_{I_2}\cdots B_{I_{n}}) - P(A_N)\cdot P(B_{I_1}B_{I_2}\cdots B_{I_{n}})\big]
\end{eqnarray*}
for $n\geq 1.$ By Lemma \ref{akxbcvgf},
\begin{align}\lbl{back_campus}
\lim_{N\to\infty}\zeta(N,n)=0
\end{align}
for each $n\geq 1$. The assertion \eqref{Upper_bound} implies that
\begin{align}\lbl{639475}
&P\Big(\bigcup_{I\in \Lambda_N}A_NB_{I}\Big)\nonumber\\
 \leq & P(A_N)\Big[\sum_{I_1\in \Lambda_N}P(B_{I_1})-\sum_{I_1< I_2\in \Lambda_N}P(B_{I_1}B_{I_2})+\cdots-  \nonumber\\
& \sum_{I_1< I_2< \cdots < I_{2k}\in \Lambda_N}P(B_{I_1}B_{I_2}\cdots B_{I_{2k}})\Big]+ \Big[\sum_{n=1}^{2k}\zeta(N,n)\Big] + H(N, 2k+1)  \nonumber\\
\leq & P(A_N)\cdot P\Big(\bigcup_{I\in \Lambda_N}B_{I}\Big)+ \Big[\sum_{n=1}^{2k}\zeta(N,n)\Big] + H(N, 2k+1),
\end{align}
where the inclusion-exclusion formula is used again in the last inequality, that is,
\begin{eqnarray*}
P\Big(\bigcup_{I\in \Lambda_N}B_{I}\Big) &\geq & \Big[\sum_{I_1\in \Lambda_N}P(B_{I_1})-\sum_{I_1< I_2\in \Lambda_N}P(B_{I_1}B_{I_2})+\cdots - \nonumber\\
&&~~~~~~~~~~~~~~~~ \sum_{I_1< I_2< \cdots < I_{2k}\in \Lambda_N}P(B_{I_1}B_{I_2}\cdots B_{I_{2k}})\Big]
\end{eqnarray*}
for all $k\geq 1$.
 By the definition of $l_N$ and \eqref{jian_jiao},
\begin{eqnarray*}
 P\Big(\bigcup_{I\in \Lambda_N}B_{I}\Big)=P(\tilde{L}_N>l_N)=P(T \tilde{L}_N^2-4\log N +\log\log N> y) \to 1-F(y)
\end{eqnarray*}
as $N\to\infty$. By \eqref{gan_doufu}, $P(A_N)\to \Phi(x)$ as $N\to\infty.$ From \eqref{abci}, by fixing $k$ first and sending $N\to \infty$ we get from \eqref{back_campus} that
\begin{eqnarray*}
& &\limsup_{N\to\infty}P\Big(\frac{1}{N}(S_N-\mu_N)\leq x,\ \tilde{L}_N>l_N\Big)\\
&\leq& \Phi(x)\cdot [1-F(y)] +\limsup_{N\to\infty}H(N, 2k+1).
\end{eqnarray*}
Now, let $k\to \infty$ and use \eqref{Maya} to see
\begin{align}\lbl{vskdnti}
\limsup_{N\to\infty}P\Big(\frac{1}{N}(S_N-\mu_N)\leq x,\ \tilde{L}_N>l_N\Big)\leq \Phi(x)\cdot [1-F(y)].
\end{align}
By applying the same argument to \eqref{Lower_bound}, we see that the counterpart of \eqref{639475} becomes
\begin{eqnarray*}
P\Big(\bigcup_{I\in \Lambda_N}A_NB_{I}\Big)
& \geq & P(A_N)\Big[\sum_{I_1\in \Lambda_N}P(B_{I_1})-\sum_{I_1< I_2\in \Lambda_N}P(B_{I_1}B_{I_2})+\cdots + \nonumber\\
&& \sum_{I_1< I_2< \cdots < I_{2k-1}\in \Lambda_N}P(B_{I_1}B_{I_2}\cdots B_{I_{2k-1}})\Big]+\nonumber\\
&& \Big[\sum_{n=1}^{2k-1}\zeta(N,n)\Big] - H(N, 2k)  \nonumber\\
&\geq & P(A_N)\cdot P\Big(\bigcup_{I\in \Lambda_N}B_{I}\Big) + \Big[\sum_{n=1}^{2k-1}\zeta(N,n)\Big] - H(N, 2k).
\end{eqnarray*}
where in the last step we use the inclusion-exclusion principle such that
\begin{eqnarray*}
P\Big(\bigcup_{I\in \Lambda_N}B_{I}\Big) &\leq & \Big[\sum_{I_1\in \Lambda_N}P(B_{I_1})-\sum_{I_1< I_2\in \Lambda_N}P(B_{I_1}B_{I_2})+\cdots + \nonumber\\
&&~~~~~~~~~~~~~~~~ \sum_{I_1< I_2< \cdots < I_{2k-1}\in \Lambda_N}P(B_{I_1}B_{I_2}\cdots B_{I_{2k-1}})\Big]
\end{eqnarray*}
for all $k\geq 1$. Review \eqref{abci} and repeat the earlier procedure to see
\begin{eqnarray*}
\liminf_{N\to\infty}P\Big(\frac{1}{N}(S_N-\mu_N)\leq x,\ \tilde{L}_N>l_N\Big)\geq \Phi(x)\cdot [1-F(y)]
\end{eqnarray*}
by sending $N\to \infty$ and then sending $k\to\infty.$
This and \eqref{vskdnti} yield \eqref{wealth_no}. The proof is completed. \hfill$\Box$

\section*{Acknowledgment}

Professors Feng and Liu thank NSFC grants 11501092 and 11571068 for partially support.
Professor  Jiang thanks NSF Grants  DMS-1406279 and DMS-1916014 for partially support.

\appendix

\section*{Appendix}

In this part we will prove the technical results stated in previous sections. We create same number of sections to accumulate the proofs of the claims in the corresponding section.

\subsection{Proofs of auxiliary results in Section \ref{Pre_123}}\lbl{att_123} We prove the lemmas in the order of their numerations.

\smallskip

\noindent\textbf{Proof of Lemma \ref{mac_pencil}}. The conclusions are about even functions of $\xi$. So, without loss of generality, assume $a\geq 0.$  Set $E\xi^2=b^2$ for some $b>0$. Then $b\geq a$. Note that
\begin{eqnarray*}
E[|\xi^2-E\xi^2|^{\tau}] &=&E[|\xi-b|^{\tau}\cdot|\xi+b|^{\tau}]\\
&\leq & \big[E|\xi-b|^{2\tau}\big]^{1/2}\cdot \big[E|\xi+b|^{2\tau}\big]^{1/2}
\end{eqnarray*}
by the Cauchy-Schwartz inequality. Notice
\begin{eqnarray*}
 E|\xi-b|^{2\tau}\leq 2^{2\tau-1}\cdot [|a-b|^{2\tau}+E|\xi-a|^{2\tau}]
\end{eqnarray*}
and
\begin{eqnarray*}
 E[|\xi+b|^{2\tau}]
 &\leq & 2^{2\tau-1}\cdot \big(|a+b|^{2\tau}+E|\xi-a|^{2\tau}\big)\\
 &\leq & 2^{4\tau-2}\cdot \big[|a-b|^{2\tau}+(2a)^{2\tau}+E|\xi-a|^{2\tau}\big]\\
 & \leq & 2^{6\tau-2}\cdot \big(|a-b|^{2\tau}+a^{2\tau}+E|\xi-a|^{2\tau}\big).
\end{eqnarray*}
Use the inequality $\sqrt{x+y+z}\leq \sqrt{x}+\sqrt{y}+\sqrt{z}$ for all $x\geq 0$, $y\geq 0$ and $z\geq 0$ to see
\begin{align}\lbl{skfwo}
&E[|\xi^2-E\xi^2|^{\tau}] \nonumber\\
\leq & 2^{4\tau-1.5}\cdot \big[|b-a|^{\tau}+\sqrt{E|\xi-a|^{2\tau}}\big]\cdot \big[|b-a|^{\tau}+a^{\tau}+\sqrt{E|\xi-a|^{2\tau}}\big].
\end{align}
If $a=0$, then
\begin{eqnarray*}
E[|\xi^2-E\xi^2|^{\tau}]
&\leq &  2^{4\tau-1.5}\cdot\Big[ {\rm Var}(\xi)^{\tau/2}+\sqrt{E(|\xi|^{2\tau})}\,\Big]^2\\
&\leq & 16^{\tau}\cdot
\big[ \mbox{Var}(\xi)^{\tau}+E(|\xi|^{2\tau})\,\big]
\end{eqnarray*}
This leads to (i) since $\mbox{Var}(\xi)^{\tau}=[E(\xi^2)]^{\tau}\leq E(|\xi|^{2\tau})$ by the H\"{o}lder inequality. Now, if $a\ne 0$, we continue from \eqref{skfwo} to see
\begin{eqnarray*}
|b-a|^{\tau}=\frac{(b^2-a^2)^{\tau}}{(a+b)^{\tau}}=\frac{1}{(a+b)^{\tau}}\cdot {\rm Var}(\xi)^{\tau}\leq \frac{1}{a^{\tau}}\cdot {\rm Var}(\xi)^{\tau}.
\end{eqnarray*}
We get (ii). The proof is finished. \hfill$\Box$

\subsection{Proofs of auxiliary results in Section \ref{beauty_ugly}}\lbl{att_proof}

In this part we develop some identities and inequalities regarding moments of random vectors with the uniform distribution on high-dimensional spheres. We will focus on developing basic tools. They are of independent interest. Review  notations $\bd{P}_i$ and $\bd{\epsilon}_i=(\epsilon_{i1},\cdots,\epsilon_{iT})^{'}\in \mathbb{R}^T$ in \eqref{shadongxi} and \eqref{park_1}.

\smallskip

\noindent\textbf{Proof of Lemma \ref{trivial12}}. (i) Notice  $\bd{M}_{ij}=(\bd{U}_j'\bd{U}_i)'(\bd{U}_j'\bd{U}_i)$. Automatically $\bd{M}_{ij}$ is non-negative definite. To show $\bd{I}-\bd{M}_{ij}$ is non-negative definite, it is enough to prove
\begin{align}\lbl{beer_large}
\bd{x}'\bd{M}_{ij}\bd{x}\leq \bd{x}'\bd{x}
\end{align}
for any $\bd{x}\in \mathbb{R}^m$. In fact, let $\bd{z}=\bd{U}_i\bd{x}$. Then
\begin{eqnarray*}
\bd{x}'\bd{M}_{ij}\bd{x}=\|\bd{U}_j'\bd{z}\|^2
=\bd{z}'\bd{U}_j\bd{U}_j'\bd{z}=\bd{z}'\bd{P}_j\bd{z}
\end{eqnarray*}
by \eqref{happy_moon}. By \eqref{park_1} and \eqref{happy_moon} again,
\begin{eqnarray*}
\bd{z}'\bd{P}_j\bd{z}\leq \bd{z}'\bd{z}=\bd{x}'\bd{U}_i'\bd{U}_i\bd{x}=\bd{x}'\bd{x}.
\end{eqnarray*}
 The above two assertions lead to \eqref{beer_large}.

(ii) Since $\bd{I}-\bd{M}$ is  non-negative definite, then all of the eigenvalues of $\bd{M}$ are in the interval $[0, 1]$. So is $\bd{M}^2$. This gives the conclusion.
 \hfill$\Box$

\smallskip


\noindent\textbf{Proof of Lemma \ref{trivial}}. Note
\begin{align}\lbl{Science}
\mbox{tr}(\bd{F}_1\bd{F}_2)=\mbox{tr}(\bd{F}_2\bd{F}_1)
\end{align}
for any matrices $\bd{F}_1$ and $\bd{F}_2.$  Write $\bd{M}_1=\bd{H}'\,\mbox{diag}(1, \cdots, 1, 0, \cdots, 0)\bd{H}$ where $\bd{H}$ is an orthogonal matrix, and the number of $1$'s is equal to $r:= \mbox{rank} (\bd{M}_1)$. Recall \eqref{Science}.
Then
\begin{eqnarray*}
\mbox{tr}\,(\bd{M}_1\bd{M}_2) &= & \mbox{tr}\,\big[\mbox{diag}(1, \cdots, 1, 0, \cdots, 0)(\bd{H}\bd{M}_2\bd{H}')\big]\\
&= & \mbox{tr}\,(\mbox{the upper-left $r\times r$ submatrix of $\bd{H}\bd{M}_2\bd{H}'$})\\
& \leq & \mbox{tr}\,(\bd{H}\bd{M}_2\bd{H}')\\
& = & \mbox{tr}\,(\bd{M}_2)
\end{eqnarray*}
by \eqref{Science},  where the inequality is obtained because $\bd{H}\bd{M}_2\bd{H}'$ is nonnegative definite and hence all of its diagonal entries are non-negative. The conclusion follows.  \hfill$\Box$

\smallskip

%

\noindent\textbf{Proof of Lemma \ref{brother_cat}}. Pick a non-negative matrix $\bd{M}_1^{1/2}$ such that $\bd{M}_1^{1/2}\cdot \bd{M}_1^{1/2}=\bd{M}_1.$ Recall the fact that $\bd{A}\bd{B}$ and $\bd{B}\bd{A}$ have the same eigenvalues for any square matrices $\bd{A}$ and $\bd{B}$. Then $\bd{M}_1\bd{M}_2$ and $\bd{M}_1^{1/2}\bd{M}_2\bd{M}_1^{1/2}$ have the same eigenvalues. Since the latter one is readily seen to be a non-negative definite matrix, we know that all of the eigenvalues of  $\bd{M}_1\bd{M}_2$ are non-negative. In particular, $\mbox{tr}(\bd{M}_1\bd{M}_2) \geq 0$.

The second conclusion holds trivially for $r=0$. We next assume $r\geq 1.$  Let $\bd{M}\ne \bd{0}$ be an $n\times n$ real matrix. Assume all eigenvalues are real and the non-zero eigenvalues are  $\lambda_1, \cdots,   \lambda_v$ with $1\leq v \leq n.$ Then,
\begin{align}\label{shine_animal}
\mbox{tr}(\bd{M}^2)=\lambda_1^2+ \cdots+ \lambda_v^2\geq \frac{1}{v}(\lambda_1+\cdots +\lambda_v)^2=\frac{1}{v}[\mbox{tr}(\bd{M})]^2.
\end{align}
From the singular value decomposition theorem (see e.g., p. 150 from \cite{HJ2ndED}), we see $v$ is the same as the number of non-zero singular values of $\bd{M}.$ Let $s_1\geq \cdots \geq s_n$ be the singular values of $\bd{M}$, that is, the eigenvalues of $(\bd{M}'\bd{M})^{1/2}$. Assume  $|\lambda_1|\geq \cdots\geq  |\lambda_v|$ without loss of generality.  We then have from the Weyl inequality that $|\lambda_1\cdots\lambda_k|\leq s_1\cdots s_k$ for all $1\leq k \leq n$; see, for example, p. 454 from \cite{HJ2ndED}. This implies that $v$ is no more than the number of non-zero eigenvalues of $(\bd{M}'\bd{M})^{1/2}$, which is the same as the number of non-zero eigenvalues of $\bd{M}'\bd{M}$, which is again equal to $\textrm{rank}(\bd{M}'\bd{M})=\textrm{rank}(\bd{M})$. That is, $v\leq \textrm{rank}(\bd{M}).$ This and \eqref{shine_animal} yield the desired conclusion by taking $\bd{M}=\bd{M}_1\bd{M}_2.$ \hfill$\square$


\smallskip

\noindent\textbf{Proof of Lemma \ref{sun_quiet}}. First,
\begin{eqnarray*}
\mbox{Var}(\bd{d}'\bd{M}\bd{d})
=E[(\bd{d}'\bd{M}\bd{d})^2]-
\big[E(\bd{d}'\bd{M}\bd{d})\big]^2.
\end{eqnarray*}
Then (iii) follows if (i) and (ii) are valid. Let us prove (i) and (ii) next.

Write $\bd{M}=\bd{H}'\,\mbox{diag}(\lambda_1, \cdots, \lambda_m)\bd{H}$ where $\bd{H}$ is an orthogonal matrix. Set $\bd{\eta}=(Z_1, \cdots, Z_m)'$. Observe $\bd{H}\bd{d}=\frac{\bd{H}\bd{\eta}}{\|\bd{H}\bd{\eta}\|}.$ By the orthogonal invariance of Gaussian distributions, $\bd{H}\bd{\eta}$ and $\bd{\eta}$ have the same distribution, so are $\bd{H}\bd{d}$ and $\bd{d}$. As a consequence, $\bd{d}'\bd{M}\bd{d}$ and $\frac{\lambda_1Z_1^2+\cdots \lambda_mZ_m^2}{Z_1^2+\cdots +Z_m^2}$ have a common distribution.
Easily,
\begin{eqnarray*}
E\frac{\lambda_1Z_1^2+\cdots \lambda_mZ_m^2}{Z_1^2+\cdots +Z_m^2}=E\frac{Z_1^2}{Z_1^2+\cdots +Z_m^2}\cdot\sum_{i=1}^m\lambda_i=\frac{1}{m}\cdot \mbox{tr}(\bd{M})
\end{eqnarray*}
by Lemma \ref{geese} with $a_1=1$ and other $a_i$'s being equal to zero. We get (i).
Now, use the formula $(a_1+\cdots +a_m)^2=\sum_{i=1}^ma_i^2+2\sum_{1\leq i<j \leq m}a_ia_j$ to see that
\begin{eqnarray*}
&& E\frac{(\lambda_1Z_1^2+\cdots \lambda_mZ_m^2)^2}{(Z_1^2+\cdots +Z_m^2)^2}\\
&= &E\frac{Z_1^4}{(Z_1^2+\cdots +Z_m^2)^2}\cdot\Big(\sum_{i=1}^m\lambda_i^2\Big)+ E\frac{2Z_1^2Z_2^2}{(Z_1^2+\cdots +Z_m^2)^2}\cdot\sum_{1\leq i<j \leq m}\lambda_i\lambda_j.
\end{eqnarray*}
By Lemma \ref{geese},
\begin{align}\lbl{fall_flower}
E\frac{Z_1^4}{(Z_1^2+\cdots +Z_m^2)^2}=\frac{3}{m(m+2)}\ \ \mbox{and}\ \ E\frac{Z_1^2Z_2^2}{(Z_1^2+\cdots +Z_m^2)^2}=\frac{1}{m(m+2)}.
\end{align}
Hence,
\begin{eqnarray*}
&&E[(\bd{d}'\bd{M}\bd{d})^2]\\
&=&
\frac{3}{m(m+2)}\cdot\Big(\sum_{i=1}^m\lambda_i^2\Big)+\frac{2}{m(m+2)}\cdot\Big(\sum_{1\leq i<j \leq m}\lambda_i\lambda_j\Big).
\end{eqnarray*}
Write $2\sum_{1\leq i<j \leq m}\lambda_i\lambda_j=(\sum_{i=1}^m\lambda_i)^2-\sum_{i=1}^m\lambda_i^2$. Use the relations  $\mbox{tr}(\bd{M})=\sum_{i=1}^m\lambda_i$ and $\mbox{tr}(\bd{M}^2)=\sum_{i=1}^m\lambda_i^2$ to see that
\begin{align}\lbl{jelly}
E[(\bd{d}'\bd{M}\bd{d})^2]=\frac{2}{m(m+2)}\cdot \mbox{tr}(\bd{M}^2)+\frac{1}{m(m+2)}\cdot [\mbox{tr}(\bd{M})]^2.
\end{align}
We get (ii). \hfill$\Box$

\smallskip

\noindent\textbf{Proof of Lemma \ref{dust_bone}}. First, by Lemma \ref{sun_quiet}, $E(\bd{d}'\bd{M}\bd{d})=\frac{1}{m}\mbox{tr}(\bd{M})$. As shown in the proof of Lemma \ref{sun_quiet}, without loss of generality, we assume $\bd{M}=\mbox{diag}(\lambda_1, \cdots, \lambda_m).$ Write
\begin{eqnarray*}
\bd{d}'\bd{M}\bd{d}-\frac{1}{m}\mbox{tr}(\bd{M})
&=&\frac{\lambda_1Z_1^2+\cdots + \lambda_mZ_m^2}{Z_1^2+ \cdots +Z_m^2}-\frac{\lambda_1+\cdots + \lambda_m}{m}\\
& = & \frac{(\lambda_1-\bar{\lambda})(Z_1^2-1)+\cdots + (\lambda_m-\bar{\lambda})(Z_m^2-1)}{Z_1^2+ \cdots +Z_m^2}
\end{eqnarray*}
where $\bar{\lambda}=\frac{\lambda_1+\cdots + \lambda_m}{m}.$ For clarity, set $a_i=\lambda_i-\bar{\lambda}$ and $\xi_i=Z_i^2-1$ for $i=1,\cdots, m$. By H\"{o}lder's inequality,
\begin{align}\lbl{marriage}
&E\Big[\bd{d}'\bd{M}\bd{d}-\frac{1}{m}\mbox{tr}(\bd{M})\Big]^{\tau}\nonumber\\
\leq & \big(E|a_1\xi_1+\cdots +a_m\xi_m|^{2\tau}\big)^{1/2}\cdot \big[E(Z_1^2+\cdots + Z_m^2)^{-2\tau}\big]^{1/2}.
\end{align}
From \eqref{hot_sleep1},  there exists  a constant $K_{\tau}>0$  depending on $\tau$ only such that
\begin{eqnarray*}
E|a_1\xi_1+\cdots +a_m\xi_m|^{2\tau}&\leq & K_{\tau}\cdot E(a_1^2\xi_1^2+\cdots +a_m^2\xi_m^2)^{\tau}.
\end{eqnarray*}
Set $b_i=a_i^2(a_1^2+\cdots + a_m^2)^{-1}$ for $i=1,\cdots, m$. Then  $b_1+\cdots + b_m=1.$ Notice $\varphi(x):=x^{\tau}$ is convex over $[0, \infty)$ since $\tau>1.$ Then
\begin{eqnarray*}
(b_1\xi_1^2+\cdots +b_m\xi_m^2)^{\tau}\leq b_1|\xi_1|^{2\tau}+\cdots +b_m|\xi_m|^{2\tau}.
\end{eqnarray*}
This implies that
\begin{eqnarray*}
(a_1^2\xi_1^2+\cdots +a_m^2\xi_m^2)^{\tau}
\leq (a_1^2+\cdots + a_m^2)^{\tau-1}
\cdot \big(a_1^2|\xi_1|^{2\tau}+\cdots +a_m^2|\xi_m|^{2\tau}\big).
\end{eqnarray*}
Hence
\begin{align}\lbl{bike_walk}
&E|a_1\xi_1+\cdots +a_m\xi_m|^{2\tau} \nonumber\\
 \leq & K_{\tau}\cdot(a_1^2+\cdots + a_m^2)^{\tau-1}
\cdot\big[a_1^2(E|\xi_1|^{2\tau})+\cdots +a_m^2E(|\xi_m|^{2\tau})\big] \nonumber\\
=& K_{\tau}\cdot (a_1^2+\cdots + a_m^2)^{\tau} \cdot E(|\xi_1|^{2\tau}).
\end{align}
Now we bound the last term in \eqref{marriage}. Since $Z_1^2+\cdots + Z_m^2$ has the $\chi^2$ distribution with $m$-degree of freedom,
\begin{eqnarray*}
E(Z_1^2+\cdots + Z_m^2)^{-2\tau}
&=& \frac{1}{2^{m/2}\Gamma(m/2)}\int_0^{\infty}x^{-2\tau}\cdot x^{(m/2)-1}e^{-x/2}\,dx\\
&=& \frac{2^{(m/2)-2\tau}\Gamma((m/2)-2\tau)}{2^{m/2}\Gamma(m/2)}.
\end{eqnarray*}
It is known that
\begin{eqnarray*}
\lim_{x\to\infty}\frac{\Gamma(x+a)}{x^a\Gamma(x)}=1
\end{eqnarray*}
for any $a \in \mathbb{R}$, see, e.g., Lemma 2.4 from \cite{DJL12}. Therefore, there exists a constant $K_{\tau}'$ such that
\begin{align}\lbl{know_lap}
E(Z_1^2+\cdots + Z_m^2)^{-2\tau}\leq K_{\tau}'\cdot \frac{1}{m^{2\tau}}
\end{align}
for every $m\geq 4\tau+1$ in which case $\Gamma((m/2)-2\tau)$ is finite. This, \eqref{marriage} and \eqref{bike_walk} conclude
\begin{eqnarray*}
E\Big[\bd{d}'\bd{M}\bd{d}-\frac{1}{m}\mbox{tr}(\bd{M})\Big]^{\tau} \leq C_{\tau}\cdot(a_1^2+\cdots + a_m^2)^{\tau/2}\cdot \frac{1}{m^{\tau}}
\end{eqnarray*}
for all $m\geq 4\tau+1$, where $C_{\tau}$ is a constant depending on $\tau$ only. Trivially,
\begin{eqnarray*}
\sum_{i=1}^ma_i^2=\sum_{i=1}^m(\lambda_i-\bar{\lambda})^2
&=&\sum_{i=1}^m\lambda_i^2-\frac{1}{m}\big(\sum_{i=1}^m\lambda_i\big)^2\\
&=&\mbox{tr}(\bd{M}^2)-\frac{1}{m}[\mbox{tr}(\bd{M})]^2.
\end{eqnarray*}
The lemma is proved. \hfill$\Box$

\smallskip

\noindent\textbf{Proof of Lemma \ref{gan_geming}}. From Lemma \ref{sun_quiet}, $E(\bd{d}'\bd{M}\bd{d})=\frac{1}{m}\mbox{tr}(\bd{M})$. The second conclusion comes from Lemma \ref{dust_bone} directly by using the formula $(x+y)^{\tau}\leq 2^{\tau-1}(x^{\tau}+ y^{\tau})$ for all $x\geq 0$ and $y\geq 0.$ Since  $\bd{b}:=\bd{a}/\|\bd{a}\|$ is a unit vector, then
\begin{eqnarray*}
\bd{b}'\bd{d}\ \mbox{and}\  \frac{Z_1}{(Z_1^2+ \cdots +Z_m^2)^{1/2}}\ \mbox{have the same distribution};
\end{eqnarray*}
 see, for instance, Theorem 1.5.7 (i) and (5) on p. 147 from Muirhead (1982).  It follows that
\begin{eqnarray*}
E(|\bd{a}'\bd{d}|^{2\tau})
&=& \|\bd{a}\|^{2\tau}\cdot E\frac{|Z_1|^{2\tau}}{(Z_1^2+ \cdots +Z_m^2)^{\tau}}\\
& \leq & \|\bd{a}\|^{2\tau}\cdot (E|Z_1|^{4\tau})^{1/2}\cdot \big[E(Z_1^2+ \cdots +Z_m^2)^{-2\tau}\big]^{1/2}.
\end{eqnarray*}
The first conclusion then follows from \eqref{know_lap}. \hfill$\Box$

\smallskip

\noindent\textbf{Proof of Lemma \ref{R_V}}. (i) Trivially,
\begin{align}\lbl{lsu}
(\bd{h}'\bd{A}\bd{h})(\bd{h}'\bd{B}\bd{h})=
\frac{1}{2}\big\{\big[\bd{h}'(\bd{A}+\bd{B})\bd{h}\big]^2-
(\bd{h}'\bd{A}\bd{h})^2-(\bd{h}'\bd{B}\bd{h})^2\big\}.
\end{align}
From (ii) of Lemma \ref{sun_quiet}, we know
\begin{eqnarray*}
E[(\bd{h}'\bd{M}\bd{h})^2]=\frac{1}{m(m+2)}\cdot \big[2\,\mbox{tr}(\bd{M}^2)+(\mbox{tr}(\bd{M}))^2\big]
\end{eqnarray*}
for any symmetric matrix $\bd{M}.$ Then, by \eqref{lsu},
\begin{eqnarray*}
& & 2m(m+2)\cdot E\big[(\bd{h}'\bd{A}\bd{h})(\bd{h}'\bd{B}\bd{h})\big]\\
&=&2\,\mbox{tr}((\bd{A}+\bd{B})^2)+(\mbox{tr}(\bd{A}+\bd{B}))^2-\\
& &\big[2\,\mbox{tr}(\bd{A}^2)+(\mbox{tr}(\bd{A}))^2\big]-\big[2\,\mbox{tr}(\bd{B}^2)+(\mbox{tr}(\bd{B}))^2\big].
\end{eqnarray*}
A simple manipulation leads to (i).

(ii) By singular value decomposition, write $\bd{C}=\bd{H}_1'\mbox{diag}(\lambda_1, \cdots, \lambda_m)\bd{H}_2$, where $\bd{H}_1$ and $\bd{H}_2$ are orthogonal matrices, and where $\lambda_1^2, \cdots, \lambda_m^2$ are the eigenvalues of $\bd{C}\bd{C}'$. Now $\bd{h}_1'\bd{C}\bd{h}_2=(\bd{H}_1\bd{h}_1)'\mbox{diag}(\lambda_1, \cdots, \lambda_m)(\bd{H}_2\bd{h}_2)$. Since $\bd{h}_1$ and $\bd{h}_2$ are i.i.d. and orthogonal-invariant, we know $\bd{H}_1\bd{h}_1$ and $\bd{H}_2\bd{h}_2$ are also i.i.d. and have the same distribution as that of $\bd{h}.$ So we are able to write
\begin{eqnarray*}
\bd{h}_1'\bd{C}\bd{h}_2=\frac{1}{\|\bd{v}\|\cdot\|\bd{w}\|}\cdot \sum_{i=1}^m\lambda_iv_iw_i
\end{eqnarray*}
where $\bd{v}=(v_1, \cdots, v_m)$, $\bd{w}=(w_1, \cdots, w_m)$ and $\{v_i, w_i;\, 1\leq i \leq m\}$ are i.i.d. $N(0,1)$-distributed random variables. By the definition of variance and the Cauchy-Schwartz inequality,
\begin{eqnarray*}
\mbox{Var}[(\bd{h}_1'\bd{C}\bd{h}_2)^2]
&\leq & E[(\bd{h}_1'\bd{C}\bd{h}_2)^4]\\
& \leq & \Big(E\frac{1}{\|\bd{v}\|^8\cdot\|\bd{w}\|^8}\Big)^{1/2}\cdot \Big[E\Big(\sum_{i=1}^m\lambda_iv_iw_i\Big)^8\Big]^{1/2}\\
&= &\big[E(\|\bd{v}\|^{-8})\big]\cdot \Big[E\Big(\sum_{i=1}^m\lambda_iv_iw_i\Big)^8\Big]^{1/2},
\end{eqnarray*}
where the last step follows from independence. By \eqref{hot_sleep2},
\begin{eqnarray*}
E\Big(\sum_{i=1}^m\lambda_iv_iw_i\Big)^8 \leq Km^3\sum_{i=1}^m\lambda_i^8=Km^3\cdot \mbox{tr}[(\bd{C}\bd{C}')^4],
\end{eqnarray*}
where $K$ is a constant. Take $\tau=2$ from \eqref{know_lap}, we have
\begin{eqnarray*}
E(\|\bd{v}\|^{-8}) \leq \frac{K'}{m^4}
\end{eqnarray*}
where $K'$ is a constant. This concludes
\begin{eqnarray*}
\mbox{Var}[(\bd{h}_1'\bd{C}\bd{h}_2)^2]\leq \frac{K'\sqrt{K}}{m^{5/2}}\cdot \sqrt{\mbox{tr}[(\bd{C}\bd{C}')^4]}.
\end{eqnarray*}

(iii). Notice
\begin{eqnarray*}
&& \mbox{Cov}\big[(\bd{h}'\bd{A}\bd{h}_1)^2, (\bd{h}'\bd{B}\bd{h}_2)^2\big]\\
&=& E\big[(\bd{h}'\bd{A}\bd{h}_1)^2 (\bd{h}'\bd{B}\bd{h}_2)^2\big]-E\big[(\bd{h}'\bd{A}\bd{h}_1)^2]\cdot  E[(\bd{h}'\bd{B}\bd{h}_2)^2\big].
\end{eqnarray*}
Observe $E(\bd{h}\bd{h}')=E(\bd{h}_1\bd{h}_1')=\frac{1}{m}\bd{I}_m$ because of the structure of $\bd{h}$ appeared in  Lemma \ref{sun_quiet}. Then, use the fact $\bd{h}'\bd{A}\bd{h}_1=\bd{h}_1'\bd{A}\bd{h}$ and independence to have
\begin{eqnarray*}
E\big[(\bd{h}'\bd{A}\bd{h}_1)^2]
&=&
E\,\mbox{tr}\big[\bd{A}(\bd{h}_1\bd{h}_1')\bd{A}'(\bd{h}\bd{h}')\big]\\
&= & \mbox{tr}\big\{\bd{A}[E(\bd{h}_1\bd{h}_1')]\bd{A}'[E(\bd{h}\bd{h}')]\big\}\\
& = & \frac{1}{m^2}\,\mbox{tr}(\bd{A}\bd{A}').
\end{eqnarray*}
The above is also true if $\bd{A}$ is replaced by $\bd{B}.$ For a vector $\bd{a}\in \mathbb{R}^m$, we see that $E(\bd{a}'\bd{h}_2)^2=E(\bd{h}_2'\bd{a}\bd{a}'\bd{h}_2)=\frac{1}{m}\|\bd{a}\|^2$ by (i) of Lemma \ref{sun_quiet}. Conditioning on $\bd{h}$, using independence and by the proved (i), we obtain
\begin{eqnarray*}
E\big[(\bd{h}'\bd{A}\bd{h}_1)^2 (\bd{h}'\bd{B}\bd{h}_2)^2\big]
&=&\frac{1}{m^2}\cdot E\big(\|\bd{A}'\bd{h}\|^2\|\bd{B}'\bd{h}\|^2\big)\\
&= & \frac{1}{m^2}\cdot E\big\{[\bd{h}'(\bd{A}\bd{A}')\bd{h}]\cdot [\bd{h}'(\bd{B}\bd{B}')\bd{h}]\big\}\\
& = & \frac{1}{m^3(m+2)}\big[2\,\mbox{tr}(\bd{A}\bd{A}'\bd{B}\bd{B}')+
\mbox{tr}(\bd{A}\bd{A}')\cdot\mbox{tr}(\bd{B}\bd{B}') \big].
\end{eqnarray*}
Combing all of the above equalities, we have
\begin{eqnarray*}
&& \mbox{Cov}\big[(\bd{h}'\bd{A}\bd{h}_1)^2, (\bd{h}'\bd{B}\bd{h}_2)^2\big]\\
&=&\frac{1}{m^3(m+2)}\big[2\,\mbox{tr}(\bd{A}\bd{A}'\bd{B}\bd{B}')+
\mbox{tr}(\bd{A}\bd{A}')\cdot\mbox{tr}(\bd{B}\bd{B}') \big]-\frac{1}{m^4}\,\mbox{tr}(\bd{A}\bd{A}')\cdot\mbox{tr}(\bd{B}\bd{B}')\\
& = & \frac{2}{m^3(m+2)}\cdot \mbox{tr}(\bd{A}\bd{A}'\bd{B}\bd{B}')-\frac{2}{m^4(m+2)}\,
\mbox{tr}(\bd{A}\bd{A}')\cdot\mbox{tr}(\bd{B}\bd{B}').
\end{eqnarray*}
The proof is completed. \hfill$\Box$

\smallskip

\noindent\textbf{Proof of Lemma \ref{use_tea}}. From Lemma \ref{long_for} and the fact that $\bd{P}_i\bd{\epsilon}_i/\|\bd{P}_i\bd{\epsilon}_i\|=\bd{U}_i\bd{s}_i$,  $\{\hat{\rho}_{ij};\, 1\leq i< j \leq N\}$ has the same distribution as that of
\begin{eqnarray*}
\{\bd{s}_i'\bd{U}_i'\bd{U}_j\bd{s}_j,\, 1\leq i < j \leq N \}.
\end{eqnarray*}
We will use this fact repeatedly to prove the results next.

By independence, $E[\hat{\rho}_{ij}|\bd{s}_i]=\bd{s}_i'\bd{U}_i'\bd{U}_j\cdot E\bd{s}_j=0$ for $i \ne j$. Hence $E\hat{\rho}_{ij}=0$. Since $\bd{s}_i'\bd{U}_i'\bd{U}_j\bd{s}_j=(\bd{s}_i'\bd{U}_i'\bd{U}_j\bd{s}_j)'
=\bd{s}_j'\bd{U}_j'\bd{U}_i\bd{s}_i\in \mathbb{R}$, we have
\begin{align}\lbl{old_bear}
\hat{\rho}_{ij}^2=\bd{s}_j'(\bd{U}_j'\bd{U}_i\bd{s}_i
\bd{s}_i'\bd{U}_i'\bd{U}_j)\bd{s}_j.
\end{align}
Let $\bd{B}=\bd{U}_j'\bd{U}_i\bd{s}_i
\bd{s}_i'\bd{U}_i'\bd{U}_j.$ Conditioning on $\bd{s}_i$, we see from independence that
\begin{align}\lbl{patience}
E\big[\bd{s}_j'(\bd{U}_j'\bd{U}_i\bd{s}_i
\bd{s}_i'\bd{U}_i'\bd{U}_j)\bd{s}_j|\bd{s}_i\big]
= E(\bd{s}_j'\bd{B}\bd{s}_j|\bd{s}_i)
= \frac{1}{m}\,\mbox{tr}(\bd{B})
\end{align}
by Lemma \ref{sun_quiet}.
By \eqref{Science}, $\mbox{tr}(\bd{B})=
\bd{s}_i'(\bd{U}_i'\bd{U}_j\bd{U}_j'\bd{U}_i)\bd{s}_i.$ The above assertions conclude that
\begin{align}\lbl{no_worry}
E\big[\bd{s}_j'(\bd{U}_j'\bd{U}_i\bd{s}_i
\bd{s}_i'\bd{U}_i'\bd{U}_j)\bd{s}_j|\bd{s}_i\big]=\frac{1}{m}\cdot \bd{s}_i'\bd{M}_{ij}\bd{s}_i
\end{align}
by the notation $\bd{M}_{ij}=\bd{U}_i'\bd{U}_j\bd{U}_j'\bd{U}_i.$ Combining \eqref{old_bear} and \eqref{no_worry} together, we obtain $E[\hat{\rho}_{ij}^2|\bd{s}_i]=\frac{1}{m}\cdot \bd{s}_i'\bd{M}_{ij}\bd{s}_i$. Now, taking a further expectation, we have from  Lemma \ref{sun_quiet} again that  $E\hat{\rho}_{ij}^2=\frac{1}{m^2}\cdot\mbox{tr}(\bd{M}_{ij})$.
By \eqref{happy_moon},  $\bd{U}_{i}\bd{U}_{i}'=\bd{P}_i$. By \eqref{Science},
\begin{align}\lbl{sun_dark}
\mbox{tr}(\bd{M}_{ij})=\mbox{tr}(\bd{U}_i\bd{U}_i'\bd{U}_j\bd{U}_j')
=\mbox{tr}(\bd{P}_i\bd{P}_j).
\end{align}
We get the second conclusion from (ii).  \hfill$\Box$

\smallskip

In the following we will use the conditional variance $\mbox{Var}(\xi_2|\xi_1)$, which is defined by $E(\xi_2^2|\xi_1)-[E(\xi_2|\xi_1)]^2$ for any random variables $\xi_1$ and $\xi_2.$

\smallskip

\noindent\textbf{Proof of Lemma \ref{good_tie}}. (i) Review \eqref{old_bear} and  the notation $\bd{B}=\bd{U}_j'\bd{U}_i\bd{s}_i
\bd{s}_i'\bd{U}_i'\bd{U}_j$.
Then
\begin{align}\lbl{sea_salt}
(\hat{\rho}_{ij})^4=(\bd{s}_j'\bd{B}\bd{s}_j)^2.
\end{align}
Since $\bd{s}_i
\bd{s}_i'$ is a rank-one matrix, we know the rank of $\bd{B}$ is no more than $1.$ As a consequence $\mbox{tr}(\bd{B}^2)= [\mbox{tr}(\bd{B})]^2=(\bd{s}_i'\bd{M}_{ij}\bd{s}_i)^2$ since $\mbox{tr}(\bd{B})=\bd{s}_i'\bd{M}_{ij}\bd{s}_i$ by \eqref{Science}. Use independence and Lemma \ref{sun_quiet} to yield
\begin{align}\lbl{remedy}
E\big[(\bd{s}_j'\bd{B}\bd{s}_j)^2\big|\bd{s}_i\big]
=\frac{3}{m(m+2)}\cdot\big(\bd{s}_i'\bd{M}_{ij}\bd{s}_i\big)^2,
\end{align}
and hence
\begin{align}\lbl{red_green}
E\big[(\hat{\rho}_{ij})^4\big|\bd{s}_i\big]
= \frac{3}{m(m+2)}\cdot(\bd{s}_i'\bd{M}_{ij}\bd{s}_i\big)^2
\end{align}
by \eqref{sea_salt}. We obtain (i).

(ii) Taking expectations for both sides of \eqref{red_green}, we get from Lemma \ref{sun_quiet}(ii) that
\begin{eqnarray*}
E\big[(\hat{\rho}_{ij})^4\big]&=&\frac{3}{m(m+2)}\cdot E(\bd{s}_i'\bd{M}_{ij}\bd{s}_i\big)^2\\
&= & \frac{3}{m^2(m+2)^2}\cdot \big\{2\,\mbox{tr}(\bd{M}_{ij}^2)+   [\mbox{tr}(\bd{M}_{ij})]^2\big\}.
\end{eqnarray*}
By \eqref{sun_dark}, $\mbox{tr}(\bd{M}_{ij})=\mbox{tr}(\bd{P}_i\bd{P}_j)$. Also, from \eqref{happy_moon} and \eqref{Science},
\begin{eqnarray*}
\mbox{tr}(\bd{M}_{ij}^2)=\mbox{tr}(\bd{U}_i'\bd{U}_j\bd{U}_j'
\bd{U}_i\bd{U}_i'\bd{U}_j\bd{U}_j'\bd{U}_i)=\mbox{tr}[(\bd{P}_i\bd{P}_j)^2].
\end{eqnarray*}
We have proved (ii).

(iii) Notice
\begin{eqnarray*}
\mbox{Var}(\hat{\rho}_{ij}^2|\bd{s}_i)&=&
E\big[(\hat{\rho}_{ij})^4|\bd{s}_i\big] - \big[E(\hat{\rho}_{ij}^2|\bd{s}_i)\big]^2\\
& = & \frac{3}{m(m+2)}\cdot\big(\bd{s}_i'\bd{M}_{ij}\bd{s}_i\big)^2
-\Big[\frac{1}{m}\cdot \bd{s}_i'\bd{M}_{ij}\bd{s}_i\Big]^2\\
&= & \frac{2(m-1)}{m^2(m+2)}\cdot\big(\bd{s}_i'\bd{M}_{ij}\bd{s}_i\big)^2
\end{eqnarray*}
by (i) proved above and (ii) from Lemma \ref{use_tea}. \hfill$\Box$

(iv) By (i) of Lemma \ref{use_tea} and (ii) proved above,
\begin{eqnarray*}
&& \mbox{Var}(\hat{\rho}_{ij}^2)\\
&=& \frac{3}{m^2(m+2)^2}\cdot \big\{2\,\mbox{tr}[(\bd{P}_i\bd{P}_j)^2]+   [\mbox{tr}(\bd{P}_i\bd{P}_j)]^2\big\} - \Big[\frac{1}{m^2}\cdot\mbox{tr}(\bd{P}_i\bd{P}_j)\Big]^2\\
& = & \frac{6}{m^2(m+2)^2}\cdot \mbox{tr}[(\bd{P}_i\bd{P}_j)^2]+\frac{2(m^2-2m-2)}{m^4(m+2)^2}\cdot [\mbox{tr}(\bd{P}_i\bd{P}_j)]^2.
\end{eqnarray*}
We finish the proof. \hfill$\Box$

\subsection{Proofs of auxiliary results in Section \ref{huangjingou}}\lbl{huangjingou_hei}

Review the interpretation of constant $C$ before the statement of Lemma \ref{Love_sound}.


\smallskip

\noindent\textbf{Proof of Lemma \ref{Love_sound}}. Recall \eqref{park_1}. Set $\bd{A}_i=\bd{x}_i(\bd{x}_i'\bd{x}_i)^{-1}\bd{x}_i'$ for $1\leq i \leq N.$ Then $\bd{A}_i$ is a $T\times T$ idempotent matrix with rank $p$ and $\mbox{tr}(\bd{A}_i)=p$ for each $i$. Since  $\bd{P}_i=\bd{I}_T-\bd{A}_i$, we see
\begin{align}\lbl{head_toss}
\bd{P}_i\bd{P}_j=\bd{I}_T+\bd{B}_{ij}
\end{align}
where $\bd{B}_{ij}:=\bd{A}_i\bd{A}_j-\bd{A}_i-\bd{A}_j.$ By Lemma \ref{brother_cat}, %
%
\begin{align}\lbl{acuby}
\mbox{tr}(\bd{F}_1\bd{F}_2)\geq 0
\end{align}
for  any non-negative definite matrices $\bd{F}_1$ and $\bd{F}_2$. As a result, $\mbox{tr}(\bd{A}_i\bd{A}_j) \geq 0$.
Easily, $\mbox{tr}(\bd{A}_i\bd{A}_j)\leq p$ by Lemma \ref{trivial}. Thus,
\begin{align}\lbl{trace_ralation1}
-2p\leq \mbox{tr}(\bd{B}_{ij})\leq -p.
\end{align}
Expand $[\bd{A}_i\bd{A}_j-(\bd{A}_i+\bd{A}_j)]^2$ and use \eqref{Science} to see
\begin{align}\lbl{plane_girl}
\mbox{tr}(\bd{B}_{ij}^2)=\mbox{tr}(\bd{A}_i\bd{A}_j\bd{A}_i\bd{A}_j)
-2\mbox{tr}(\bd{A}_i\bd{A}_j\bd{A}_i)-2\mbox{tr}(\bd{A}_j\bd{A}_i\bd{A}_j)+
\mbox{tr}((\bd{A}_i+\bd{A}_j)^2).
\end{align}
By \eqref{acuby}, $\mbox{tr}(\bd{A}_i\bd{A}_j\bd{A}_i\bd{A}_j)=
\mbox{tr}[\bd{A}_i(\bd{A}_j\bd{A}_i\bd{A}_j)]\geq 0$ because $\bd{A}_j\bd{A}_i\bd{A}_j$ is non-negative definite. So each trace in \eqref{plane_girl} is non-negative. Also, $\mbox{tr}((\bd{A}_i+\bd{A}_j)^2)=\mbox{tr}(\bd{A}_i^2)+
2\mbox{tr}(\bd{A}_i\bd{A}_j)+\mbox{tr}(\bd{A}_j^2)\leq 4p$ by Lemma  \ref{trivial}. Observe
\begin{center}
$\max\{\mbox{tr}(\bd{A}_i\bd{A}_j\bd{A}_i), \mbox{tr}(\bd{A}_j\bd{A}_i\bd{A}_j)\} \leq \mbox{tr}(\bd{A}_i\bd{A}_j)$,
\end{center}
$\mbox{tr}(\bd{A}_i\bd{A}_j\bd{A}_i) \leq \mbox{tr}(\bd{A}_i^2)$ and $\mbox{tr}(\bd{A}_j\bd{A}_i\bd{A}_j) \leq \mbox{tr}(\bd{A}_j^2)$ by Lemma  \ref{trivial}. Therefore,
\begin{eqnarray*}
2\mbox{tr}(\bd{A}_i\bd{A}_j\bd{A}_i)+2\mbox{tr}(\bd{A}_j\bd{A}_i\bd{A}_j)\leq \mbox{tr}((\bd{A}_i+\bd{A}_j)^2) \leq 4p.
\end{eqnarray*}
It follows that
\begin{align}\lbl{trace_ralation2}
0\leq \mbox{tr}(\bd{B}_{ij}^2) \leq 5p.
\end{align}
With the above preparation, we now derive the conclusions. In fact, from \eqref{head_toss},
\begin{eqnarray*}
[\mbox{tr}(\bd{P}_i\bd{P}_j)]^2
= T^2 + 2\,\mbox{tr}(\bd{B}_{ij})\cdot T+[\mbox{tr}(\bd{B}_{ij})]^2.
\end{eqnarray*}
This implies (i) by \eqref{trace_ralation1}. Now, from \eqref{head_toss} again,
\begin{align}\lbl{Rain_5}
\mbox{tr}((\bd{P}_i\bd{P}_j)^2)=T+2\,\mbox{tr}(\bd{B}_{ij})+ \mbox{tr}(\bd{B}_{ij}^2).
\end{align}
Then (ii) follows from \eqref{trace_ralation1} and \eqref{trace_ralation2}. Let us show the remaining two claims next.

By the definition of $\bd{B}_{ij}$ and the notation $q=|S|$,
\begin{align}\lbl{star_butter}
\bd{P}_S\bd{P}_j
=& \Big[q\bd{I}_T-\sum_{i\in S}\bd{A}_i\Big]\cdot (\bd{I}_T-\bd{A}_j)\nonumber\\
= & q\bd{I}_T-q\bd{A}_j-\Big[\sum_{i\in S}\bd{A}_i\Big]+
\sum_{i\in S}\bd{A}_i\bd{A}_j\nonumber\\
= & q\bd{I}_T+\sum_{i\in S}\bd{B}_{ij}.
\end{align}
Hence,
\begin{eqnarray*}
\mbox{tr}(\bd{P}_S\bd{P}_j)=
qT+\sum_{i\in S}\mbox{tr}(\bd{B}_{ij}).
\end{eqnarray*}
The inequality from \eqref{trace_ralation1} implies that $\sum_{i\in S}\mbox{tr}(\bd{B}_{ij})$ is between $-2pq$ and $-pq$. This leads to (iii) by a trivial equality $(x+y)^2=x^2+2xy+y^2$ for all $x, y\in \mathbb{R}.$

To get (iv), we start from \eqref{star_butter} again such that
\begin{eqnarray*}
\frac{1}{q^2}(\bd{P}_S\bd{P}_j)^2
=\bd{I}_T+
\frac{2}{q}\Big(\sum_{i\in S}\bd{B}_{ij}\Big)
+\Big(\frac{1}{q}\sum_{i\in S}\bd{B}_{ij}\Big)^2
\end{eqnarray*}
Then
\begin{align}\lbl{story_tell}
\frac{1}{q^2}\mbox{tr}((\bd{P}_S\bd{P}_j)^2)
=&T+\frac{2}{q}\Big[\sum_{i\in S}\mbox{tr}(\bd{B}_{ij})\Big]+
\mbox{tr}\Big[\Big(\frac{1}{q}\sum_{i\in S}\bd{B}_{ij}\Big)^2\Big] \nonumber\\
:=& T +C_{ij}.
\end{align}
By the Cauchy-Schwartz inequality, for any $T\times T$ matrix $\bd{M}=(m_{ij})_{T\times T}$, we have $|\mbox{tr}(\bd{M}^2)|=|\sum_{1\leq i, j\leq T}m_{ij}m_{ji}|\leq \sum_{1\leq i, j\leq T}m_{ij}^2=\|\bd{M}\|_F^2$, where $\|\bd{M}\|_F:=\sqrt{\mbox{tr}(\bd{M}\bd{M}')}$ is  the Frobenius norm of $\bd{M}$. By the triangle inequality and then Lemma \ref{trivial}, $\|\bd{B}_{ij}\|_F \leq \|\bd{A}_i\bd{A}_j\|_F+\|\bd{A}_{i}\|_F+\|\bd{A}_{j}\|_F\leq 3\sqrt{p}$. It follows that
\begin{align}\lbl{chuckle}
\Big|\mbox{tr}\Big[\Big(\frac{1}{q}\sum_{i\in S}\bd{B}_{ij}\Big)^2\Big]\Big|^{1/2}
\leq  \Big\|\frac{1}{q}\sum_{i\in S}\bd{B}_{ij}\Big\|_F
\leq  \frac{1}{q}\sum_{i\in S}\|\bd{B}_{ij}\|_F\leq 3\sqrt{p}.
\end{align}
This and \eqref{trace_ralation1} conclude that $-13p\leq C_{ij}\leq 7p$ for all $i,j$. We then get (iv) from \eqref{story_tell}.

Now we prove (v). Obviously (i) and (ii) still hold if symbol``$T$" is replaced by ``$m$". On the other hand,  by the triangle inequality and the facts $T=m+p$ and $T^2-m^2=2mp+p^2$,
\begin{eqnarray*}
&&\frac{1}{mq^2}\cdot \big|[\mbox{tr}(\bd{P}_S\bd{P}_j)]^2-m^2q^2\big|\\
&\leq & \frac{1}{mq^2}\cdot
\big\{\big|[\mbox{tr}(\bd{P}_S\bd{P}_j)]^2-T^2q^2\big|+ (2mp+p^2)q^2\big\}\\
& = & \Big(1+\frac{p}{m}\Big)\cdot\frac{1}{Tq^2}\cdot\big|[\mbox{tr}(\bd{P}_S\bd{P}_j)]^2-T^2q^2\big|+ 2p+\frac{p^2}{m}.
\end{eqnarray*}
Since $p$ is fixed,  (iii) is also true if ``$T$" is replaced by ``$m$". The remaining part of (v) is  obtained similarly.

The constant $K$ is taken to be the maximum of the five bounds in (i)-(v). \hfill$\Box$

\smallskip

\noindent\textbf{Proof of Lemma \ref{shadow_1}}.  By \eqref{happy_moon}, $\bd{U}_i\bd{U}_i'=\bd{P}_i$. Use this fact and \eqref{Science} to see
\beaa
\mbox{tr}(\bd{M}_{ij})=\mbox{tr}(\bd{P}_i\bd{P}_j)\ \ \ \mbox{and}\ \ \ \mbox{tr}(\bd{M}_{ij}^2)=\mbox{tr}\big[(\bd{P}_{i}\bd{P}_{j})^2\big].
\eeaa
Let $\xi=\bd{e}'\bd{M}_{ij}\bd{e}$. By Lemma \ref{sun_quiet}, $E\xi=\frac{1}{m}\mbox{tr}(\bd{P}_i\bd{P}_j)$ and
\beaa
\mbox{Var}(\xi)=\frac{2}{m(m+2)}\Big\{\mbox{tr}\big[(\bd{P}_{i}\bd{P}_{j})^2\big]
-\frac{1}{m}\cdot \big[\mbox{tr}(\bd{P}_i\bd{P}_j)\big]^2\Big\}.
\eeaa
By taking $\alpha=4$ in Lemma \ref{dust_bone}, we get
\beaa
E\big[(\xi-E\xi)^4\big]\leq \frac{C}{m^4}\Big\{\mbox{tr}\big[(\bd{P}_{i}\bd{P}_{j})^2\big]
-\frac{1}{m}\cdot \big[\mbox{tr}(\bd{P}_i\bd{P}_j)\big]^2\Big\}^2,
\eeaa
as $m\geq 4\alpha+1$, that is,  $T\geq p+17$ since $m=T-p$ and $\alpha=4$. Notice $\mbox{rank}(\bd{P}_i\bd{P}_j)\leq \mbox{rank}(\bd{P}_i)=m$.
By Lemma \ref{brother_cat}, Lemma \ref{Love_sound}(v) and the triangle inequality,
\beaa
0\leq \mbox{tr}\big[(\bd{P}_{i}\bd{P}_{j})^2\big]
-\frac{1}{m}\cdot \big[\mbox{tr}(\bd{P}_i\bd{P}_j)\big]^2\leq C.
\eeaa
This says that
\beaa
\max\Big\{\mbox{Var}(\xi),\, \sqrt{E\big[(\xi-E\xi)^4\big]}\Big\} \leq \frac{C}{m^2}.
\eeaa
From Lemma \ref{trivial}, $E\xi\leq 1$. Recall $T=m+p$. By (i) of Lemma \ref{Love_sound}, there exists a constant $K>0$ such that
\beaa
E\xi=\frac{1}{m}\mbox{tr}(\bd{P}_i\bd{P}_j)\geq \frac{1}{m}\sqrt{T^2-TK}\geq \frac{1}{2}
\eeaa
as $N \geq C$ since $T=T_N\to \infty$ as $N\to \infty$, where $C>0$ is a constant free of $N$, $p$ and $T$. By using the above two inequalities and Lemma \ref{mac_pencil}, we see
\beaa
&& E[ (\xi^2-E\xi^2)^2 ]\\
&\leq & 16^2\cdot \Big[4\cdot {\rm Var}(\xi)^{2}+\sqrt{E(|\xi-E\xi|^{4})}\Big]\cdot
\Big[1+4\cdot {\rm Var}(\xi)^{2}+\sqrt{E(|\xi-E\xi|^{4})}\Big]\\
& \leq & \frac{C}{m^2}.
\eeaa
The proof is completed. \hfill$\square$

\smallskip

\noindent\textbf{Proof of Lemma \ref{snow_sound}}. The second expression of $X_j$ follows from Lemma \ref{use_tea}.
%
Now we start to compute $E(X_j^2)$.

Evidently, $\bd{s}_i'\bd{U}_i'\bd{U}_j\bd{s}_j=\bd{s}_j'\bd{U}_j'\bd{U}_i\bd{s}_i
\in \mathbb{R}$. We have from  \eqref{old_bear} that
\begin{eqnarray*}
\hat{\rho}_{ij}^2=\bd{s}_j'(\bd{U}_j'\bd{U}_i\bd{s}_i
\bd{s}_i'\bd{U}_i'\bd{U}_j)\bd{s}_j=\bd{s}_i'(\bd{U}_i'\bd{U}_j\bd{s}_j
\bd{s}_j'\bd{U}_j'\bd{U}_i)\bd{s}_i.
\end{eqnarray*}
Let $\bd{H}_{ij}=\bd{U}_i'\bd{U}_j\bd{s}_j
\bd{s}_j'\bd{U}_j'\bd{U}_i$. Recall $\bd{M}_{ij}=\bd{U}_i'\bd{U}_j\bd{U}_j'\bd{U}_i$.  Write
\begin{align}\lbl{soon_depart}
\frac{1}{T}X_j
= \sum_{i=1}^{j-1}\big[\bd{s}_i'\bd{H}_{ij}\bd{s}_i
-    \frac{1}{m}
\bd{s}_i'\bd{M}_{ij}\bd{s}_i\big]
\end{align}
for $2\leq j \leq N.$ Given $\bd{s}_j$, the conditional mean of the term in the sum above is equal to
\begin{eqnarray*}
& & E\big(\bd{s}_i'\bd{H}_{ij}\bd{s}_i
\big|\bd{s}_j\big)-\frac{1}{m}
E\big(\bd{s}_i'\bd{M}_{ij}\bd{s}_i\big)\\
&= & \frac{1}{m}\cdot \mbox{tr}(\bd{H}_{ij})-\frac{1}{m^2}\cdot \mbox{tr}(\bd{M}_{ij})\\
& = & \frac{\bd{s}_j'\bd{M}_{ji}\bd{s}_j}{m}-\frac{1}{m^2}\cdot \mbox{tr}(\bd{P}_i\bd{P}_j)
\end{eqnarray*}
by Lemma \ref{sun_quiet} and \eqref{happy_moon}. Observe that, given $\bd{s}_j$, the terms in the sum from \eqref{soon_depart} are independent. Also, it is true that $E(\psi+c)^2=\mbox{Var}(\psi)+c^2$ for any random variable $\psi$ with mean zero and constant $c$.  Thus,
\begin{align}\lbl{light_play}
 E\Big(\frac{1}{T^2}X_j^2\Big|\bd{s}_j\Big)
 = & \sum_{i=1}^{j-1} \mbox{Var}\Big[\big(\bd{s}_i'\bd{H}_{ij}\bd{s}_i
-    \frac{1}{m}
\bd{s}_i'\bd{M}_{ij}\bd{s}_i\big)\big|\bd{s}_j\Big]\nonumber\\
&+ \Big[\sum_{i=1}^{j-1}\Big(\frac{\bd{s}_j'\bd{M}_{ji}\bd{s}_j}{m}-\frac{1}{m^2}\cdot \mbox{tr}(\bd{P}_i\bd{P}_j)\Big)\Big]^2.
\end{align}
In what follows we exam the last two terms carefully. Write
\begin{eqnarray*}
\bd{s}_i'\bd{H}_{ij}\bd{s}_i-
\frac{1}{m}
\bd{s}_i'\bd{M}_{ij}\bd{s}_i
=\bd{s}_i'\bd{D}_{ij}\bd{s}_i
\end{eqnarray*}
where
\begin{eqnarray*}
\bd{D}_{ij}:=\bd{H}_{ij}- \frac{\bd{M}_{ij}}{m}.
\end{eqnarray*}
Define
\begin{eqnarray*}
\Upsilon_{ij}=\mbox{Var}\big[\big(\bd{s}_i'\bd{H}_{ij}\bd{s}_i
 -    \frac{1}{m}
\bd{s}_i'\bd{M}_{ij}\bd{s}_i\big)\big|\bd{s}_j\big].
\end{eqnarray*}
Therefore, we get from Lemma \ref{sun_quiet} that
\begin{align}\lbl{mean_hello}
\Upsilon_{ij}=\frac{2}{m(m+2)}\cdot \mbox{tr}(\bd{D}_{ij}^2)-\frac{2}{m^2(m+2)} [\mbox{tr}(\bd{D}_{ij})]^2.
\end{align}
First,
\begin{eqnarray*}
\mbox{tr}(\bd{D}_{ij})=\bd{s}_j'\bd{M}_{ji}\bd{s}_j
-\frac{1}{m}\cdot\mbox{tr}(\bd{P}_i\bd{P}_j),
\end{eqnarray*}
hence
\begin{align}\lbl{sea_sea}
[\mbox{tr}(\bd{D}_{ij})]^2=(\bd{s}_j'\bd{M}_{ji}\bd{s}_j)^2
-\frac{2\,\mbox{tr}(\bd{P}_i\bd{P}_j)}{m}\cdot \bd{s}_j'\bd{M}_{ji}\bd{s}_j+\frac{1}{m^2}
\cdot[\mbox{tr}(\bd{P}_i\bd{P}_j)]^2.
\end{align}
Second,
\begin{eqnarray*}
\mbox{tr}(\bd{D}_{ij}^2)=\mbox{tr}(\bd{H}_{ij}^2)-2\cdot \frac{\mbox{tr}(\bd{H}_{ij}\bd{M}_{ij})}{m} + \frac{\mbox{tr}(\bd{M}_{ij}^2)}{m^2}.
\end{eqnarray*}
Observe the rank of $\bd{s}_j\bd{s}_j'$ is at most one, since   $\bd{H}_{ij}=\bd{U}_i'\bd{U}_j(\bd{s}_j
\bd{s}_j')\bd{U}_j'\bd{U}_i$, we know  $\mbox{rank}(\bd{H}_{ij})\leq 1$. As a consequence,
\begin{eqnarray*}
\mbox{tr}(\bd{H}_{ij}^2)=[\mbox{tr}(\bd{H}_{ij})]^2=[\bd{s}_j'\bd{M}_{ji}\bd{s}_j]^2.
\end{eqnarray*}
Now, by the definition of $\bd{M}_{ij}$ and the fact $\bd{U}_{i}\bd{U}_{i}'=\bd{P}_i$ in  \eqref{happy_moon},
\begin{align}
 & \mbox{tr}(\bd{H}_{ij}\bd{M}_{ij})=\mbox{tr}(\bd{U}_i'\bd{U}_j\bd{s}_j
\bd{s}_j'\bd{U}_j'\bd{U}_i\bd{U}_i'\bd{U}_j\bd{U}_j'\bd{U}_i)
=\bd{s}_j'\bd{M}_{ji}^2\bd{s}_j;\nonumber\\
 & \mbox{tr}(\bd{M}_{ij}^2)
=\mbox{tr}(\bd{U}_i'\bd{U}_j\bd{U}_j'\bd{U}_i\bd{U}_i'\bd{U}_j\bd{U}_j'\bd{U}_i)
=\mbox{tr}((\bd{P}_i\bd{P}_j)^2),\lbl{mo_fact}
\end{align}
where \eqref{Science} is used above. Combining the above identities to see
\begin{eqnarray*}
\mbox{tr}(\bd{D}_{ij}^2)=[\bd{s}_j'\bd{M}_{ji}\bd{s}_j]^2
-2\cdot \frac{\bd{s}_j'\bd{M}_{ji}^2\bd{s}_j}{m} + \frac{\mbox{tr}((\bd{P}_i\bd{P}_j)^2)}{m^2}.
\end{eqnarray*}
This together with \eqref{mean_hello} and \eqref{sea_sea} implies that
\begin{eqnarray*}
&& \Upsilon_{ij}\\
& =& \frac{2}{m(m+2)}\big[(\bd{s}_j'\bd{M}_{ji}\bd{s}_j)^2
-2\cdot \frac{\bd{s}_j'\bd{M}_{ji}^2\bd{s}_j}{m} + \frac{\mbox{tr}((\bd{P}_i\bd{P}_j)^2)}{m^2}\Big] -\nonumber\\
& &\nonumber\\
& & \frac{2}{m^2(m+2)}\Big\{(\bd{s}_j'\bd{M}_{ji}\bd{s}_j)^2
-2\,\frac{\mbox{tr}(\bd{P}_i\bd{P}_j)}{m}\cdot (\bd{s}_j'\bd{M}_{ji}\bd{s}_j)+\frac{1}{m^2}
\cdot[\mbox{tr}(\bd{P}_i\bd{P}_j)]^2\Big\}.
\end{eqnarray*}
By a trivial sorting, we obtain
\begin{align}
&\Upsilon_{ij} \nonumber\\
=&\frac{2m-2}{m^2(m+2)}\cdot(\bd{s}_j'\bd{M}_{ji}\bd{s}_j)^2
-\frac{4}{m^2(m+2)}\cdot \bd{s}_j'\bd{M}_{ji}^2\bd{s}_j
+\frac{4\,\mbox{tr}(\bd{P}_i\bd{P}_j)}{m^3(m+2)}\cdot \bd{s}_j'
\bd{M}_{ji}\bd{s}_j + \nonumber\\
&\frac{2}{m^3(m+2)}\cdot \mbox{tr}((\bd{P}_i\bd{P}_j)^2)-\frac{2}{m^4(m+2)}\cdot[\mbox{tr}(\bd{P}_i\bd{P}_j)]^2. \lbl{paper_window}
\end{align}
Now we analyze the expectation of each term above in order to compute the mean of the conditional variance. It is easy to check
\begin{align}\lbl{liu_man}
\mbox{tr}(\bd{M}_{ji}\bd{M}_{jk})=\mbox{tr}(\bd{P}_{i}\bd{P}_{j}\bd{P}_{k}\bd{P}_{j})
\end{align}
 for any $1\leq i, j, k\leq N$. Now, by Lemma \ref{sun_quiet},
\begin{eqnarray*}
E(\bd{s}_j'\bd{M}_{ji}\bd{s}_j)^2
& = & \frac{1}{m(m+2)}\cdot \big\{2\,\mbox{tr}(\bd{M}_{ji}^2)+[\mbox{tr}(\bd{M}_{ji})]^2\big\}\\
&= & \frac{1}{m(m+2)}\cdot \big\{2\,\mbox{tr}((\bd{P}_i\bd{P}_j)^2)+[\mbox{tr}(\bd{P}_i\bd{P}_j)]^2\big\}
\end{eqnarray*}
since $\mbox{tr}(\bd{M}_{ji})=\mbox{tr}(\bd{P}_i\bd{P}_j)$. By Lemma \ref{sun_quiet} again,
\begin{eqnarray*}
& & E(\bd{s}_j'\bd{M}_{ji}^2\bd{s}_j)=\frac{1}{m}\cdot \mbox{tr}(\bd{M}_{ji}^2)=\frac{1}{m}\cdot \mbox{tr}((\bd{P}_i\bd{P}_j)^2);\\
& & E(\bd{s}_j'
\bd{M}_{ji}\bd{s}_j)
=\frac{1}{m}\cdot\mbox{tr}(\bd{M}_{ji})=\frac{1}{m}\cdot\mbox{tr}(\bd{P}_i\bd{P}_j).
\end{eqnarray*}
Take expectations for both sides of \eqref{paper_window} and use the above facts to see
\begin{align}\lbl{tedious_lion}
E\Upsilon_{ij} =& \mbox{tr}((\bd{P}_i\bd{P}_j)^2)\cdot \Big[\frac{4m-4}{m^3(m+2)^2}-\frac{4}{m^3(m+2)}+\frac{2}{m^3(m+2)}\Big]+ \nonumber\\
 & [\mbox{tr}(\bd{P}_i\bd{P}_j)]^2\cdot\Big[\frac{2m-2}{m^3(m+2)^2}
+\frac{4}{m^4(m+2)}-\frac{2}{m^4(m+2)}\Big] \nonumber\\
=& \mbox{tr}((\bd{P}_i\bd{P}_j)^2)\cdot\frac{2m-8}{m^3(m+2)^2}+ [\mbox{tr}(\bd{P}_i\bd{P}_j)]^2\cdot\frac{2m^2+4}{m^4(m+2)^2}.
\end{align}

Now we turn to study the mean of the last term from \eqref{light_play}. Write
\begin{eqnarray*}
\Xi_{ij}:&=&\sum_{i=1}^{j-1}\Big(\frac{\bd{s}_j'\bd{M}_{ji}\bd{s}_j}{m}
-\frac{1}{m^2}\cdot \mbox{tr}(\bd{P}_i\bd{P}_j)\Big)\\
& = & \frac{\bd{s}_j'\bd{M}_{j\blacktriangle}\bd{s}_j}{m}-\frac{1}{m^2}\cdot \mbox{tr}(\bd{P}_{j\blacktriangle}\bd{P}_j),
\end{eqnarray*}
for $2\leq j \leq N,$ where we define
\begin{align}\lbl{teacher_ego}
\bd{M}_{j\blacktriangle}=\sum_{i=1}^{j-1}\bd{M}_{ji}\ \ \mbox{and}\ \
\bd{P}_{j\blacktriangle}=\sum_{i=1}^{j-1}\bd{P}_i.
\end{align}
By using Lemma \ref{sun_quiet},
\begin{eqnarray*}
E\frac{\bd{s}_j'\bd{M}_{j\blacktriangle}\bd{s}_j}
{m}
=\frac{\mbox{tr}(\bd{M}_{j\blacktriangle})}{m^2}
=\frac{1}{m^2}\sum_{i=1}^{j-1}\mbox{tr}(\bd{M}_{ji}).
\end{eqnarray*}
Since $\mbox{tr}(\bd{M}_{ji})=\mbox{tr}(\bd{P}_i\bd{P}_j)$,
the above is equal to
\begin{eqnarray*}
\frac{1}{m^2}\sum_{i=1}^{j-1}\mbox{tr}(\bd{P}_i\bd{P}_j)=\frac{1}{m^2}\cdot \mbox{tr}(\bd{P}_{j\blacktriangle}\bd{P}_j).
\end{eqnarray*}
As a byproduct,
\begin{align}\lbl{force_seat}
\mbox{tr}(\bd{M}_{j\blacktriangle})=\mbox{tr}(\bd{P}_{j\blacktriangle}\bd{P}_j).
\end{align}
Therefore
\begin{align}\lbl{yellow_fa}
 E(\Xi_{ij}^2)
=&\frac{1}{m^2}\cdot\mbox{Var}\big(\bd{s}_j'\bd{M}_{j\blacktriangle}\bd{s}_j\big)\nonumber\\
 = & \frac{2}{m^3(m+2)}\cdot\Big\{\mbox{tr}(\bd{M}_{j\blacktriangle}^2)-\frac{1}{m}\cdot \big[\mbox{tr}(\bd{M}_{j\blacktriangle})\big]^2\Big\}
\end{align}
by Lemma \ref{sun_quiet}. 
Note that
\begin{align}\lbl{stay_tone}
\mbox{tr}(\bd{M}_{j\blacktriangle}^2)=\mbox{tr}\Big[\Big(\sum_{i=1}^{j-1}\bd{M}_{ji}\Big)^2\Big]
=\sum_{1\leq i, k \leq j-1}\mbox{tr}(\bd{P}_{i}\bd{P}_{j}\bd{P}_{k}\bd{P}_{j})
=  \mbox{tr}\big(\big(\bd{P}_{j\blacktriangle}\bd{P}_{j}\big)^2\big)
\end{align}
by \eqref{liu_man}. This, \eqref{force_seat} and \eqref{yellow_fa} conclude
\begin{align}\lbl{smile_butt}
E(\Xi_{ij}^2)=\frac{2}{m^3(m+2)}
\cdot\Big\{\mbox{tr}\big((\bd{P}_{j\blacktriangle}\bd{P}_{j})^2\big)-
\frac{1}{m}\cdot\big[\mbox{tr}(\bd{P}_{j\blacktriangle}\bd{P}_{j})\big]^2\Big\}.
\end{align}
Review the notations $\Xi_{ij}$ and $\Upsilon_{ij}$, the conclusion follows from \eqref{light_play}, \eqref{tedious_lion} and \eqref{smile_butt}.\hfill$\Box$

\smallskip

\noindent\textbf{Proof of Lemma \ref{kiss_rain}}. By Lemma \ref{snow_sound},
\begin{eqnarray*}
 &&\frac{1}{T^2}E(X_j^2)\\
 &=&\frac{2m-8}{m^3(m+2)^2}\sum_{i=1}^{j-1}\mbox{tr}((\bd{P}_i\bd{P}_j)^2)+ \frac{2m^2+4}{m^4(m+2)^2}\sum_{i=1}^{j-1}[\mbox{tr}(\bd{P}_i\bd{P}_j)]^2+\\
&&\frac{2}{m^3(m+2)}
\Big\{\mbox{tr}\big((\bd{P}_{j\blacktriangle}\bd{P}_{j})^2\big)-
\frac{1}{m}\cdot\big[\mbox{tr}(\bd{P}_{j\blacktriangle}\bd{P}_{j})\big]^2\Big\}
\end{eqnarray*}
for $2\leq j \leq N.$ We next analyze the above three terms.

Review $m=T-p.$ From Lemma \ref{Love_sound}(v), there exists a constant $K$ not depending on $T$ or $N$ such that
\begin{eqnarray*}
&& \frac{1}{mj^2}\cdot \big|[\mbox{tr}(\bd{P}_{j\blacktriangle}\bd{P}_j)]^2-m^2(j-1)^2\big|\leq K,\\
&& \frac{1}{j^2}\cdot \big|\mbox{tr}((\bd{P}_{j\blacktriangle}\bd{P}_j)^2)-m(j-1)^2\big|\leq K
\end{eqnarray*}
for all $2\leq j \leq N.$ By the triangle inequality,
\begin{align}\lbl{henan_xu}
\Big|\mbox{tr}\big((\bd{P}_{j\blacktriangle}\bd{P}_{j})^2\big)-
\frac{1}{m}\cdot\big[\mbox{tr}(\bd{P}_{j\blacktriangle}\bd{P}_{j})\big]^2\Big|\leq 2Kj^2
\end{align}
for $2\leq j \leq N.$ This, (i) and (ii) from Lemma \ref{Love_sound} imply
\begin{eqnarray*}
\frac{1}{T^2}E(X_j^2)&\leq & \frac{2m-8}{m^3(m+2)^2}(T+K)(j-1)\\
& &+ \frac{2m^2+4}{m^4(m+2)^2}(T^2+KT)(j-1)+
\frac{4Kj^2}{m^3(m+2)}.
\end{eqnarray*}
It follows that
\begin{eqnarray*}
&&\frac{1}{T^2}\sum_{j=2}^NE(X_j^2)\\
&\leq & \frac{2m-8}{m^3(m+2)^2}(T+K)\cdot \frac{1}{2}(N-1)N+\frac{2m^2+4}{m^4(m+2)^2}(T^2+KT)\cdot\\
& & \frac{1}{2}(N-1)N+\frac{4K}{m^3(m+2)}\cdot \frac{1}{6}N(N+1)(2N+1).
\end{eqnarray*}
Similarly, by the lower bound from \eqref{henan_xu},
\begin{eqnarray*}
&&\frac{1}{T^2}\sum_{j=2}^NE(X_j^2)\\
&\geq & \frac{2m-8}{m^3(m+2)^2}(T-K)\cdot \frac{1}{2}(N-1)N+\frac{2m^2+4}{m^4(m+2)^2}(T^2-KT)\cdot \\
& &\frac{1}{2}(N-1)N-\frac{4K}{m^3(m+2)}\cdot \Big[\frac{1}{6}N(N+1)(2N+1) -1\Big].
\end{eqnarray*}
Inspecting the above two bounds carefully, the dominating term is
\begin{eqnarray*}
\frac{2m^2+4}{m^4(m+2)^2}T^2\cdot \frac{1}{2}(N-1)N=\frac{T^2N^2}{m^4}(1+o(1))
\end{eqnarray*}
provided
\begin{eqnarray*}
\frac{4K}{m^3(m+2)}\cdot \frac{1}{6}N(N+1)(2N+1)=o\Big(\frac{T^2N^2}{m^4}\Big).
\end{eqnarray*}
This is equivalent to that $N=o(T^2).$ Therefore
\begin{eqnarray*}
\frac{1}{T^2}\sum_{j=2}^NE(X_j^2) =\frac{T^2N^2}{m^4}(1+o(1))
\end{eqnarray*}
as $N\to\infty$. Consequently,
\begin{eqnarray*}
\frac{1}{N^2}\sum_{j=2}^NE(X_j^2)\to 1
\end{eqnarray*}
as $N\to \infty$ since $m=T-p$ and $p$ is fixed. \hfill$\Box$

\smallskip

\noindent\textbf{Proof of Lemma \ref{cousin_basket}}. Set
\begin{eqnarray*}
R_j=\Big(\sum_{i=1}^{j-1}\bd{s}_i'\bd{M}_{ij}\bd{s}_i\Big)^2,\ \ 2\leq j \leq N.
\end{eqnarray*}
Then
\begin{align}\lbl{957}
{\rm Var}\Big[\sum_{j=2}^N\Big(\sum_{i=1}^{j-1}\bd{s}_i'\bd{M}_{ij}\bd{s}_i\Big)^2\Big] =& {\rm Var}\Big(\sum_{j=2}^NR_j\Big) \nonumber\\
 \leq & (N-1)\Big[\sum_{j=2}^N{\rm Var}(R_j)\Big]
\end{align}
by the convexity of function $f(x):=x^2$ for $x \in \mathbb{R}$. We next calculate ${\rm Var}(R_j)$ for each $j$.

Fix $2\leq j \leq N.$  For simplicity of notation,  set $\xi=\sum_{i=1}^{j-1}\bd{s}_i'\bd{M}_{ij}\bd{s}_i.$ Then $R_j=\xi^2$.
Note that $\mbox{tr}(\bd{M}_{ij})=\mbox{tr}(\bd{P}_i\bd{P}_j)$ by \eqref{sun_dark} and $\mbox{tr}(\bd{M}_{ij}^2)=\mbox{tr}[(\bd{P}_{i}\bd{P}_{j})^2]$ from \eqref{liu_man}. Then
\begin{align}\lbl{123}
E\xi
=\frac{1}{m}\sum_{i=1}^{j-1}\mbox{tr}(\bd{M}_{ij})
=\frac{1}{m}\sum_{i=1}^{j-1}\mbox{tr}(\bd{P}_i\bd{P}_j)= \frac{1}{m}\mbox{tr}(\bd{P}_{j\blacktriangle}\bd{P}_j),
\end{align}
where $\bd{P}_{j\blacktriangle}$ is defined in \eqref{teacher_ego}. By independence and  Lemma \ref{sun_quiet}, we obtain
\begin{align}\lbl{317}
\mbox{Var}(\xi)=&\sum_{i=1}^{j-1}\mbox{Var}
\big(\bd{s}_i'\bd{M}_{ij}\bd{s}_i\big) \nonumber\\
 = & \frac{2}{m(m+2)}\cdot\sum_{i=1}^{j-1}\Big\{\mbox{tr}(\bd{M}_{ij}^2)-\frac{1}{m}\cdot \big[\mbox{tr}(\bd{M}_{ij})\big]^2\Big\} \nonumber\\
 = & \frac{2}{m(m+2)}\cdot\sum_{i=1}^{j-1}\Big\{\mbox{tr}[(\bd{P}_i\bd{P}_j)^2]-\frac{1}{m}\cdot \big[\mbox{tr}(\bd{P}_i\bd{P}_j)\big]^2\Big\}.
\end{align}
Furthermore, by \eqref{hot_sleep2} and then Lemma \ref{dust_bone},
\begin{align}\lbl{273}
E[(\xi-E\xi)^4]
\leq & (K j)\cdot \sum_{i=1}^{j-1}E\Big[\bd{s}_i'\bd{M}_{ij}\bd{s}_i-
E(\bd{s}_i'\bd{M}_{ij}\bd{s}_i)\Big]^4 \nonumber\\
 \leq & \frac{K'j}{m^4}\sum_{i=1}^{j-1}
\Big\{\mbox{tr}(\bd{M}_{ij}^2)-\frac{1}{m}[\mbox{tr}(\bd{M}_{ij})]^2\Big\}^{2} \nonumber\\
= & \frac{K'j}{m^4}\sum_{i=1}^{j-1}
\Big\{\mbox{tr}[(\bd{P}_i\bd{P}_j)^2]-\frac{1}{m}[\mbox{tr}(\bd{P}_i\bd{P}_j)]^2\Big\}^{2}
\end{align}
as $m \geq 4\cdot 4+1=17$. Now we estimate the terms from \eqref{123}-\eqref{273}.

First, from \eqref{head_toss} and \eqref{trace_ralation1} we see that
\begin{eqnarray*}
T-2p\leq \mbox{tr}(\bd{P}_i\bd{P}_j)\leq T-p
\end{eqnarray*}
for any $1\leq i< j\leq N$, which implies that
\begin{eqnarray*}
\Big(1-\frac{p}{m}\Big)(j-1)\leq E\xi\leq (j-1)
\end{eqnarray*}
by \eqref{123}, the definition of $\bd{P}_{j\blacktriangle}$ and the notation $m=T-p$. Now, by (i) and (ii) from Lemma \ref{Love_sound},
\begin{align}\lbl{Victor_use}
\mbox{tr}[(\bd{P}_i\bd{P}_j)^2]-\frac{1}{m}[\mbox{tr}(\bd{P}_i\bd{P}_j)]^2
\leq  \mbox{tr}[(\bd{P}_i\bd{P}_j)^2]-\frac{1}{T}[\mbox{tr}(\bd{P}_i\bd{P}_j)]^2\leq C
\end{align}
for $1\leq i < j \leq N$. Hence
\begin{eqnarray*}
&& \frac{j}{2} \leq E\xi\leq j, \ \ \mbox{Var}(\xi)\leq \frac{Cj}{m^2}\ \  \mbox{and}\ \  E[(\xi-E\xi)^4]\leq \frac{Cj^2}{m^4}
\end{eqnarray*}
uniformly for all $2\leq j \leq N$ as $N$ is sufficiently large, where the ``$\frac{1}{2}$" appeared in the lower bound of $E\xi$ is not essential, it can be any positive number less than one. We then have from (ii) of Lemma \ref{mac_pencil} (taking $\alpha=2$) that
\begin{eqnarray*}
{\rm Var}(R_j)={\rm Var}(\xi^2) \leq C\cdot \frac{j^3}{m^2}
\end{eqnarray*}
uniformly for all $2\leq j \leq N$ as $N$ is sufficiently large.  This implies that
\begin{eqnarray*}
\sum_{j=2}^N{\rm Var}(R_j)=O\Big(\frac{N^4}{T^2}\Big)
\end{eqnarray*}
as $N\to\infty.$ 

\smallskip

\noindent\textbf{Proof of Lemma \ref{Lose_Ren}}. Similar to the last inequality from \eqref{957}, we have
\begin{align}\lbl{xiyouji}
&{\rm Var}\Big\{\sum_{j=2}^N{\rm tr}\Big[\Big(\sum_{i=1}^{j-1}\bd{U}_j'\bd{U}_i\bd{s}_i
\bd{s}_i'\bd{U}_i'\bd{U}_j\Big)^2\Big]\Big\} \nonumber\\
 \leq & (N-1)\cdot \sum_{j=2}^N{\rm Var}\Big\{{\rm tr}\Big[\Big(\sum_{i=1}^{j-1}\bd{U}_j'\bd{U}_i\bd{s}_i
\bd{s}_i'\bd{U}_i'\bd{U}_j\Big)^2\Big]\Big\}.
\end{align}
Use the formula that $(a_1+\cdots +a_n)^2=\sum_{1\leq i, k \leq n}a_ia_k$ for any real numbers $a_i$'s to see
\begin{eqnarray*}
&& \mbox{tr}\Big[\Big(\sum_{i=1}^{j-1}\bd{U}_j'\bd{U}_i\bd{s}_i
\bd{s}_i'\bd{U}_i'\bd{U}_j\Big)^2\Big]\\
&=& \sum_{1\leq i, k\leq j-1}\mbox{tr}\big(\bd{U}_j'\bd{U}_i\bd{s}_i
\bd{s}_i'\bd{U}_i'\bd{U}_j  \bd{U}_j'\bd{U}_k\bd{s}_k
\bd{s}_k'\bd{U}_k'\bd{U}_j\big)\\
& = & \sum_{1\leq i, k\leq j-1}\big(\bd{s}_i'\bd{U}_i'\bd{U}_j  \bd{U}_j'\bd{U}_k\bd{s}_k\big)^2
\end{eqnarray*}
since
\begin{eqnarray*}
& &\mbox{tr}\big(\bd{U}_j'\bd{U}_i\bd{s}_i\bd{s}_i'\bd{U}_i'\bd{U}_j  \bd{U}_j'\bd{U}_k\bd{s}_k
\bd{s}_k'\bd{U}_k'\bd{U}_j\big)\\
&=&\mbox{tr}\big[(\bd{s}_i'\bd{U}_i'\bd{U}_j  \bd{U}_j'\bd{U}_k\bd{s}_k)(\bd{s}_k'\bd{U}_k'\bd{U}_j\bd{U}_j'\bd{U}_i\bd{s}_i)\big]\\
&= & (\bd{s}_i'\bd{U}_i'\bd{U}_j  \bd{U}_j'\bd{U}_k\bd{s}_k)^2
\end{eqnarray*}
by \eqref{Science}. Set
\begin{align}\lbl{qingchun}
\bd{J}_{ijk}:=\bd{U}_i'\bd{U}_j  \bd{U}_j'\bd{U}_k
\end{align}
for all $1\leq i, j, k \leq N$. Of course, $J_{iji}=\bd{M}_{ij}$ which appears in Lemma \ref{trivial12}. Furthermore, $J_{ijk}'=J_{kji}$ for all $i,j,k$. Then $\bd{s}_i'\bd{J}_{ijk}\bd{s}_k=(\bd{s}_i'\bd{J}_{ijk}\bd{s}_k)'
=\bd{s}_k'\bd{J}_{kji}\bd{s}_i$. Thus,
\begin{eqnarray*}
\sum_{1\leq i, k\leq j-1}\big(\bd{s}_i'\bd{U}_i'\bd{U}_j  \bd{U}_j'\bd{U}_k\bd{s}_k\big)^2=\sum_{i=1}^{j-1}\big(\bd{s}_i'\bd{M}_{ij}\bd{s}_i\big)^2
+ 2\sum_{1\leq i< k\leq j-1}\big(\bd{s}_i'\bd{J}_{ijk}\bd{s}_k\big)^2.
\end{eqnarray*}
Thus,
\begin{align}\lbl{wcpl}
& {\rm Var}\Big\{{\rm tr}\Big[\Big(\sum_{i=1}^{j-1}\bd{U}_j'\bd{U}_i\bd{s}_i
\bd{s}_i'\bd{U}_i'\bd{U}_j\Big)^2\Big]\Big\} \nonumber\\
 \leq & 2\cdot {\rm Var}\Big[\sum_{i=1}^{j-1}\big(\bd{s}_i'\bd{M}_{ij}\bd{s}_i\big)^2\Big] + 8\cdot {\rm Var}\Big[\sum_{1\leq i< k\leq j-1}\big(\bd{s}_i'\bd{J}_{ijk}\bd{s}_k\big)^2\Big]
\end{align}
by the formula ${\rm Var}(g_1+g_2)\leq 2{\rm Var}(g_1)+2{\rm Var}(g_2)$ for any random variables $g_1$ and $g_2$. Now we estimate the last two terms one by one.

First, by independence,
\begin{eqnarray*}
{\rm Var}\Big[\sum_{i=1}^{j-1}\big(\bd{s}_i'\bd{M}_{ij}\bd{s}_i\big)^2\Big]
=\sum_{i=1}^{j-1}{\rm Var}\big[\big(\bd{s}_i'\bd{M}_{ij}\bd{s}_i\big)^2\big].
\end{eqnarray*}
By Lemma \ref{shadow_1}, $\mbox{Var}((\bd{s}_i'\bd{M}_{ij}\bd{s}_i)^2) \leq \frac{C}{m^2}$.   Therefore,
\begin{align}\lbl{cwuri}
{\rm Var}\Big[\sum_{i=1}^{j-1}\big(\bd{s}_i'\bd{M}_{ij}\bd{s}_i\big)^2\Big]
\leq \frac{Cj}{m^2}
\end{align}
uniformly for all $2\leq j \leq N.$

Second,
\begin{align}\lbl{73094}
& {\rm Var}\Big[\sum_{1\leq i< k\leq j-1}\big(\bd{s}_i'\bd{J}_{ijk}\bd{s}_k\big)^2\Big] \nonumber\\
 = & \sum_{1\leq i< k\leq j-1}{\rm Var}\big[\big(\bd{s}_i'\bd{J}_{ijk}\bd{s}_k\big)^2\big] + \sum{\rm Cov}\big(\big(\bd{s}_i'\bd{J}_{ijk}\bd{s}_k\big)^2, \big(\bd{s}_r'\bd{J}_{rjs}\bd{s}_s\big)^2\big),
\end{align}
where the last sum runs over all pairs $\{i,k\}$ and $\{r,s\}$ in the set $\{(i, k);\, 1\leq i<k\leq j-1\}$ satisfying
$(i, k)\ne (r, s)$ and $\{i,k\}\cap\{r,s\}\ne \emptyset$. Our remaining tasks are to evaluate the terms in the above two sums.

By (ii) of Lemma \ref{R_V},
\begin{align}\lbl{ma_po22}
\mbox{Var}\big[(\bd{s}_i'\bd{J}_{ijk}\bd{s}_k)^2\big] \leq \frac{C}{m^{5/2}}\big[\mbox{tr}\big((\bd{J}_{ijk}\bd{J}_{ijk}')^4\big)\big]^{1/2}.
\end{align}
 Write
\begin{align}\lbl{Jen_wu22}
\bd{J}_{ijk}\bd{J}_{ijk}'=\bd{U}_i'\bd{U}_j  \bd{U}_j'\bd{U}_k\bd{U}_k'\bd{U}_j\bd{U}_j'\bd{U}_i.
\end{align}
Now we need a fact from linear algebra that
\begin{align}\lbl{glade22}
\mbox{tr}(\bd{U}'\bd{A}\bd{U})\leq \mbox{tr}(\bd{A})
\end{align}
for any $T\times T$ nonnegative-definite matrix $\bd{A}$ and any $T\times m$ matrix $\bd{U}$ satisfying $\bd{U}'\bd{U}=\bd{I}_m$. In fact, take $\bar{\bd{U}}$ to be an $T\times (T-m)$ matrix such that $(\bd{U}, \bar{\bd{U}})$ is orthogonal. Easily,
\begin{eqnarray*}
(\bd{U}, \bar{\bd{U}})'\bd{A}(\bd{U}, \bar{\bd{U}})
=\begin{pmatrix}
\bd{U}'\bd{A}\bd{U} & *\\
* & \bar{\bd{U}}'\bd{A}\bar{\bd{U}}
\end{pmatrix}
.
\end{eqnarray*}
This leads to $\mbox{tr}(\bd{A})=\mbox{tr}(\bd{U}'\bd{A}\bd{U})+ \mbox{tr}(\bar{\bd{U}}'\bd{A}\bar{\bd{U}})$, and then we get \eqref{glade22} since $\bar{\bd{U}}'\bd{A}\bar{\bd{U}}$ is non-negative definite matrix.

Looking at \eqref{Jen_wu22}, we have from  \eqref{glade22} that
\begin{align}\lbl{cheng}
\mbox{tr}(\bd{J}_{ijk}\bd{J}_{ijk}')\leq \mbox{tr}(\bd{U}_j  \bd{U}_j'\bd{U}_k\bd{U}_k'\bd{U}_j\bd{U}_j')\leq \cdots \leq \mbox{tr}(\bd{U}_k\bd{U}_k')=m
\end{align}
since $\bd{U}_k'\bd{U}_k=\bd{I}_m$ from \eqref{happy_moon}. By the same argument, $\mbox{tr}\big((\bd{J}_{ijk}\bd{J}_{ijk}')^4\big) \leq m$. This and \eqref{ma_po22} imply
\begin{eqnarray*}
\mbox{Var}\big[(\bd{s}_i'\bd{J}_{ijk}\bd{s}_k)^2\big] \leq \frac{C}{m^{2}}.
\end{eqnarray*}
Hence, by \eqref{73094},
\begin{align}\lbl{ucnbkk}
{\rm Var}\Big[\sum_{1\leq i< k\leq j-1}\big(\bd{s}_i'\bd{J}_{ijk}\bd{s}_k\big)^2\Big] \leq \frac{Cj^2}{m^{2}} + \sum{\rm Cov}\big(\big(\bd{s}_i'\bd{J}_{ijk}\bd{s}_k\big)^2, \big(\bd{s}_r'\bd{J}_{rjs}\bd{s}_s\big)^2\big),\ \ \ \ \ \ \ \ \ \
\end{align}
where the last sum runs over all pairs $\{i,k\}$ and $\{r,s\}$ as stated below \eqref{73094}. Since there are only three free indices among those two pairs, it is easy to see that the total number of those pairs is no more than $\binom{j-1}{2}\cdot (j-3)\cdot 2\leq j^3$. Review the fact aforementioned  that $\bd{s}_i'\bd{J}_{ijk}\bd{s}_k=\bd{s}_k'\bd{J}_{kji}\bd{s}_i$ for all $1\leq i,j,k\leq N.$
By (iii) of Lemma \ref{R_V},
\begin{eqnarray*}
&&\mbox{Cov}\big[(\tilde{\bd{s}}'\bd{A}\tilde{\bd{s}}_1)^2, (\tilde{\bd{s}}'\bd{B}\tilde{\bd{s}}_2)^2\big]\\
&&=\frac{2}{m^3(m+2)}\cdot \mbox{tr}(\bd{A}\bd{A}'\bd{B}\bd{B}')-\frac{2}{m^4(m+2)}\,
\mbox{tr}(\bd{A}\bd{A}')\cdot\mbox{tr}(\bd{B}\bd{B}')
\end{eqnarray*}
for any $m\times m$ matrices $\bd{A}$ and $\bd{B}$, where $\tilde{\bd{s}}$, $\tilde{\bd{s}}_1$ and $\tilde{\bd{s}}_2$ are i.i.d. random vectors uniformly distributed over the $m$-dimensional sphere $\mathbb{S}^{m-1}$. Since $\tilde{\bd{s}}'\bd{A}\tilde{\bd{s}}_1=\tilde{\bd{s}}_1'\bd{A}'\tilde{\bd{s}}$, the above also implies that
\begin{eqnarray*}
&&\mbox{Cov}\big[(\tilde{\bd{s}}_1'\bd{A}\tilde{\bd{s}})^2, (\tilde{\bd{s}}'\bd{B}\tilde{\bd{s}}_2)^2\big]\\
&=&\frac{2}{m^3(m+2)}\cdot \mbox{tr}(\bd{A}'\bd{A}\bd{B}\bd{B}')-\frac{2}{m^4(m+2)}\,
\mbox{tr}(\bd{A}\bd{A}')\cdot\mbox{tr}(\bd{B}\bd{B}').
\end{eqnarray*}
By the same argument, we get similar bounds for $\mbox{Cov}\big[(\tilde{\bd{s}}_1'\bd{A}\tilde{\bd{s}})^2, (\tilde{\bd{s}}_2'\bd{B}\tilde{\bd{s}})^2\big]$ and $\mbox{Cov}\big[(\tilde{\bd{s}}'\bd{A}\tilde{\bd{s}}_1)^2, (\tilde{\bd{s}}_2'\bd{B}\tilde{\bd{s}})^2\big]$.
Hence, the maximum of the absolute values of the four covariances
is dominated by
\begin{align}\label{hei_an22}
&\frac{2}{m^4}\big[\mbox{tr}(\bd{A}\bd{A}'\bd{B}\bd{B}')+
\mbox{tr}(\bd{A}'\bd{A}\bd{B}\bd{B}')+\mbox{tr}(\bd{A}'\bd{A}\bd{B}'\bd{B})
\nonumber\\
&+\mbox{tr}(\bd{A}\bd{A}'\bd{B}'\bd{B})+\,\frac{1}{m}
\mbox{tr}(\bd{A}\bd{A}')\cdot\mbox{tr}(\bd{B}\bd{B}')\big]
\end{align}
(the bound above is an easy choice and we may choose a different one).

%
So the maximum of $|{\rm Cov}\big(\big(\bd{s}_i'\bd{J}_{ijk}\bd{s}_k\big)^2, \big(\bd{s}_r'\bd{J}_{rjs}\bd{s}_s\big)^2\big)|$ from \eqref{ucnbkk} is bounded by the maxima of the quantity in \eqref{hei_an22} with $\bd{A}=\bd{J}_{ijk}$ and $\bd{B}=\bd{J}_{rjs}$.  Now, recalling \eqref{Science} and \eqref{qingchun}, by using the same procedure as those in \eqref{glade22} and \eqref{cheng}, we know that each trace from \eqref{hei_an22} is bounded by $m$. Therefore,
\begin{eqnarray*}
|{\rm Cov}\big(\big(\bd{s}_i'\bd{J}_{ijk}\bd{s}_k\big)^2, \big(\bd{s}_r'\bd{J}_{rjs}\bd{s}_s\big)^2\big)| \leq \frac{10}{m^3}.
\end{eqnarray*}
This joining \eqref{ucnbkk} says that
\begin{eqnarray*}
&& {\rm Var}\Big[\sum_{1\leq i< k\leq j-1}\big(\bd{s}_i'\bd{J}_{ijk}\bd{s}_k\big)^2\Big] \leq \frac{Cj^2}{m^{2}} + \frac{Cj^3}{m^{3}}
\end{eqnarray*}
uniformly for all $3\leq j\leq N.$ Combining \eqref{wcpl}, \eqref{cwuri} and the above, we see
\begin{eqnarray*}
 {\rm Var}\Big\{{\rm tr}\Big[\Big(\sum_{i=1}^{j-1}\bd{U}_j'\bd{U}_i\bd{s}_i
\bd{s}_i'\bd{U}_i'\bd{U}_j\Big)^2\Big]\Big\} &\leq & C\cdot \Big(\frac{j}{m^2} + \frac{j^2}{m^{2}} + \frac{j^3}{m^{3}}\Big)\\
& \leq & (2C)\cdot \Big(\frac{j^2}{m^{2}} + \frac{j^3}{m^{3}}\Big).
\end{eqnarray*}
Recall $T=m+p$ with $p$ being fixed. By \eqref{xiyouji}, we arrive at
\begin{eqnarray*}
&&{\rm Var}\Big\{\sum_{j=2}^N{\rm tr}\Big[\Big(\sum_{i=1}^{j-1}\bd{U}_j'\bd{U}_i\bd{s}_i
\bd{s}_i'\bd{U}_i'\bd{U}_j\Big)^2\Big]\Big\} =O\Big(\frac{N^4}{T^2}+\frac{N^5}{T^3}\Big)
\end{eqnarray*}
as $N\to \infty$. The proof is finished. \hfill$\square$

\subsection{Proofs of auxiliary results in Section \ref{Pre_linear}}\lbl{Pre_linear_app} Although the results stated in  Section \ref{Pre_linear} serve the understanding of sample correlation coefficients $\hat{\rho}_{ij}$, their proofs have their own merits.

\smallskip

\noindent\textbf{Proof of Lemma \ref{shock}}. First, by the Chebyshev inequality,
\begin{align}
P(|a_1\xi_1+\cdots+ a_m\xi_m|\geq x)\leq & \frac{1}{x^2}\cdot E(a_1\xi_1+\cdots+ a_m\xi_m)^2 \nonumber\\
 = & \frac{1}{x^2}\lbl{rong}
\end{align}
since the last expectation is equal to $a_1^2+\cdots + a_m^2=1$. Let $\{\bar{\xi}_i;\, 1\leq i \leq m\}$ be an independent copy of $\{\xi_i;\, 1\leq i \leq m\}$. Then, we see $P(|a_1\bar{\xi}_1+\cdots+ a_m\bar{\xi}_m|\geq x/2)\leq \frac{4}{x^2} \leq \frac{1}{2}$ for $x\geq 3$, and hence $P(|a_1\bar{\xi}_1+\cdots+ a_m\bar{\xi}_m|< x/2)\geq \frac{1}{2}.$ Consequently,
\begin{align}\lbl{so_much}
& P\big(|a_1\xi_1+\cdots+ a_m\xi_m|\geq x)\nonumber\\
 \leq & 2P(|a_1\xi_1+\cdots+ a_m\xi_m|\geq x,\, |a_1\bar{\xi}_1+\cdots+ a_m\bar{\xi}_m|< \frac{x}{2}\big) \nonumber\\
 \leq & 2P\big(|a_1\eta_1+\cdots+ a_m\eta_m|\geq \frac{x}{2}\big),
\end{align}
where $\eta_i=\xi_i-\bar{\xi}_i$ for $1\leq i \leq m.$ The advantage in doing so is that $\eta_i$'s are symmetric and i.i.d. random variables with mean $0$, variance $2$  and $E|\eta_1|^\tau<\infty.$ Set $S_m=a_1\eta_1+\cdots+ a_m\eta_m.$ By a different version of the Hoffmann-J{\o}gensen inequality (Lemma 2.2 from \cite{LRJW}), for any integer $j\geq 1$, there exist positive constants $C_j$ and $D_j$ such that
\begin{align}\lbl{ajfycb}
P\Big(|S_m|\geq \frac{x}{2}\Big)\leq C_j\cdot P\Big(\max_{1\leq i \leq m}|a_i\eta_i|\geq \frac{x}{2j}\Big) + D_j\cdot P\Big(|S_m|\geq \frac{x}{4j}\Big)^j
\end{align}
for any $x>0.$ Similar to \eqref{rong},
\begin{align}\lbl{kcsur}
P\Big(|S_m|\geq \frac{x}{4j}\Big)\leq \frac{32j^2}{x^2}.
\end{align}
Furthermore,
\begin{align}\lbl{30975}
 P\Big(\max_{1\leq i \leq m}|a_i\eta_i|\geq \frac{x}{2j}\Big)
 \leq & \sum_{i=1}^m P\Big(|a_i\eta_i|\geq \frac{x}{2j}\Big) \nonumber\\
 \leq & \frac{(2j)^\tau}{x^\tau}\cdot E|\eta_1|^\tau\cdot \sum_{i=1}^m |a_i|^\tau \nonumber\\
 \leq & \frac{(2j)^\tau E|\eta_1|^\tau}{x^\tau}
\end{align}
since $\sum_{i=1}^m |a_i|^\tau\leq 1$ as $\tau\geq 2$. Combing \eqref{ajfycb}-\eqref{30975}, we have
\begin{eqnarray*}
P\Big(|S_m|\geq \frac{x}{2}\Big)\leq \frac{C_j'}{x^\tau}+\frac{D_j'}{x^{2j}}
\end{eqnarray*}
for all $x>0$, where $C_j'$ and $D_j'$ are constants depending on $j$ and $\tau.$ Taking integer $j\geq \tau/2$, we have
\begin{eqnarray*}
P\Big(|S_m|\geq \frac{x}{2}\Big)\leq \frac{K}{x^\tau}
\end{eqnarray*}
for $x\geq 3$, where $K$ is constant depending on $\tau.$   The desired conclusion follows from \eqref{so_much}. \hfill$\Box$

\smallskip

\noindent\textbf{Proof of Lemma \ref{Root}}. By the Taylor expansion, $e^y=1+y+\frac{1}{2}y^2+\frac{1}{6}y^3e^{\rho}$ for any $y\in \mathbb{R}$, where $\rho$ is between $0$ and $y.$ It follows that
\begin{eqnarray*}
e^{\theta \xi_1}
&=& 1+\theta \xi_1+\frac{1}{2}\theta^2\xi_1^2+\frac{1}{6}\theta^3\xi_1^3e^{\rho}\\
& \leq & 1+\theta \xi_1+\frac{1}{2}\theta^2\xi_1^2+\frac{1}{6}|\theta|^3|\xi_1|^3e^{\omega|\xi_1|/2}
\end{eqnarray*}
for all $\theta \in [-\omega/2, \omega/2].$ Set $\lambda=\frac{1}{6}E(|\xi_1|^3e^{\omega|\xi_1|/2})$. Then $\lambda<\infty$ since $Ee^{\omega|\xi_1|}<\infty$. It follows that
\begin{eqnarray*}
Ee^{\theta \xi_1}\leq 1+\frac{1}{2}\theta^2+\lambda|\theta|^3\leq \exp\Big(\frac{1}{2}\theta^2+\lambda|\theta|^3\Big)
\end{eqnarray*}
for all $\theta \in [-\omega/2, \omega/2].$ Now, notice $|a_i|\leq 1$ for each $i$, by the Markov inequality and the above,
\begin{eqnarray*}
P(a_1\xi_1+\cdots+ a_m\xi_m\geq x)
&\leq & e^{-\tau x}Ee^{\tau(a_1\xi_1+\cdots+ a_m\xi_m)}\\
& = & e^{-\tau x}\prod_{i=1}^mEe^{a_i\tau\xi_i}\\
& \leq & e^{-\tau x}\prod_{i=1}^m\exp\Big(\frac{1}{2}a_i^2\tau^2+\lambda|a_i|^2|\tau|^3\Big)
\end{eqnarray*}
for any $x\geq 0$ and  $\tau \in [0, \omega/2].$ From the assumption that $a_1^2+\cdots + a_m^2=1$ we see
\begin{eqnarray*}
P(a_1\xi_1+\cdots+ a_m\xi_m\geq x) \leq e^{-\tau x}\cdot \exp\Big(\frac{1}{2}\tau^2+\lambda|\tau|^3\Big)
\end{eqnarray*}
for all $\tau \in [0, \omega/2].$ By taking $\tau=\omega/2$, we get
\begin{eqnarray*}
P(a_1\xi_1+\cdots+ a_m\xi_m\geq x) \leq e^{\omega^2+\lambda \omega^3}\cdot e^{-\omega x/2}
\end{eqnarray*}
for all $x\geq 0$. Obviously, the above also holds if ``$a_i$" is replaced by ``$-a_i$". By taking $K=2e^{\omega^2+\lambda \omega^3}+\frac{2}{\omega}$, we have that
\begin{eqnarray*}
P(a_1\xi_1+\cdots+ a_m\xi_m\geq x) \leq K\cdot e^{-x/K}.
\end{eqnarray*}
The proof is completed.  \hfill$\Box$

\smallskip

\noindent\textbf{Proof of Lemma \ref{sleepless}}. Since $\xi_1$ is a subgaussian random variable, there exists $\sigma>0$ such that $Ee^{t\xi}\leq e^{\sigma^2t^2/2}$ for all $t>0.$ Hence,
\begin{eqnarray*}
P(a_1\xi_1+\cdots+ a_m\xi_m\geq x)
&\leq & e^{-tx}E\exp(t(a_1\xi_1+\cdots+ a_m\xi_m))\\
& = & e^{-tx}\prod_{i=1}^mEe^{ta_i\xi_i}\\
& \leq & e^{-tx}\cdot e^{(a_1^2+\cdots +a_m^2)\sigma^2t^2/2}\\
& = & e^{-tx+\sigma^2t^2/2}
\end{eqnarray*}
for all $t>0.$ Take $t=\frac{x}{\sigma^2}$ to get
\begin{eqnarray*}
P(a_1\xi_1+\cdots+ a_m\xi_m\geq x)\leq e^{-x^2/(2\sigma^2)}.
\end{eqnarray*}
Similarly, $P((-a_1)\xi_1+\cdots+ (-a_m)\xi_m\geq x)\leq e^{-x^2/(2\sigma^2)}$. The results then follows by taking $K=1/(2\sigma^2)$. \hfill$\Box$

\subsection{Proofs of auxiliary results in Section \ref{Pre_Asym}}\lbl{Pre_Asym_last}
Review  $\mathbb{S}^{T-1}$ stands for the unit sphere in the $T$-dimensional Euclidean space.

\smallskip

\noindent\textbf{Proof of Lemma \ref{Lili}}. By Theorem 1.5.7(i) and the argument for (5) on p.147 of \cite{mui}, the density of  $\bd{s}_1'\bd{s}_2$ is given by
\begin{eqnarray*}
g(\rho)=\frac{1}{\sqrt{\pi}}
\frac{\Gamma(\frac{T}{2})}{\Gamma(\frac{T-1}{2})}(1-\rho^2)^{(T-3)/2},\ \ |\rho| <1.
\end{eqnarray*}
Hence
\begin{eqnarray*}
P(\bd{s}_1'\bd{s}_2\geq l_N)=\frac{1}{\sqrt{\pi}}
\frac{\Gamma(\frac{T}{2})}{\Gamma(\frac{T-1}{2})}\int_{l_N}^1(1-\rho^2)^{(T-3)/2}\,d\rho.
\end{eqnarray*}
Let $t=t_q\in (0, 1)$ for each $q\geq 1$ satisfying $qt_q^2\to \infty$ as $q\to\infty$. By Lemma 6.2 from \cite{Cai_J},
\begin{eqnarray*}
\int_{t}^1(1-\rho^2)^{q/2}\,d\rho=\frac{1}{qt}(1-t^2)^{(q+2)/2}(1+o(1))
\end{eqnarray*}
as $q\to\infty.$
Now, by taking $q=T-3$ we have $ql_N^2=(T-3)\cdot\frac{4\log N}{T}(1+o(1))=4(\log N)(1+o(1))\to \infty.$ By (33) from \cite{Cai_J},
\begin{eqnarray*}
\frac{\Gamma(\frac{T}{2})}{\Gamma(\frac{T-1}{2})}=\sqrt{\frac{T}{2}}(1+o(1))
\end{eqnarray*}
as $N\to\infty.$ Consequently,
\begin{eqnarray*}
P(\bd{s}_1'\bd{s}_2\geq l_N)&=&\frac{1}{\sqrt{\pi}}\cdot \sqrt{\frac{T}{2}}\frac{1}{(T-3)l_N}(1-l_N^2)^{(T-1)/2}(1+o(1))\\
& = & \frac{1}{2\sqrt{2\pi}}\cdot \frac{1}{\sqrt{\log N}}\cdot\exp\Big[\frac{T}{2}\log \big(1-l_N^2\big)\Big]\cdot (1+o(1))
\end{eqnarray*}
since
\begin{eqnarray*}
l_N^2=\frac{4\log N -\log\log N+y}{T}.
\end{eqnarray*}
By the Taylor expansion, $\log (1-x)=-x+O(x^2)$ as $x\to 0.$ Then
\begin{eqnarray*}
\frac{T}{2}\log \big(1-l_N^2\big)
&=&\frac{T}{2}\big[-l_N^2+O\big(l_N^4\big)\big]\\
&=&-2\log N +\frac{1}{2}(\log\log N)-\frac{1}{2}y +O\Big(\frac{\log^2 N}{T}\Big)
\end{eqnarray*}
as $N\to\infty$. Then the conclusion follows from the assumption $\log N=o(\sqrt{T})$. \hfill$\square$

\smallskip

\noindent\textbf{Proof of Lemma \ref{Gaussian_inter}}. For any vector $\bd{a}\in \mathbb{S}^{T-1}$, the distribution of $\bd{a}'\bd{s}_2$ is independent of $\bd{a}$; see, e.g., Theorem 1.5.7(i) and the argument for (5) on p.147 from \cite{mui}. Hence, by taking $\bd{a}=(1, 0, \cdots, 0)'\in \mathbb{S}^{T-1}$ and using independence, we see $\bd{s}_1'\bd{s}_2$ has the same distribution as that of $Z_1(Z_1^2+\cdots + Z_T^2)^{-1/2}$, where $Z_1, \cdots, Z_T$ are i.i.d. $N(0, 1)$-distributed random variables.  Then
\begin{eqnarray*}
P\big(\max_{1\leq i \leq k}|\xi_i|\geq t\big)
&\leq & k\cdot P(|\xi_1|\geq t)\\
&= & k\cdot P\Big(\frac{|Z_1|}{\sqrt{Z_1^2+\cdots + Z_T^2}}\geq t\Big).
\end{eqnarray*}
By the large deviation bound for the sum of i.i.d. random variables (see, e.g., page 27 from \cite{DZ98}),
\begin{eqnarray*}
P\Big(\frac{1}{k}\sum_{i=1}^k Z_i^2\in A\Big) \leq 2\cdot \exp\big\{-k\inf_{x\in A}I(x)\big\}
\end{eqnarray*}
where $A\subset \mathbb{R}$ is any Borel set and $\Lambda(x)=\sup_{\theta\in \mathbb{R}}\{\theta x-\log E e^{\theta \xi^2} \}$, where $\xi$ is a $N(0,1)$ random variable. Since $\log Ee^{\theta \xi^2}=-\frac{1}{2}\log(1-2\theta)$ for $\theta<1/2,$
it is easy to check that
\begin{eqnarray*}
\Lambda(x)=\begin{cases}
\frac{1}{2}(x-1-\log x), & \text{if $x>0$;}\\
\infty, & \text{if $x\leq 0$.}
\end{cases}
\end{eqnarray*}
Observe that $\Lambda(x)$ is decreasing for $x\in (0, 1)$,  $\Lambda(1)=0$ and $\Lambda(0+)=\infty.$ Hence, for any $r\in (0, 1)$,
\begin{align}\lbl{times_let}
P(Z_1^2+\cdots + Z_T^2\leq rT)	\leq 2\cdot e^{-cT}
\end{align}
where $c=\Lambda(r)>0$. Thus,
\begin{eqnarray*}
&& P\big(\max_{1\leq i \leq k}|\xi_i|\geq t\big)\\
& \leq & k\cdot P\Big(\frac{|Z_1|}{\sqrt{Z_1^2+\cdots + Z_T^2}}\geq t, Z_1^2+\cdots + Z_T^2> rT\Big) +\\
& & k\cdot P(Z_1^2+\cdots + Z_T^2\leq rT)\\
& \leq & k\cdot P\big(|Z_1|\geq t\sqrt{rT}\big)+ (2k)\cdot e^{-cT}.
\end{eqnarray*}
Take $r=\frac{1}{2}$ and the result follows by the well-known inequality that 
$P(|Z_1|>x)\leq e^{-x^2/2}$ for $x \geq 1.$
\hfill$\Box$

\smallskip

\noindent\textbf{Proof of Lemma \ref{hua_dan}}. First,
\begin{eqnarray*}
\min_{1\leq i \leq k}v_i=\sqrt{\delta}Z + \sqrt{1-\delta}\cdot\min_{1\leq i \leq k}Z_i.
\end{eqnarray*}
If $\min_{1\leq i \leq k}v_i>x$ and $\sqrt{\delta}Z\leq y$, then
\begin{eqnarray*}
\min_{1\leq i \leq k}Z_i> \frac{x-y}{\sqrt{1-\delta}}.
\end{eqnarray*}
Then, for the event $\{\min_{1\leq i \leq k}v_i>x\}$, considering if $\sqrt{\delta}Z\leq y$ occurs or not, we have from independence that
\begin{eqnarray*}
P\big(\min_{1\leq i \leq k}v_i>x\big)
&\leq & P\big(\sqrt{\delta}Z>y\big) + P\Big(\min_{1\leq i \leq k}Z_i> \frac{x-y}{\sqrt{1-\delta}}\Big)\\
& \leq & P\Big(Z>\frac{y}{\sqrt{\delta}}\Big) + P\Big(Z_1> \frac{x-y}{\sqrt{1-\delta}}\Big)^k.
\end{eqnarray*}
Use the inequality that $P(Z_1>t)\leq \frac{1}{\sqrt{2\pi}\,t}e^{-t^2/2}$ for any $t>0$ to have
\begin{eqnarray*}
P\big(\min_{1\leq i \leq k}v_i>x\big)\leq \frac{1}{y}\exp\Big(-\frac{y^2}{2\delta}\Big)+ \frac{1}{(x-y)^k}\cdot \exp\Big[-\frac{k(x-y)^2}{2(1-\delta)}\Big].
\end{eqnarray*}
The proof is completed. \hfill$\Box$

\smallskip

\noindent\textbf{Proof of Lemma \ref{Zheng_zhou}}.  Let $Z_1, \cdots, Z_T$ be i.i.d. standard normals. Write  $\bd{Z}=(Z_1, \cdots, Z_T)'\in \mathbb{R}^T$. Then, $\bd{s}$ has the same distribution as that of $\frac{\bd{Z}}{\|\bd{Z}\|}$. Therefore, for each $r \in (0,1)$,
\begin{eqnarray*}
P\big(\min_{1\leq i \leq k}|\bd{a}_i'\bd{s}|>z\big) &= & P\big(\min_{1\leq i \leq k}|\bd{a}_i'\bd{Z}|>z\cdot \|\bd{Z}\|\big)\\
& \leq  & P\big(\min_{1\leq i \leq k}|\bd{a}_i'\bd{Z}|>z\cdot \|\bd{Z}\|,\, \|\bd{Z}\|>\sqrt{rT} \big) + P(\|\bd{Z}\|\leq \sqrt{rT} \big)\\
& \leq & P\big(\min_{1\leq i \leq k}|\bd{a}_i'\bd{Z}|>z\sqrt{rT} \big) + 2\cdot e^{-cT}
\end{eqnarray*}
where $c=c_r>0$ is a constant and the inequality in \eqref{times_let} is used in the last step. Observe that
\begin{eqnarray*}
\Big\{\min_{1\leq i \leq k}|\bd{a}_i'\bd{Z}|>z\sqrt{rT}\Big\}\subset \bigcup\Big\{\min_{1\leq i \leq k}\bd{\epsilon}_i\bd{a}_i'\bd{Z}>z\sqrt{rT}\Big\}
\end{eqnarray*}
where the union is taken over $2^k$ many events such that  $\bd{\epsilon}_i=\pm 1$ for each  $1\leq i \leq k$. Hence,
\begin{align}\lbl{where_pink}
P\Big(\min_{1\leq i \leq k}|\bd{a}_i'\bd{Z}|>z\sqrt{rT}\Big)\leq  \sum P\Big(\min_{1\leq i \leq k}\bd{\epsilon}_i\bd{a}_i'\bd{Z}>z\sqrt{rT}\Big)
\end{align}
where the sum runs over all possible $\epsilon_i=\pm 1$ for all $1\leq i \leq k.$
Easily,  the $k$-dimensional centered Gaussian random vector
\begin{eqnarray*}
\bd{u}:=\begin{pmatrix}
\epsilon_1\bd{a}_1'\\
 \vdots\\
 \epsilon_k\bd{a}_k'
\end{pmatrix}
_{k\times T}\cdot
\bd{Z}
\end{eqnarray*}
has covariance matrix
\begin{eqnarray*}
\bd{\Sigma}=E(\bd{u}\bd{u}')&=&\begin{pmatrix}
\epsilon_1\bd{a}_1'\\
 \vdots\\
 \epsilon_k\bd{a}_k'
\end{pmatrix} E(\bd{Z}\bd{Z}')\cdot (\epsilon_1\bd{a}_1, \cdots,
 \epsilon_k\bd{a}_k)\\
&=& (\bd{\epsilon}_i\bd{\epsilon}_j\bd{a}_i'\bd{a}_j)_{k\times k}.
\end{eqnarray*}
Obviously, the diagonal entries of $\bd{\Sigma}$ are all equal to $1$ because $\bd{a}_i$'s are unit vectors. By assumption, we have
\begin{eqnarray*}
\max_{1\leq i< j \leq k}(\bd{\epsilon}_i\bd{a}_i)'(\bd{\epsilon}_j\bd{a}_j)\leq \delta.
\end{eqnarray*}
By Lemma  \ref{Kahane_ineq}, we have that for all possible values of $\bd{\epsilon}_i$'s
\begin{eqnarray*}
P\Big(\min_{1\leq i \leq k}\bd{\epsilon}_i\bd{a}_i'\bd{Z}>z\sqrt{rT}\Big)\leq P\Big(\min_{1\leq i \leq k}v_i>z\sqrt{rT}\Big),
\end{eqnarray*}
where $(v_1, \cdots, v_k)'$ is a centered Gaussian random vector such that $E(v_i^2)=1$ for each $i$  and $E(v_iv_j)=\delta$ for all $i \ne j.$ Consequently, it is seen from \eqref{where_pink} that
\begin{eqnarray*}
P\Big(\min_{1\leq i \leq k}|\bd{a}_i'\bd{Z}|>z\sqrt{rT}\Big)\leq 2^k\cdot P\Big(\min_{1\leq i \leq k}v_i>z\sqrt{rT}\Big).
\end{eqnarray*}
 Without loss of generality, we are able to write
\begin{eqnarray*}
v_i=\sqrt{\delta}Z+\sqrt{1-\delta}Z_i
\end{eqnarray*}
for $1\leq i \leq k$, where $\{Z, Z_1, \cdots, Z_k\}$ are i.i.d. standard normals. We get the inequality by Lemma \ref{hua_dan}. \hfill$\Box$

\smallskip

\noindent\textbf{Proof of Lemma \ref{wen_fei}}. Review Lemma \ref{long_for}, we know  $\{\hat{\rho}_{ij};\, 1\leq i< j \leq N\}$ has the same distribution as that of
\begin{align}\lbl{waiting_5}
\{\bd{s}_i'\bd{U}_i'\bd{U}_j\bd{s}_j,\, 1\leq i < j \leq N \},
\end{align}
where $\bd{s}_1, \cdots, \bd{s}_N$ be i.i.d. random vectors uniformly distributed on $\mathbb{S}^{m-1}.$
Set $\bd{M}_{ij}=\bd{U}_i'\bd{U}_j\bd{U}_j'\bd{U}_i$ for any $1\leq i<j \leq N.$ By Lemma \ref{use_tea},
\begin{align}
& E[\hat{\rho}_{ij}|\bd{s}_i]=0\ \ \ \mbox{and}\ \ \   E[\hat{\rho}_{ij}^2|\bd{s}_i]=\frac{1}{m}\cdot \bd{s}_i'\bd{M}_{ij}\bd{s}_i;  \lbl{qianqian_jie}\\
& E\hat{\rho}_{ij}^2=\frac{1}{m^2}\cdot\mbox{tr}(\bd{P}_i\bd{P}_j). \lbl{jinggang}
\end{align}
Observe from \eqref{waiting_5} that  $\{\hat{\rho}_{ij};\, 1\leq j \leq N,\, j\ne i\}$ are conditionally  independent random variables given $\bd{s}_i.$ Denote by $P_1$,  $E_1$ and $\mbox{Var}_1$ the conditional probability, the conditional expectation and the conditional variance given $\bd{s}_i.$
Take $\xi=\hat{\rho}_{ij}$ in Lemma \ref{mac_pencil}(i). Then, by \eqref{qianqian_jie},
\begin{eqnarray*}
E_1\xi=0\ \ \ \mbox{and}\ \ \ \mbox{Var}_1(\xi)=\frac{1}{m}\cdot \bd{s}_i'\bd{M}_{ij}\bd{s}_i
\end{eqnarray*}
and
\begin{eqnarray*}
E_1(|\xi|^{2\tau})=E_1(|\bd{a}'\bd{s}_j|^{2\tau})\leq C\cdot \frac{\|\bd{a}\|^{2\tau}}{m^{\tau}}=C\cdot \frac{(\bd{s}_i'\bd{M}_{ij}\bd{s}_i)^{\tau}}{m^{\tau}}
\end{eqnarray*}
by Lemma \ref{gan_geming}, where $\bd{a}'=\bd{s}_i'\bd{U}_i'\bd{U}_j$. Lemma \ref{mac_pencil}(i) says that
\begin{eqnarray*}
E_1\big(\big| \hat{\rho}_{ij}^2-E_1\hat{\rho}_{ij}^2\big|^{\tau}\big)
\leq C\cdot \frac{(\bd{s}_i'\bd{M}_{ij}\bd{s}_i)^{\tau}}{m^{\tau}}.
\end{eqnarray*}
Therefore,
\begin{align}\lbl{974520}
E\big(\big| \hat{\rho}_{ij}^2-E_1\hat{\rho}_{ij}^2\big|^{\tau}\big)
\leq \frac{C}{m^{\tau}}\cdot E\big[(\bd{s}_i'\bd{M}_{ij}\bd{s}_i)^{\tau}\big].
\end{align}
Now we estimate the last expectation. By \eqref{sun_dark} and Lemma \ref{trivial},
\begin{eqnarray*}
\mbox{tr}(\bd{M}_{ij})=\mbox{tr}(\bd{P}_i\bd{P}_j)\leq m.
\end{eqnarray*}
By \eqref{mo_fact} and Lemma \ref{trivial} again,
\begin{eqnarray*}
\mbox{tr}(\bd{M}_{ij}^2)= \mbox{tr}((\bd{P}_i\bd{P}_j)^2)\leq m.
\end{eqnarray*}
It then follows from Lemma \ref{gan_geming} that
\begin{align}\lbl{4850heal}
E\big[(\bd{s}_i'\bd{M}_{ij}\bd{s}_i)^{\tau}\big]
\leq & \frac{C}{m^{\tau}}\cdot \Big\{[\mbox{tr}(\bd{M}_{ij})]^{\tau} + \Big[\mbox{tr}(\bd{M}_{ij}^2)-\frac{1}{m}[\mbox{tr}(\bd{M}_{ij})]^2\Big]^{\tau/2}\Big\}\nonumber\\
\leq & \frac{C}{m^{\tau}}\cdot \Big\{[\mbox{tr}(\bd{M}_{ij})]^{\tau} + \big[\mbox{tr}(\bd{M}_{ij}^2)\big]^{\tau/2}\Big\}\nonumber\\
\leq &  C
\end{align}
by Lemma \ref{brother_cat}. Then, \eqref{974520} and \eqref{4850heal} lead to the first conclusion. Now we prove the second one. Notice
\begin{align}\lbl{wei_tand}
& E\Big|\sum_{j=i+1}^N (\hat{\rho}_{ij}^2-E\hat{\rho}_{ij}^2)\Big|^{\tau} \nonumber\\
 \leq & 2^{\tau-1}\cdot E\Big|\sum_{j=i+1}^N \big[\hat{\rho}_{ij}^2-E(\hat{\rho}_{ij}^2|\bd{s}_i)\big]\Big|^{\tau} + 2^{\tau-1}\cdot E\Big|\sum_{j=i+1}^N[E(\hat{\rho}_{ij}^2|\bd{s}_i)-E\hat{\rho}_{ij}^2]\Big|^{\tau}.\ \ \ \ \
\end{align}
By \eqref{hot_sleep2} and the fact that $\{\hat{\rho}_{ij};\, 1\leq j \leq N,\, j\ne i\}$ are conditionally  independent random variables given $\bd{s}_i$, we see that
\begin{eqnarray*}
 E_1\Big|\sum_{j=i+1}^N \big[\hat{\rho}_{ij}^2-E(\hat{\rho}_{ij}^2|\bd{s}_i)\big]\Big|^{\tau}
 \leq   K_{\tau}\cdot (N-i)^{(\tau/2)-1}\cdot \sum_{j=i+1}^NE_1\big| \hat{\rho}_{ij}^2-E_1\hat{\rho}_{ij}^2\big|^{\tau}.
\end{eqnarray*}
Take expectation for both sides of the above and use the first conclusion to see that
\begin{align}\lbl{xcmo8rc}
E\Big|\sum_{j=i+1}^N \big[\hat{\rho}_{ij}^2-E(\hat{\rho}_{ij}^2|\bd{s}_i)\big]\Big|^{\tau}
 \leq   C\cdot \frac{(N-i)^{\tau/2}}{m^{\tau}}.
\end{align}
Now we estimate the last term from \eqref{wei_tand}. By \eqref{qianqian_jie} and \eqref{jinggang},
\begin{align}\lbl{leaksd}
 E\Big|\sum_{j=i+1}^N\big[E(\hat{\rho}_{ij}^2|\bd{s}_i)-E\hat{\rho}_{ij}^2\big]\Big|^{\tau}
 = &\frac{1}{m^{\tau}}\cdot E\Big|\sum_{j=i+1}^N \big[\bd{s}_i'\bd{M}_{ij}\bd{s}_i -E(\bd{s}_i'\bd{M}_{ij}\bd{s}_i)\big]\Big|^{\tau}\nonumber\\
 = & \frac{1}{m^{\tau}}\cdot E\big|\bd{s}_i'\bd{M}_{i\blacktriangledown}\bd{s}_i-E(\bd{s}_i'\bd{M}_{i\blacktriangledown}\bd{s}_i)\big|^{\tau},
\end{align}
where
\begin{eqnarray*}
\bd{M}_{i\blacktriangledown}:=\sum_{j=i+1}^N\bd{M}_{ij}.
\end{eqnarray*}
By Lemma \ref{dust_bone} again,
\begin{align}\lbl{yi_large}
 E\big|\bd{s}_i'\bd{M}_{i\blacktriangledown}\bd{s}_i-E(\bd{s}_i'
\bd{M}_{i\blacktriangledown}\bd{s}_i)\big|^{\tau} \leq  \frac{C}{m^{\tau}}\cdot \Big[\mbox{tr}(\bd{M}_{i\blacktriangledown}^2)-\frac{1}{m}[\mbox{tr}(\bd{M}_{i\blacktriangledown})]^2\Big]^{\tau/2}.
\end{align}
First,
\begin{eqnarray*}
\mbox{tr}(\bd{M}_{i\blacktriangledown})=\sum_{j=i+1}^N\mbox{tr}(\bd{P}_i\bd{P}_j).
\end{eqnarray*}
Easily, $\mbox{tr}(\bd{M}_{ij}\bd{M}_{ik})=\mbox{tr}(\bd{P}_{j}\bd{P}_i\bd{P}_k\bd{P}_i)$ since $\bd{M}_{ij}=\bd{U}_i'\bd{U}_j\bd{U}_j'\bd{U}_i$ and $\bd{U}_i\bd{U}_i'=\bd{P}_i$ for each $i=1\cdots N$ as stated in \eqref{happy_moon}. Hence,
\begin{eqnarray*}
\mbox{tr}(\bd{M}_{i\blacktriangledown}^2)=\sum_{i<j, k\leq N}\mbox{tr}(\bd{M}_{ij}\bd{M}_{ik})
=\sum_{i<j, k\leq N}\mbox{tr}(\bd{P}_{j}\bd{P}_i\bd{P}_k\bd{P}_i)=
\mbox{tr}\Big[\Big(\sum_{j=i+1}^N\bd{P}_{j}\bd{P}_i\Big)^2\Big].
\end{eqnarray*}
Define $\bd{P}_{i\bullet}=\sum_{j=i+1}^N\bd{P}_j.$
Then $\mbox{tr}(\bd{M}_{i\blacktriangledown}^2)=
\mbox{tr}(\big(\bd{P}_{i\bullet}\bd{P}_i\big)^2).$ It follows that
\begin{align}\lbl{eagle_peach}
\Big|\mbox{tr}(\bd{M}_{i\blacktriangledown}^2)-\frac{1}{m}[\mbox{tr}(\bd{M}_{i\blacktriangledown})]^2\Big|
=\Big|\mbox{tr}((\bd{P}_{i\bullet}\bd{P}_i)^2)-
\frac{1}{m}[\mbox{tr}(\bd{P}_{i\bullet}\bd{P}_i)]^2\Big|.
\end{align}
By taking $S=\{i+1, i+2, \cdots, N\}$, we get from (v) of Lemma \ref{Love_sound} that
\begin{eqnarray*}
\frac{1}{(N-i)^2}\cdot \Big|\mbox{tr}((\bd{P}_{i\bullet}\bd{P}_i)^2)-
\frac{1}{m}[\mbox{tr}(\bd{P}_{i\bullet}\bd{P}_i)]^2\Big|\leq C.
\end{eqnarray*}
This together with \eqref{yi_large} and \eqref{eagle_peach} concludes that
\begin{eqnarray*}
 E\big[\bd{s}_i'\bd{M}_{i\blacktriangledown}\bd{s}_i-E(\bd{s}_i'
\bd{M}_{i\blacktriangledown}\bd{s}_i)\big]^{\tau} \leq  \frac{C}{m^{\tau}}(N-i)^{\tau},
\end{eqnarray*}
which joins \eqref{leaksd} to yield
\begin{eqnarray*}
E\Big|\sum_{j=i+1}^N\big[E(\hat{\rho}_{ij}^2|\bd{s}_i)-
E\hat{\rho}_{ij}^2\big]\Big|^{\tau}
\leq  \frac{C}{m^{2\tau}}(N-i)^{\tau}.
\end{eqnarray*}
The second inequality then follows from the above, \eqref{wei_tand} and \eqref{xcmo8rc}. The third inequality is similarly obtained by simply replacing ``$\sum_{j=i+1}^N$" to ``$\sum_{i=1}^{j-1}$" in the above argument.  \hfill$\Box$

\end{document}